\documentclass[sn-mathphys]{sn-jnl}

\usepackage{mathrsfs}
\usepackage{mathtools}
\usepackage{amsmath}
\usepackage{amsfonts}
\usepackage{amsthm}
\usepackage{color}
\usepackage{xcolor}
\usepackage{graphicx}
\usepackage{tabularx}
\usepackage{makecell}

\usepackage{bm}

\usepackage{caption}
\usepackage{subcaption}

\usepackage{float}
\usepackage{url}
\usepackage{multicol}
\usepackage{txfonts}
\usepackage{setspace}
\usepackage{arydshln}
\usepackage{enumitem}
\usepackage{lipsum}
\usepackage{siunitx,booktabs}
\usepackage{nth}
\usepackage{accents} 

\labelformat{enumi}{#1}
\labelformat{enumii}{#1}
\labelformat{enumiii}{#1}
\labelformat{enumiv}{#1}

\jyear{2021}%

\theoremstyle{thmstyleone}%
\theoremstyle{thmstyletwo}%
\newtheorem{remark}{Remark}%

\theoremstyle{thmstylethree}%

\graphicspath{{figures/}}

\usepackage{listings}
\definecolor{light-gray}{gray}{0.95}
\begin{document}

\title[Design optimization of non-matching shells using free-form deformation]{Automated shape and thickness optimization for non-matching isogeometric shells using free-form deformation}


\author[1]{\mbox{\fnm{Han} \sur{Zhao}}}
\author[1]{\mbox{\fnm{David} \sur{Kamensky}}} 
\author[1]{\mbox{\fnm{John T.} \sur{Hwang}}} 
\author*[2,1]{\mbox{\fnm{Jiun-Shyan} \sur{Chen}}}\email{jsc137@ucsd.edu}

\affil[1]{\orgdiv{Department of Mechanical and Aerospace Engineering}, \orgname{University of California San Diego}, \orgaddress{\street{9500 Gilman Drive, Mail Code 0411}, \city{La Jolla}, \state{CA} \postcode{92093}, \country{USA}}}

\affil[2]{\orgdiv{Department of Structural Engineering}, \orgname{University of California San Diego}, \orgaddress{\street{9500 Gilman Drive, Mail Code 0085}, \city{La Jolla}, \state{CA} \postcode{92093}, \country{USA}}}

\abstract{Isogeometric analysis (IGA) has emerged as a promising approach in the field of structural optimization, benefiting from the seamless integration between the computer-aided design (CAD) geometry and the analysis model by employing non-uniform rational B-splines (NURBS) as basis functions. However, structural optimization for real-world CAD geometries consisting of multiple non-matching NURBS patches remains a challenging task. In this work, we propose a unified formulation for shape and thickness optimization of separately-parametrized shell structures by adopting the free-form deformation (FFD) technique, so that continuity with respect to design variables is preserved at patch intersections during optimization. Shell patches are modeled with isogeometric Kirchhoff--Love theory and coupled using a penalty-based method in the analysis. We use Lagrange extraction to link the control points associated with the B-spline FFD block and shell patches, and we perform IGA using the same extraction matrices by taking advantage of existing finite element assembly procedures in the FEniCS partial differential equation (PDE) solution library.  Moreover, we enable automated analytical derivative computation by leveraging advanced code generation in FEniCS, thereby facilitating efficient gradient-based optimization algorithms. The framework is validated using a collection of benchmark problems, demonstrating its applications to shape and thickness optimization of aircraft wings with complex shell layouts.}


\maketitle

\section*{Statements and declarations}
\subsection*{Conflict of interest}
The authors have no relevant financial or non-financial interests to disclose.

\section*{Article highlights}
\begin{itemize}
    \item Shape and thickness optimization for complex shell structures using free-form deformation.
    \item Lagrange extraction for automated isogeometric analysis and sensitivity computation of free-form deformation.
    \item Application to aircraft wings to streamline design through analysis, and optimization.
    \item Open-source implementation of proposed optimization approach using code generation in FEniCS.
\end{itemize}


\section{Introduction}\label{sec:intro}
Shell structures exhibit exceptional stiffness and strength-to-self-weight ratios, and are extensively employed in various engineering fields, such as aerospace, automotive, and marine engineering \cite{farshad2013design}. The performance of such structures is greatly influenced by geometric and material properties. Thus, structural optimization plays a vital role in obtaining superior designs for shell structures. In this paper, we present an optimization approach based on free-form deformation (FFD) \cite{Sederberg1986} to achieve the optimal shape and thickness distribution for isogeometric shell structures. 

Structural analysis is involved in the optimization process to evaluate the response of the current design and guide the subsequent iterations. The finite element (FE) method \cite{Johnson2012} is a well-established approach used to approximate PDEs with Lagrange polynomial basis functions. However, the discretization of the computational domain through interconnected simple elements, known as meshing, for complicated geometries is the primary challenge in FE analysis. The FE mesh generation and related process can account for up to 80\% of total analysis time \cite{Hardwick2005}. Alternatively, isogeometric analysis (IGA) \cite{Hughes05a, CoHuBa09} offers the possibility to bypass FE mesh generation by approximating the solution using the smooth non-uniform rational B-spline (NURBS) \cite{piegl1996nurbs} basis functions. NURBS is the industrial standard widely used to represent computer-aided design (CAD) models, making IGA an ideal method to streamline the design-through-analysis process.

IGA has gained growing interest since its introduction not only due to the unified description between design and analysis models but also the regularity provided by NURBS basis functions. The smoothness in splines allows for direct discretization of the Kirchhoff--Love shell model \cite{reddy2006theory}, a fourth-order PDE that requires at least $C^1$ continuity with the Galerkin method. Various applications including wind turbines \cite{Bazilevs2011, Johnson2020, Herrema2019a}, bioprosthetic heart valves \cite{Kamensky2015, kamensky2018contact, Xu2017b}, and car hoods \cite{zareh2019kirchhoff} have employed the isogeometric Kirchhoff--Love shell formulations \cite{Kiendl2009, Kiendl2011, Kiendl2015, Casquero2017, Casquero2020} and demonstrated exceptional results. However, practical CAD geometries are often too complex to be represented by a single tensor product NURBS surface. To make the CAD models with multiple NURBS patches directly available for analysis, \cite{Bazilevs10c} introduces a fictitious strip to add bending stiffness at patch interfaces with conforming discretizations. Additionally, various methods, such as Nitsche's method \cite{Guo2018, Benzaken2021}, penalty method \cite{herrema2019penalty, leonetti2020robust}, and super-penalty method \cite{coradello2021coupling, coradello2021projected} have been applied to CAD models with non-matching NURBS surfaces, further expanding the applicability of IGA in dealing with complex geometries.

The seamless integration between CAD and analysis models in IGA makes it a natural choice for design optimization. The updated design in the optimization process can be precisely captured in the analysis, which in turn ensures accurate responses due to the exact geometry representation and excellent approximation capabilities of spline basis functions \cite{Evans2009}. Shape optimization using IGA has been investigated extensively \cite{wall2008isogeometric, cho2009isogeometric, ha2010numerical, qian2010full, li2011isogeometric, azegami2013shape}. Many applications, such as beam structures \cite{nagy2010isogeometric}, vibrating membranes \cite{manh2011isogeometric}, shell structures \cite{nagy2013isogeometric, KIENDL2014148}, and complex photonic crystals \cite{qian2011isogeometric} show superior design by employing IGA in optimization. Moreover, topology optimization \cite{seo2010shape, qian2013topology} also benefits from the same spline basis in design models and analysis. Shape optimization for shell structures with stiffeners has been explored in \cite{hirschler2019embedded, hao2023isogeometric} using the FFD concept, a B-spline solid is extruded from a ``master'' part, which is stiffened with several ``slave'' stiffeners, to modify the shape of the whole shell structure. 

In this work, we employ the open-source framework PENGoLINS \cite{Zhao2022} for automated IGA of non-matching shell structures using the penalty method in \cite{herrema2019penalty}. The shape of the shell structure is updated through a trivariate B-spline FFD block, which encompasses the entire shell structure, without differentiating the ``master'' and ``slave'' parts. The FFD block modifies the Lagrange control points of all shell patches concurrently to preserve the surface--surface intersections. Subsequently, we obtain the resulting NURBS surfaces of shells using the Lagrange extraction technique \cite{Schillinger2016}, which is also implemented in the IGA using FE subroutines. Moreover, this approach is also applicable to shell thickness optimization where the thickness distribution is continuous at patch intersections. By integrating these two design variables, simultaneous shape and thickness optimization for non-matching shells can be effectively achieved. This combined optimization approach enables the exploration of complex design spaces while preserving the geometric integrity of the non-matching shell structure. To demonstrate the capability of the proposed method, we apply it to the design optimization of aircraft wings, effectively navigating the unconventional design space.

The structure of this paper is outlined as follows, we introduce commonly used notations and terminologies in Section \ref{sec:notation-terminology} for reference. Section \ref{sec:shell-coupling} reviews the penalty-based formulation for coupling of non-matching isogeometric Kirchhoff--Love shells and the computational algorithm for automated IGA of non-matching shell structures. We present the FFD-based shape and thickness optimization approach and formulate the sensitivities using the Lagrange extraction technique in Section \ref{sec:ffd-based-optimization}. The optimization approach is validated using a suite of benchmark problems in Section \ref{sec:benchmark} and is applied to aircraft wing design optimization in Section \ref{sec:application}, where superior design solutions are demonstrated. Lastly, we conclude the effectiveness of the proposed method and discuss potential future directions in Section \ref{sec:conclusion}.

\section{Notation and terminology} \label{sec:notation-terminology}
In this section, we provide a summary of commonly used notions and terminologies for reference, as the formulations in the following sections can become complex due to the use of the extraction concept in IGA, interval quadrature meshes for coupling separate spline patches, and FFD B-spline blocks in optimization. 
\begin{itemize}
    \singlespacing
    

    \item $p_{\text{sh}}$ : order of spline surfaces for shell structures.
    
    \item $\mathcal{V}^\text{I,IGA}$ : IGA function space for shell patch $S^\text{I}$. Superscript $\text{I}$  indicates shell patch index. For single patch formulations, index $\text{I}$ is neglected.

    \item $n^\text{I,IGA}$ : number of degrees of freedom (DoFs) in $\mathcal{V}^\text{I,IGA}$.
    
    \item $\mathbf{N}^{\text{I,IGA}}$ : spline basis functions in $\mathcal{V}^{\text{I,IGA}}$.

    \item $\mathbf{u}^{\text{I,IGA}}$ : displacement in IGA DoFs for shell patch $S^{\text{I}}$.

    \item $\mathbf{P}^{\text{I,IGA}}$ : NURBS control points for shell patch $S^{\text{I}}$.

    \item $\mathcal{V}^{\text{I,FE}}$ : FE function space for shell patch $S^{\text{I}}$.

    \item $n^\text{I,FE}$: number of DoFs in $\mathcal{V}^{\text{I,FE}}$. 

    \item $\mathbf{N}^{\text{I,FE}}$ : Lagrange polynomial basis functions in $\mathcal{V}^{\text{I,FE}}$.

    \item $\mathbf{u}^{\text{I,FE}}$ : displacement in FE DoFs for shell patch $S^{\text{I}}$.

    \item $\mathbf{P}^{\text{I,FE}}$ : Lagrange control points for shell patch $S^{\text{I}}$.

    \item $\mathbf{M}^{\text{I}}$ : extraction matrix for $\text{I}$-th shell patch that represents spline basis functions with Lagrange basis functions.

    \item $\boldsymbol{\xi}^{\text{MI}}$ : parametric coordinates of intersection $\mathcal{L}$ with respect to shell patch $S^{\text{I}}$.

    \item $\mathcal{V}^{\gamma \text{M}}$ : function space of the interval quadrature mesh. Superscript $\gamma \in \{0,1\}$ indicates the derivative order of the functions interpolated from $\mathcal{V}^{\text{I,FE}}$  to $\mathcal{V}^{\gamma \text{M}}$.

    \item $\mathbf{u}^{\gamma \text{MI}}$ : interpolated displacement with $\gamma$-th derivative from $\mathcal{V}^{\text{I,FE}}$ to $\mathcal{V}^{\gamma \text{M}}$.

    \item $\mathbf{P}^{\gamma \text{MI}}$ : interpolated control point functions with $\gamma$-th derivative from $\mathcal{V}^{\text{I,FE}}$ to $\mathcal{V}^{\gamma \text{M}}$.

    \item $\mathbf{T}^{\gamma \text{MI}}$ : interpolation or transfer matrix that interpolates $\gamma$-th derivative of functions from $\mathcal{V}^{\text{I,FE}}$ to $\mathcal{V}^{\gamma \text{M}}$.

    \item $p_{\text{FFD}}$ : order of B-spline basis functions for the FFD block.
    
    \item $\mathbf{N}_{\text{FFD}}$ : basis functions of 3D B-spline FFD block.
    
    \item $\mathbf{P}_{\text{FFD}}$: B-spline control points of FFD block.
    \singlespacing
\end{itemize}

In Section \ref{subsec:penalty-shell-coupling}, we review the basic Kirchhoff--Love shell and penalty formulations, and the IGA and FE function spaces are not differentiated in the subsection.


\section{Coupling of non-matching isogeometric shells}\label{sec:shell-coupling}
In aerospace structural applications, many components such as aircraft wings, empennage, and fuselage, can be modeled using shell theory. The Kirchhoff--Love shell model \cite{Kiendl2009} is employed in this work with IGA discretization. Typically, aircraft CAD geometries are composed of a collection of NURBS patches, and we adopt a penalty-based coupling method \cite{herrema2019penalty} to perform analysis for shells composed of multiple patches with arbitrary intersections and assume St. Venant--Kirchhoff material model. We first review the basic kinematics of Kirchhoff--Love shell theory and the coupling formulation for a pair of intersecting shell patches. We then elucidate the computational procedures for the analysis of non-matching shell structures.


\subsection{Penalty-based shell coupling }\label{subsec:penalty-shell-coupling}
The Kirchhoff--Love shell model neglects transverse shear strains, with straight lines normal to the midsurface remaining straight and normal to the midsurface before and after the deformation, and shell thickness kept unchanged. Thus the displacement field of the 3D shell can be uniquely represented by the displacement of its midsurface. In the reference configuration, the midsurface of a single-patch shell can be represented by $\mathbf{X}(\boldsymbol{\xi})$ 
, where $\boldsymbol{\xi}=\{\xi_1, \xi_2\}$, are parametric coordinates of the midsurface. In the deformed configuration with displacement $\mathbf{u}(\boldsymbol{\xi})$, the shell midsurface is expressed as 
\begin{align}
    \mathbf{x}(\boldsymbol{\xi}) = \mathbf{X}(\boldsymbol{\xi}) + \mathbf{u}(\boldsymbol{\xi}) \text{ .} 
    \label{eq:deformed-shell-config}
\end{align}
Taking derivatives of midsurface with respect to parametric coordinates, we can obtain the covariant basis vectors in reference and deformed configurations on the tangent plane as
\begin{align}
    \mathbf{A}_\alpha  = \frac{\partial \mathbf{X}}{\partial \xi_\alpha} \quad \text{ and } \quad \mathbf{a}_\alpha = \frac{\partial \mathbf{x}}{\partial \xi_\alpha} \quad \text{ for } \alpha \in \{1,2\} \text{ ,} \label{eq:midsurface-basis-vectors}
\end{align}
and unit vectors normal to the midsurface read as
\begin{align}
    \mathbf{A}_3 = \frac{\mathbf{A}_1 \times \mathbf{A}_2}{\Vert \mathbf{A}_1 \times \mathbf{A}_2 \Vert} \quad \text{ and } \quad \mathbf{a}_3 = \frac{\mathbf{a}_1 \times \mathbf{a}_2}{\Vert \mathbf{a}_1 \times \mathbf{a}_2 \Vert} \text{ .} \label{eq:midsurface-normal-vectors}
\end{align}
Membrane strain $\boldsymbol{\varepsilon}$ and change of curvature $\boldsymbol{\kappa}$ can be derived from \eqref{eq:midsurface-basis-vectors} and \eqref{eq:midsurface-normal-vectors}. Applying the appropriate material model, we can obtain the corresponding normal forces $\mathbf{n}$ and bending moments $\mathbf{m}$. The internal virtual work of Kirchhoff--Love shell is given by
\begin{align}
    \delta W_{\text{int}} = \int_{\mathcal{S}} (\mathbf{n}:\delta \boldsymbol{\varepsilon} + \mathbf{m} : \delta \boldsymbol{\kappa} )\ d\mathcal{S} \text{ ,} \label{eq:KL-shell-int-virtual-work}
\end{align}
where $\mathcal{S}$ is the midsurface of the shell. External virtual work is defined as
\begin{align}
    \delta W_{\text{ext}} = \int_{\Omega} \rho \mathbf{B} \cdot \delta \mathbf{u}\ d\Omega + \int_{\Gamma} \mathbf{T} \cdot \delta \mathbf{u}\ d\Gamma \text{ ,} \label{eq:KL-shell-ext-virtual-work}
\end{align}
where $\rho$ is the density, $\mathbf{B}$ is the body force, and $\mathbf{T}$ is the traction. $\Omega$ is the 3D domain of the shell volume with a thickness of $t$ about the midsurface, and $\Gamma$ is the 2D surface boundary where $\mathbf{T}$ is applied.
The principle of virtual work states
\begin{align}
    \delta W = \delta W_{\text{int}} - \delta W_{\text{ext}} = 0
    \text{ .} \label{eq:total-virtual-work}
\end{align}
\eqref{eq:total-virtual-work} is a nonlinear equation and can be solved using the Newton--Raphson method,
\begin{align}
    \frac{\partial^2 W}{(\partial \mathbf{u})^2}\,\Delta \mathbf{u} = -\frac{\partial W}{\partial \mathbf{u}} \text{.} \label{eq:single-patch-Newton-iteration}
\end{align}
Here, we only go through the fundamental equations of the Kirchhoff--Love shell model, readers are referred to \cite[Section 3]{Kiendl2011} for detailed derivation.

Practical complex shell structures typically consist of multiple NURBS patches. For example, a pair of intersecting shells joined at a certain angle, which is prevalent when modeling the stiffeners of an aircraft wing, cannot be described by one NURBS surface. To make the CAD geometry directly available for IGA, separate NURBS patches need to be coupled together so that both displacement continuity and joint angle conservation on the surface--surface intersections are maintained during deformation. 

\begin{figure}[!htb]\centering
    \includegraphics[width=0.8\textwidth]{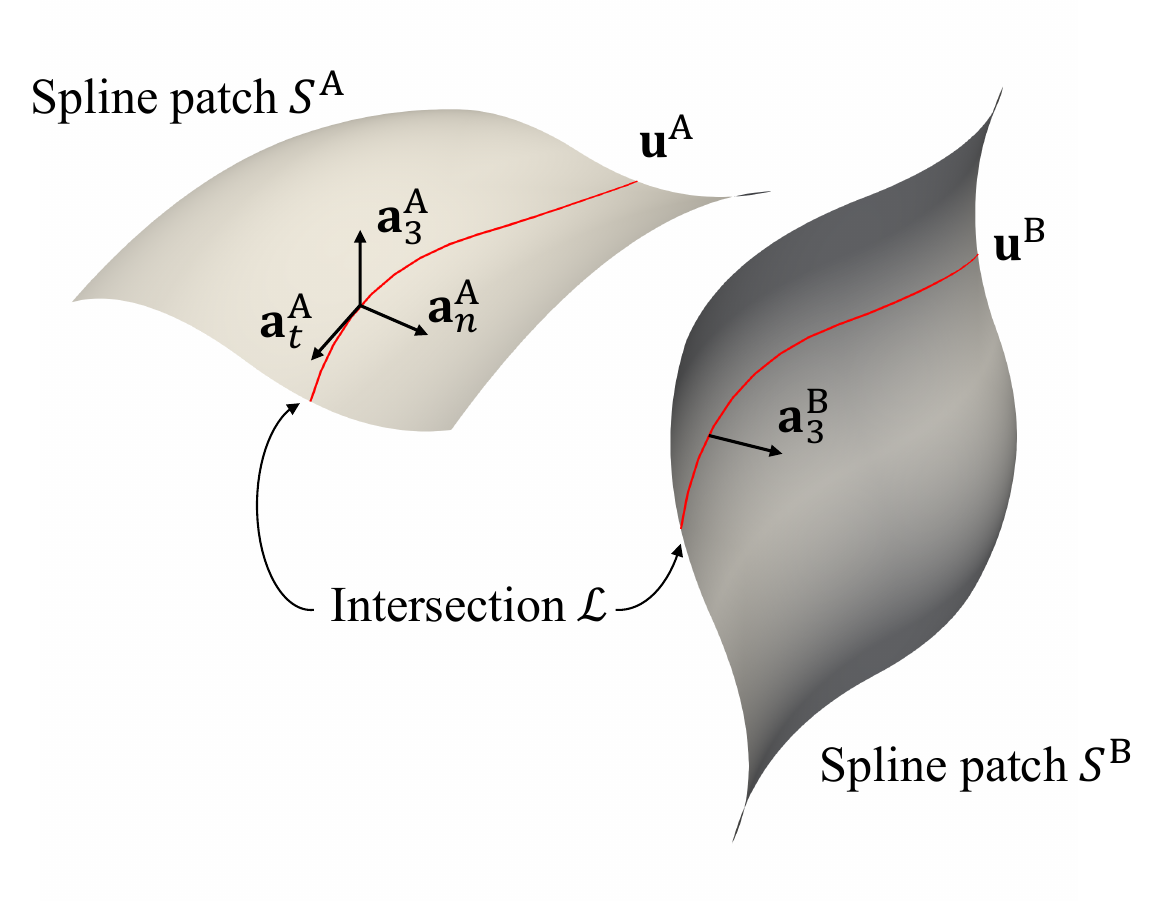}
    \caption{Spline patches $\mathcal{S}^\text{A}$ and $\mathcal{S}^\text{B}$ with one intersection $\mathcal{L}$ (indicated with red curves), where $\mathbf{u}$ and $\mathbf{a}_3$ are displacement and unit normal vector of midsurface, $\mathbf{a}_t$ and $\mathbf{a}_n$ are unit tangent vector and unit conormal vector on the intersection.}
    \label{fig:nonmatching-coupling-phy-sep}
\end{figure}

In this project, we use the penalty-based shell coupling formulation proposed by Herrema et al. \cite{herrema2019penalty}, which offers a good balance between accuracy and simplicity. Consider two shells modeled with spline patches $\mathcal{S}^\text{A}$ and $\mathcal{S}^\text{B}$ that share one intersection $\mathcal{L}$, as depicted in Figure \ref{fig:nonmatching-coupling-phy-sep}. The unit vector $\mathbf{a}^\text{A}_{t}$ is tangent to the intersection in the deformed configuration. In the reference and deformed configurations, conormal vectors are defined as  
\begin{align}
    \mathbf{A}^\text{A}_{n} = \mathbf{A}^\text{A}_{t} \times \mathbf{A}^\text{A}_{3} \quad \text{ and } \quad \mathbf{a}^\text{A}_{n} = \mathbf{a}^\text{A}_{t} \times \mathbf{a}^\text{A}_{3} \text{ .} 
    \label{eq:conormal-vector}
\end{align}
Subsequently, shell patches $\mathcal{S}^\text{A}$ and $\mathcal{S}^\text{B}$ are coupled through the following penalty energy
\begin{align}
    W^\text{AB}_\text{pen} = \frac{1}{2}\int_{\mathcal{L}}\alpha_\text{d}\left\Vert\mathbf{u}^\text{A}-\mathbf{u}^\text{B}\right\Vert^2 + \alpha_\text{r}\left(\left(\mathbf{a}^\text{A}_3\cdot\mathbf{a}_3^\text{B} - \mathbf{A}^\text{A}_3\cdot{\mathbf{A}}_3^\text{B}\right)^2 + \left(\mathbf{a}^\text{A}_n\cdot\mathbf{a}_3^\text{B} - \mathbf{A}^\text{A}_n\cdot{\mathbf{A}}_3^\text{B}\right)^2\right)\,d\mathcal{L}\text{ ,} \label{eq:penalty-energy}
\end{align}
where the $\alpha_\text{d}$ term penalizes the differences of the displacement on the intersection to maintain $C^0$ continuity, and the $\alpha_\text{r}$ term penalizes the change of the angle between $\mathcal{S}^\text{A}$ and $\mathcal{S}^\text{B}$ to keep joint angle conservation. Displacement and angle conservation penalty parameters $\alpha_\text{d}$ and $\alpha_\text{r}$ are defined based on a problem-independent, dimensionless penalty coefficient $\alpha$ scaled with element size and material properties. For isotropic materials, the penalty parameters are given by
\begin{align}
    \alpha_\text{d} = \frac{\alpha E t}{h (1-\nu^2)} \quad \text{and} \quad \alpha_\text{r} = \frac{\alpha E t^3}{12 h (1-\nu^2)} \text{ ,} \label{eq:penalty-parameters}
\end{align}
where $E$ is Young's modulus, $\nu$ is Poisson ratio, $t$ is shell thickness, and $h$ is element size. Details about the derivation of penalty parameters can be found in \cite[Section 2]{herrema2019penalty}.

\subsection{Computational procedures for non-matching shells} \label{subsec:computation-nonmatching-shell} 
In this section, we go through the computational procedures for structural analysis of isogeometric Kirchhoff--Love shells. The stiffness matrix of non-matching shells is formulated and used as the sensitivity in the design optimization in Section \ref{sec:ffd-based-optimization}.

\subsubsection{IGA of Kirchhoff--Love shells using extraction} \label{subsubsec:extraction-based-iga-kL-shell}

The concept of extraction \cite{BSEH11, SBVSH11, Schillinger2016, fromm2023interpolation} is utilized in the implementation of IGA, whose spline basis functions can be represented exactly by the linear combination of Lagrange basis functions. These Lagrange basis functions can be used in the classical FEM, allowing IGA to be implemented using finite element software with pre-defined extraction operators. An open-source IGA Python library named tIGAr is developed by Kamensky et al. \cite{Kamensky2019} using the finite element software FEniCS \cite{Logg2012}. Implementation and technical details are discussed in Section \ref{subsec:implementation}. In this section, we illustrate the basic mathematical operations and workflow of IGA using extraction. 


To perform IGA, an extraction matrix $\mathbf{M}$ is generated to represent functions defined in spline function space $\mathcal{V}^\text{IGA}$ with FE basis functions in $\mathcal{V}^\text{FE}$. The relation between these two sets of basis functions is given by 
\begin{align}
    \mathbf{N}^{\text{IGA}} = \mathbf{M}^T \mathbf{N}^{\text{FE}} \text{ ,} \label{eq:extraction-relation}
\end{align}
where $\mathbf{N}^{\text{IGA}}$ are IGA basis functions and $\mathbf{N}^{\text{FE}}$ are FE basis functions. Each column of $\mathbf{M}$ is the line combination of $\mathbf{N}^{\text{FE}}$ giving an IGA basis function.
In the analysis, we first create an extraction matrix $\mathbf{M}$ and assemble the stiffness matrix $\mathbf{K}^\text{FE}$ and force vector $\mathbf{F}^\text{FE}$ in $\mathcal{V}^\text{FE}$ using existing finite element subroutines. Then the displacement in $\mathcal{V}^\text{IGA}$ is solved as
\begin{align}
    (\mathbf{M}^T \mathbf{K}^\text{FE} \mathbf{M}) \mathbf{u}^\text{IGA} = \mathbf{M}^T \mathbf{F}^\text{FE} \text{ ,}
    \label{eq:single-patch-linear-system}
\end{align}
with problem-specific boundary conditions applied to $\mathbf{M}^T \mathbf{K}^\text{FE} \mathbf{M}$ and $\mathbf{M}^T \mathbf{F}^\text{FE}$.

\begin{remark}
    For the purpose of clarity, we assume control points of spline surfaces have unit weights. Therefore, rational spline basis functions are the same as homogeneous spline basis functions, both denoted as $\mathbf{N}^{\text{IGA}}$. In practice, weights need to be taken into consideration for correct geometric mapping and analysis.
    \label{rem:unit-weight}
\end{remark}

For single patch Kirchhoff--Love shell analysis, stiffness matrix $\mathbf{K}^\text{FE}$ is the second derivative of total work $\frac{\partial^2 W}{(\partial \mathbf{u}^\text{FE})^2}$. $\mathbf{M}^T \frac{\partial^2 W}{(\partial \mathbf{u}^\text{FE})^2} \mathbf{M}$ changes basis of $\frac{\partial^2 W}{(\partial \mathbf{u}^\text{FE})^2}$ from $\mathcal{V}^\text{FE}$ to  $\mathcal{V}^\text{IGA}$, and the formulation of IGA stiffness matrix can also be expressed as
\begin{align}
    \mathbf{M}^T \frac{\partial^2 W}{(\partial \mathbf{u}^\text{FE})^2} \mathbf{M} = \left(\frac{\partial \mathbf{u}^\text{FE}}{\partial \mathbf{u}^\text{IGA}}\right)^T \frac{\partial^2 W}{(\partial \mathbf{u}^\text{FE})^2} \frac{\partial \mathbf{u}^\text{FE}}{\partial \mathbf{u}^\text{IGA}} = \frac{\partial^2 W}{(\partial \mathbf{u}^\text{IGA})^2} \text{ .}
    \label{eq:stiffness-matrix-KL-shell}
\end{align}
The right-hand side (RHS) of \eqref{eq:single-patch-linear-system} for Kirchhoff--Love shell is equivalent to
\begin{align}
    \mathbf{M}^T \left(-\frac{\partial W}{\partial \mathbf{u}^\text{FE}} \right) = \left(\frac{\partial \mathbf{u}^\text{FE}}{\partial \mathbf{u}^\text{IGA}}\right)^T \left(-\frac{\partial W}{\partial \mathbf{u}^\text{FE}} \right) = -\frac{\partial W}{\partial \mathbf{u}^\text{IGA}} \text{ .}
    \label{eq:force-vector-KL-shell}
\end{align}
Therefore, \eqref{eq:single-patch-linear-system} is recovered as the linear system in $\mathcal{V}^\text{IGA}$ to solve for displacements increments in IGA DoFs,
\begin{align}
    \frac{\partial^2 W}{(\partial \mathbf{u}^\text{IGA})^2} \Delta \mathbf{u}^\text{IGA} = -\frac{\partial W}{\partial \mathbf{u}^\text{IGA}} \text{ .}
    \label{eq:single-patch-KL-shell-linear-system}
\end{align}

\subsubsection{Integration of penalty energy on patch intersections} \label{subsubsec:integration-penalty-energy}
\begin{figure}[!htb]\centering
    \includegraphics[width=0.9\textwidth]{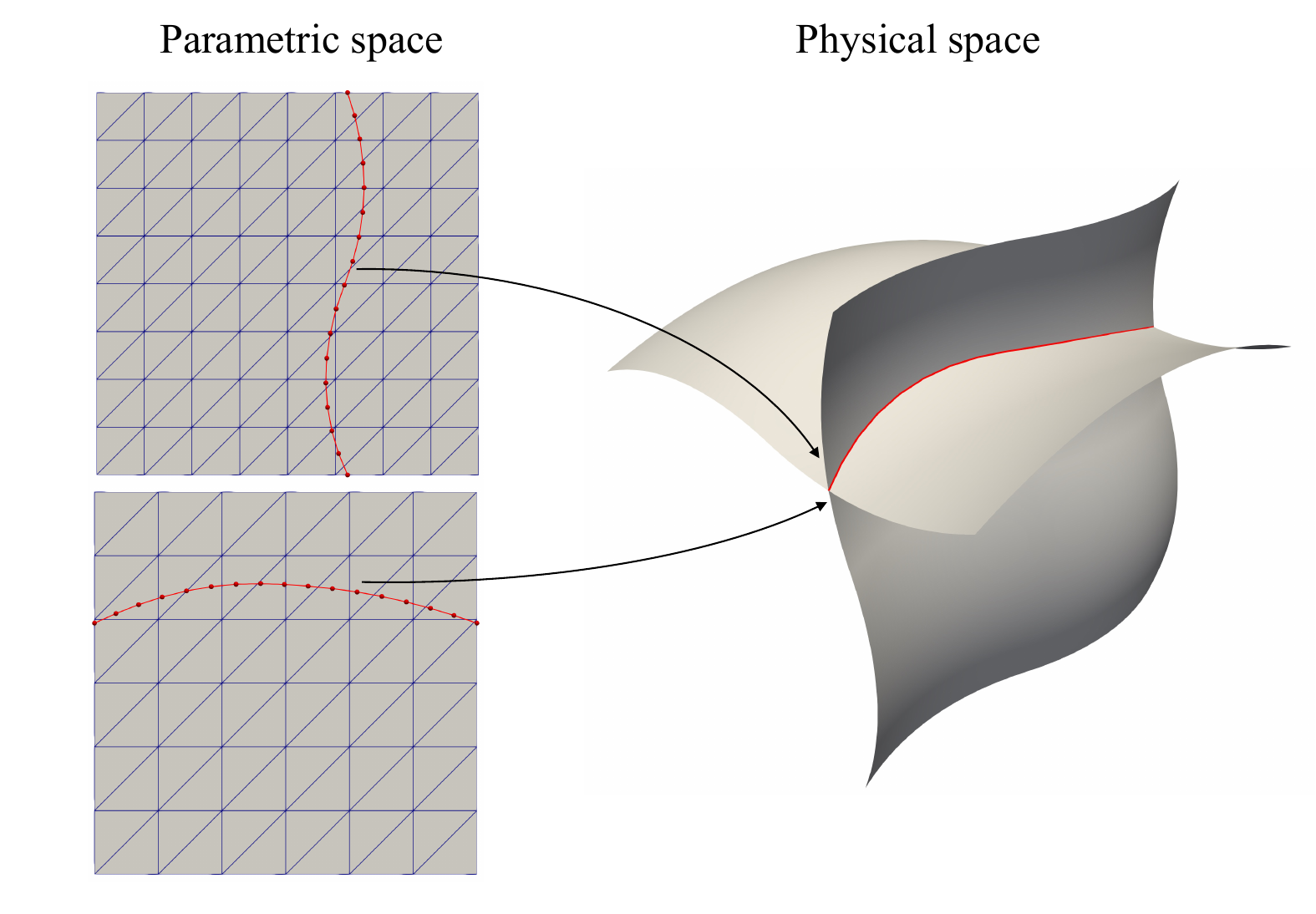}
    \caption{A pair of shell patches share one surface--surface intersection in the physical space, associated parametric surfaces are illustrated. A topologically 1D, geometrically 2D quadrature mesh (red interval mesh) is created in the parametric space and moved to parametric locations of the intersection with respect to two shells to create transfer matrices. Geometric and displacement quantities of shell patches are interpolated to the quadrature mesh, where the penalty energy is integrated.}
    \label{fig:nonmatching-coupling}
\end{figure}

For shell geometries that comprise multiple NURBS surfaces with non-matching intersections, the penalty energy is introduced to couple separate shell patches. In this section, we demonstrate the process to integrate the penalty energy \eqref{eq:penalty-energy} by using quadrature meshes. Again, taking two intersecting shells with one intersection as an illustrative example, the schematic configuration in parametric and physical spaces is shown in Figure \ref{fig:nonmatching-coupling}. We generate a geometrically 2D, topologically 1D interval mesh $\Omega^\text{M}$, which is named quadrature mesh, in the parametric space to represent the intersection curve for the penalty energy integration. The parametric coordinates of the surface--surface intersection with respect to two spline patches are denoted with $\boldsymbol{\xi}^\text{MA}$ and $\boldsymbol{\xi}^\text{MB}$. The overall algorithm to compute $W^\text{AB}_\text{pen}$ on $\Omega^\text{M}$ are outlined as follows:
\begin{enumerate}
\singlespacing
    \item Define function spaces $\mathcal{V}^\text{0M}$ and $\mathcal{V}^\text{1M}$ on $\Omega^\text{M}$ with linear FE basis functions to create displacements, control point functions, and their first derivatives.

    \item Move quadrature mesh $\Omega^\text{M}$ to parametric coordinate $\boldsymbol{\xi}^\text{MA}$ to generate transfer matrices between function spaces of shell patch $\mathcal{S}^\text{A}$ and quadrature mesh $\Omega^\text{M}$, where the configuration is shown in the upper left part of Figure \ref{fig:nonmatching-coupling}. The following transfer matrices are created:
    \begin{enumerate}
        \item $\mathbf{T}^\text{0MA} : \mathcal{V}^\text{A,FE} \rightarrow \mathcal{V}^\text{0M}$, which interpolates functions from $\mathcal{V}^\text{A,FE}$ to $\mathcal{V}^\text{0M}$. Each entry is defined as 
        \begin{align}
            \mathbf{T}_{ij}^\text{0MA} = \mathbf{N}_i^\text{A,FE}(\boldsymbol{\xi}^\text{MA}_j) \text{ ,} \label{eq:transfer-matrix-0-definition}
        \end{align}
        where $\mathbf{N}^{\text{A,FE}}_i$ is the $i$-th basis function in $\mathcal{V}^\text{A,FE}$, and $\boldsymbol{\xi}^\text{MA}_j$ is the $j$-th coordinate of $\Omega^\text{M}$.
        
        \item $\mathbf{T}^\text{1MA} : \mathcal{V}^\text{A,FE} \rightarrow \mathcal{V}^\text{1M}$, which interpolates the first derivative of functions from $\mathcal{V}^\text{A,FE}$ to $\mathcal{V}^\text{1M}$. $\mathbf{T}^\text{1MA}$ is defined as
        \begin{align}
            \mathbf{T}^\text{1MA}_{ij} = \mathbf{N}_{i,\,\boldsymbol{\xi}}^\text{A,FE}(\boldsymbol{\xi}^\text{MA}_j) \text{ ,} \label{eq:transfer-matrix-1-definition}
        \end{align}
        where $\mathbf{N}^{\text{A,FE}}_{i,\,\boldsymbol{\xi}}$ is the parametric derivative of $i$-th basis function in $\mathcal{V}^\text{A,FE}$.
    \end{enumerate}

    \item Compute displacement $\mathbf{u}^\text{0MA}$, geometric mapping $\mathbf{P}^\text{0MA}$, and their first derivatives $\mathbf{u}^\text{1MA}, \mathbf{P}^\text{1MA}$ from parametric domain of $\mathcal{S}^\text{A}$ to $\Omega^\text{M}$ using transfer matrices generated from last step,
    \begin{align}
        \mathbf{u}^{\gamma\text{MA}} = \mathbf{T}^{\gamma\text{MA}} \mathbf{u}^\text{A,FE} \quad \text{ and } \quad 
        \mathbf{P}^{\gamma\text{MA}} = \mathbf{T}^{\gamma\text{MA}} \mathbf{P}^\text{A,FE} \quad \text{ for } \gamma \in \{0,1\} \text{ .} 
        \label{eq:transfer-disp-geometric-mapping-A}
    \end{align}

    \item Substitute interpolated quantities from previous step into \eqref{eq:deformed-shell-config}--\eqref{eq:midsurface-normal-vectors} to obtain basis vectors on $\Omega^{\text{M}}$ before and after deformation: $\mathbf{A}_i^\text{MA}$ and $\mathbf{a}_i^\text{MA}$ for $i \in \{1,2,3\}$.

    \item Use \eqref{eq:conormal-vector} to compute conormal vectors: $\mathbf{A}_n^\text{MA}$ and $\mathbf{a}_n^\text{MA}$.
    
    \item Compute physical element size $h_X^{A}$ of $\mathcal{S}^\text{A}$, then interpolate it to $\Omega^\text{M}$ to obtain $h_X^\text{MA}$.

    \item Move quadrature mesh $\Omega^\text{M}$ to parametric coordinate $\boldsymbol{\xi}^\text{MB}$, which is depicted in the lower left configuration in Figure \ref{fig:nonmatching-coupling}. Then create transfer matrices $\mathbf{T}^\text{0MB}$ and $\mathbf{T}^\text{1MB}$.
    
    \item Repeat steps 3, 4 and 6 to obtain corresponding displacement $\mathbf{u}^\text{0MB}$, normal vectors $\mathbf{A}^\text{BM}_3$ and $\mathbf{a}^\text{BM}_3$, and element size $h^\text{MB}_X$ from $\mathcal{S}^\text{B}$.

    \item Calculate penalty parameters based on \eqref{eq:penalty-parameters} using averaged physical element size $h^\text{M}_X = \frac{1}{2}(h^\text{MA}_X + h^\text{MB}_X)$ for each node on $\Omega^\text{M}$,
    \begin{align}
        \alpha_\text{d}^\text{M} = \frac{\alpha E t^\text{M}}{h^\text{M}_X (1-\nu^2)} \quad \text{and} \quad \alpha_\text{r}^\text{M} = \frac{\alpha E (t^\text{M})^3}{12 h^\text{M}_X (1-\nu^2)} \text{ .}
    \end{align}
    \begin{remark}
        In this section, we assume $\mathcal{S}^\text{A}$ and $\mathcal{S}^\text{B}$ have the same thickness $t^\text{M} = t^\text{A} = t^\text{B} $ in penalty parameters, same for Young's modulus and Poisson ratio. For the thickness optimization problem discussed in Section \ref{subsec:sensitivities-thickness-opt}, shell patches may have different thicknesses, or each shell patch has variable thickness distribution. In these scenarios, we use an identical method to the element size calculation to interpolate thickness from shell patches to quadrature mesh and take the average, i.e., $t^{M} = \frac{1}{2}(t^\text{MA} + t^\text{MB})$.
        \label{rem:variable-thickness}
    \end{remark}

    \item Based on \eqref{eq:penalty-energy}, we can integrate the penalty energy on $\Omega^\text{M}$ with quantities computed from previous steps,
    \begin{align}
    \begin{split}
        W^\text{AB}_\text{pen} = \frac{1}{2}\int_{\Omega^\text{M}}  \biggl(&\alpha_\text{d}^\text{M} \Vert \mathbf{u}^\text{0MA} - \mathbf{u}^\text{0MB} \Vert^2 \\ 
        & + \alpha_\text{r}^\text{M} \Bigr((\mathbf{a}^\text{MA}_3 \cdot \mathbf{a}^\text{MB}_3 - \mathbf{A}^\text{MA}_3 \cdot \mathbf{A}^\text{MB}_3)^2 \\
        &\quad \quad \quad 
        + (\mathbf{a}^\text{MA}_n \cdot \mathbf{a}^\text{MB}_3 - \mathbf{A}^\text{MA}_n \cdot \mathbf{A}^\text{MB}_3)^2 \Bigr) \biggr)\,J^\text{M} \,d\Omega \text{ ,}
        \label{eq:penalty-energy-quadrature-mesh}
    \end{split}
    \end{align}
    where $J^\text{M}$ is the associated line Jacobian mapping $\Omega^\text{M}$ to $\mathcal{L}$.


\end{enumerate}

\subsubsection{Assembly of non-matching system} \label{subsubsec:assembly-nonmatching-system}
With the coupled formulations for intersecting shells, we can perform IGA on the non-matching shells directly. For a pair of shells with one intersection illustrated in Figure \ref{fig:nonmatching-coupling}, the total virtual energy in equilibrium is given by
\begin{align}
    \delta W_\text{t} = \delta W^\text{A} + \delta W^\text{B} + \delta W_\text{pen}^\text{AB} = 0 \text{ ,}
    \label{eq:total-virtual-energy}
\end{align}
which can be solved by the Newton--Raphson method. We can solve the linearized problem of \eqref{eq:total-virtual-energy} to obtain the full displacement $\mathbf{u}^\text{IGA} = \begin{bmatrix}
    \mathbf{u}^\text{A,IGA} & \mathbf{u}^\text{B,IGA}
\end{bmatrix}^T$ with one iteration,
\begin{align}
    \setlength\arraycolsep{6pt}
    \renewcommand\arraystretch{1.5}
    \begin{bmatrix}
        \frac{\partial^2 W_\text{t}}{(\partial \mathbf{u^\text{A,IGA}})^2} & \frac{\partial^2 W_\text{t}}{\partial \mathbf{u^\text{A,IGA}} \partial \mathbf{u^\text{B,IGA}}} \\
        \frac{\partial^2 W_\text{t}}{\partial \mathbf{u^\text{B,IGA}} \partial \mathbf{u^\text{A,IGA}}} & \frac{\partial^2 W_\text{t}}{(\partial \mathbf{u^\text{B,IGA}})^2}
    \end{bmatrix} \begin{bmatrix}
        \mathbf{u}^\text{A,IGA} \\
        \mathbf{u}^\text{B,IGA}
    \end{bmatrix} = -\begin{bmatrix}
        \frac{\partial W_\text{t}}{\partial \mathbf{u}^\text{A,IGA}} \\
        \frac{\partial W_\text{t}}{\partial \mathbf{u}^\text{B,IGA}}
    \end{bmatrix} \text{ .}
    \label{eq:nonmatching-linear-system}
\end{align}
Formulations for the derivative of $W_\text{t}$ can be obtained by means of the chain rule. Using $\mathcal{S}^\text{A}$ for demonstration, the derivative of $W_\text{t}$ with respect to $\mathbf{u}^\text{A,IGA}$ is given by
\begin{align}
    \frac{\partial W_\text{t}}{\partial \mathbf{u}^\text{A,IGA}} = (\mathbf{M}^{\text{A}})^T \left(\frac{\partial W^\text{A}}{\partial \mathbf{u}^\text{A,FE}} + \sum_{\gamma=0}^1 (\mathbf{T}^{\gamma\text{MA}})^T \frac{\partial W^\text{AB}_\text{pen}}{\partial \mathbf{u}^{\gamma\text{MA}}} \right) \text{ .}
    \label{eq:partial-wt-partial-uaiga}
\end{align}
The second derivatives of $W$ in the left-hand side (LHS) of \eqref{eq:nonmatching-linear-system} can be computed with the following formulations. Taking $\mathcal{S}^\text{A}$ for illustration, the diagonal block reads as
\begin{align}
\begin{split}
    \frac{\partial^2 W_\text{t}}{(\partial \mathbf{u^\text{A,IGA}})^2} = (\mathbf{M}^{\text{A}})^T \left( \frac{\partial^2 W^\text{A}}{(\partial \mathbf{u}^\text{A,FE})^2} + \sum_{\delta=0}^1 \sum_{\gamma=0}^1
     (\mathbf{T}^{\gamma\text{MA}})^T\frac{\partial^2 W^\text{AB}_\text{pen}}{\partial \mathbf{u}^{\gamma\text{MA}} \partial \mathbf{u}^{\delta\text{MA}}} \mathbf{T}^{\delta\text{MA}} \right) \, \mathbf{M}^\text{A} \text{ ,}
    \label{eq:ppartial-wt-ppartial-uaiga}
\end{split}
\end{align}
and the off-diagonal block can be expressed as
\begin{align}
    \frac{\partial^2 W_\text{t}}{\partial \mathbf{u^\text{A,IGA}} \partial \mathbf{u^\text{B,IGA}}} = (\mathbf{M}^{\text{A}})^T \left( \sum_{\delta=0}^1 \sum_{\gamma=0}^1 (\mathbf{T}^{\gamma\text{MA}})^T\frac{\partial^2 W^\text{AB}_\text{pen}}{\partial \mathbf{u}^{\gamma\text{MA}} \partial \mathbf{u}^{\delta\text{MB}}} \mathbf{T}^{\delta\text{MB}} \right) \, \mathbf{M}^\text{B} \text{ .}
    \label{eq:ppartial-wt-ppartial-uaiga-ubiga}
\end{align}
In \eqref{eq:partial-wt-partial-uaiga}--\eqref{eq:ppartial-wt-ppartial-uaiga-ubiga}, expressions of derivatives of $W^\text{A}$ can be found in \cite[Section 3]{Kiendl2009}, and \cite[Section 2.2]{herrema2019penalty} spells out derivatives of $W_\text{pen}^\text{AB}$. These derivatives are computed automatically by the implementation discussed in Section \ref{subsec:implementation} using computer algebra and code generation capabilities in FEniCS. Substituting \eqref{eq:partial-wt-partial-uaiga}--\eqref{eq:ppartial-wt-ppartial-uaiga-ubiga} into \eqref{eq:nonmatching-linear-system}, one can solve the displacements in IGA DoFs for $\mathcal{S}^\text{A}$ and $\mathcal{S}^\text{B}$.

\begin{remark}
    For clarity, we only demonstrate the non-matching system with two shell patches. However, \eqref{eq:nonmatching-linear-system} can be readily extended to shell geometries consisting of more than two surfaces. We refer interested readers to \cite[Section 3.2]{Zhao2022} for the assembly of shell structures with $m$ patches and arbitrary intersections. 
\end{remark}


\section{FFD-based design optimization}\label{sec:ffd-based-optimization}

CAD geometries of Kirchhoff--Love shells can be used for analysis without finite element mesh generation by employing formulations from Section \ref{sec:shell-coupling}. This provides attractive features for shape optimization of shell structures, where the discretization of shell patches stays unaltered during shape evolution. The updated geometry in optimization iterations stays consistent with the analysis model. As a result, this approach minimizes the effort required for geometry processing while simultaneously enhancing accuracy.

In the context of shape optimization for non-matching shell structures, it is crucial to ensure that updated shell patches remain properly connected. Failure to maintain this connectivity can result in separation or self-penetration of shell patches during the optimization iteration. Such issues would lead to false analysis and yield unrealistic optimal shapes. To tackle this challenge, we adopt the FFD-based \cite{Sederberg1986} technique combined with Lagrange extraction to update shell geometry and demonstrate the workflow in Section \ref{subsec:nonmatching-shell-update-ffd}. A comparable concept can be applied to thickness optimization to ensure continuous thickness distribution at the surface--surface intersections if needed. Sensitivities for both shape and thickness optimization are given in the subsequent sections.

\subsection{Non-matching shells update through FFD block} \label{subsec:nonmatching-shell-update-ffd}

\begin{figure}[!htb]\centering
    \includegraphics[width=1.0\textwidth]{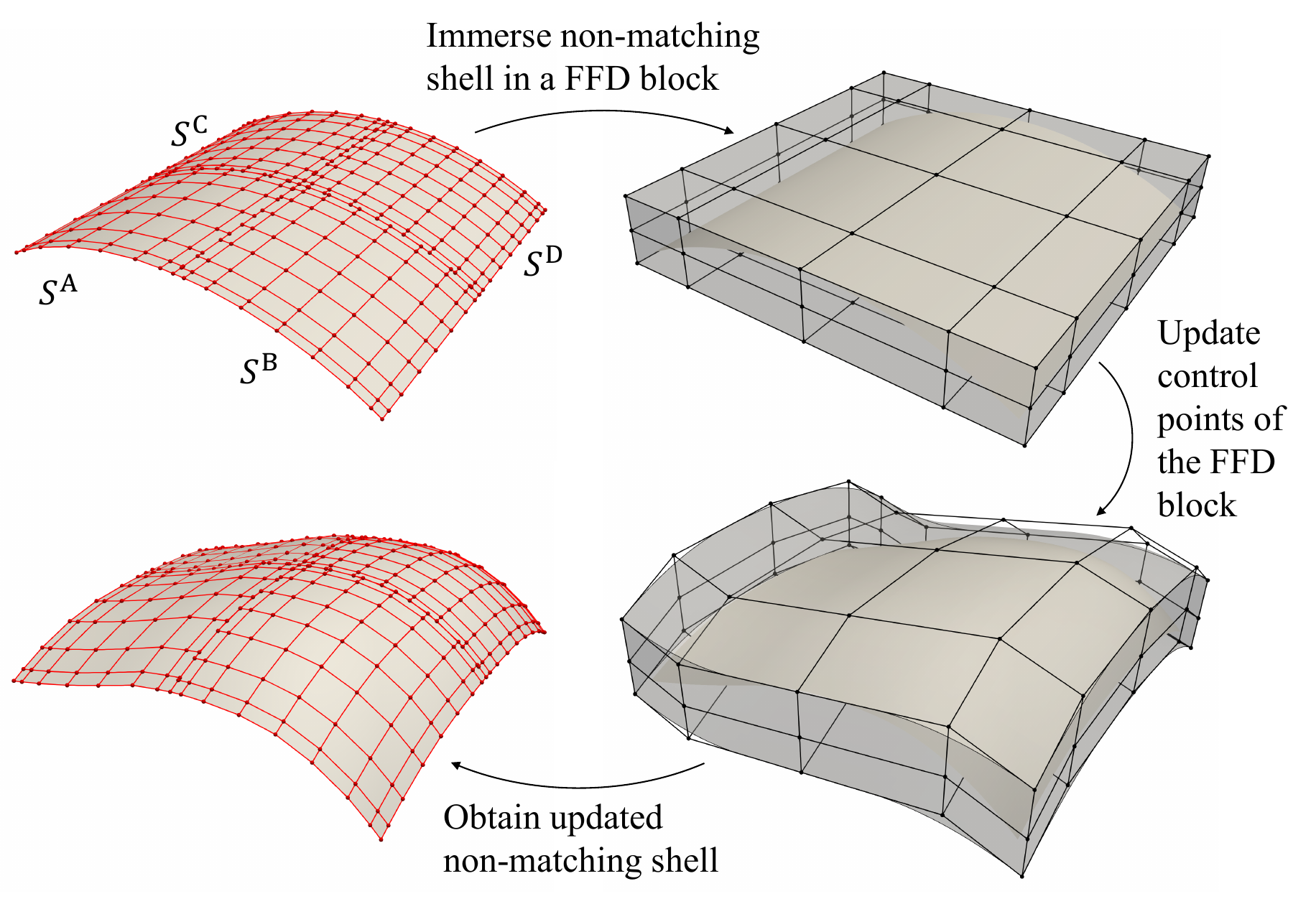}
    \caption{Workflow of FFD-based shape optimization for non-matching shell structures. A cylindrical roof consisting of four non-matching NURBS patches is first immersed in a trivariate B-spline block. We update the control points of the FFD block to deform the shape of the non-matching cylindrical roof. Control nets of the NURBS surfaces are indicated with red color, and black is used for the control net of the FFD block.}
    \label{fig:nonmatching-FFD-update}
\end{figure}

We use a cylindrical roof consisting of four non-matching shell patches that are shown in the upper-left part of Figure \ref{fig:nonmatching-FFD-update} as an example to demonstrate the FFD-based shape optimization approach. Note that this approach can be applied to shell structures with an arbitrary number of patches. 

\subsubsection{Approximating shell patches with extraction} \label{subsubsec:approx-shell-extraction}

For the initial CAD geometry consisting of $m$ Kirchhoff--Love shell patches, define a set of NURBS surface functions $\{\mathbf{S}^\text{A}(\boldsymbol{\xi}), \mathbf{S}^\text{B}(\boldsymbol{\xi}), \ldots, \mathbf{S}^{\text{I}^m}(\boldsymbol{\xi}) \}$, and the I-th shell patch $\mathbf{S}^\text{I}(\boldsymbol{\xi})$ is expressed as
\begin{align}
    \mathbf{S}^\text{I}(\boldsymbol{\xi}) = \mathbf{N}^{\text{I,IGA}}(\boldsymbol{\xi}) \mathbf{P}^{\text{I,IGA}} \text{ ,} \label{eq:surface-expression-IGA}
\end{align}
where $\mathbf{N}^{\text{I,IGA}}(\boldsymbol{\xi})$ are the spline basis functions of degree $p_\text{sh}$ in $\mathcal{V}^{\text{I,IGA}}$. We omit degree $p_\text{sh}$ in the notation for clarity. $\mathbf{P}^{\text{I,IGA}}$ are the NURBS control points for surface function $\mathbf{S}^\text{I}$.


Using the extraction concept discussed in Section \ref{subsubsec:extraction-based-iga-kL-shell}, NURBS surface function $\mathbf{S}^\text{I}(\boldsymbol{\xi})$ can be represented with Lagrange polynomials as well,
\begin{align}
    \mathbf{S}^\text{I}(\boldsymbol{\xi}) = \mathbf{N}^{\text{I,FE}}(\boldsymbol{\xi}) \mathbf{P}^{\text{I,FE}} \text{ ,} \label{eq:surface-expression-FE}
\end{align}
where $\mathbf{N}^{\text{I,FE}}(\boldsymbol{\xi})$ are basis functions in the finite element function space $\mathcal{V}^{\text{I,FE}}_\text{s}$ with nodal interpolatory property, and $\mathbf{P}^{\text{I,FE}}$ are Lagrange control points, or nodal values of $\mathbf{S}^\text{I}$. Plugging nodal coordinate $\boldsymbol{\xi}^{\text{I,FE}}$ of $\mathcal{V}^{\text{I,FE}}_\text{s}$ into \eqref{eq:surface-expression-FE}, coordinate of the NURBS surface $\mathbf{S}^\text{I}$ is represented with nodal values in the discrete setting,
\begin{align}
    \mathbf{S}^\text{I}(\boldsymbol{\xi}^{\text{I,FE}}) = \mathbf{N}^{\text{I,FE}}(\boldsymbol{\xi}^{\text{I,FE}}) \mathbf{P}^{\text{I,FE}} = \mathbf{P}^{\text{I,FE}} \text{ .} \label{eq:surface-approximate-FE-nodes}
\end{align}
Based on \eqref{eq:extraction-relation}, Lagrange control points can be obtained through the extraction matrix and NURBS control points. We have the following relation,
\begin{align}
    \mathbf{S}^\text{I}(\boldsymbol{\xi}^{\text{I,FE}}) = \mathbf{P}^{\text{I,FE}} = \mathbf{M}^\text{I} \mathbf{P}^{\text{I,IGA}} \text{ .}
    \label{eq:Lagrage-NURBS-control-points-relation}
\end{align}

\subsubsection{Relating FFD block control points and shell geometry } \label{subsubsec:relating-cp-ffd-shells}

The first step of Figure \ref{fig:nonmatching-FFD-update} illustrates the initial configuration of a collection of intersecting non-matching shell patches $S$, where red control nets are displayed. To enforce connectivity of the intersections during optimization, we immerse $S$ in a trivariate B-spline block, which is refered to as an FFD block, and use control points of the FFD block as design variables. A schematic demonstration is shown in the second step of Figure \ref{fig:nonmatching-FFD-update}. The FFD B-spline block is defined as
\begin{align}
    \mathbf{V}(\boldsymbol{\theta}) = \mathbf{N}_\text{FFD}(\boldsymbol{\theta}) \mathbf{P}_\text{FFD} \text{ ,} \label{eq:ffd-expression}
\end{align}
where $\boldsymbol{\theta}$ is the parametric coordinate of the FFD block, $\mathbf{N}_\text{FFD}(\boldsymbol{\theta})$ are B-spline solid basis functions of degree $p_\text{FFD}$ with knots vector, and $\mathbf{P}_\text{FFD}$ are B-spline block control points.

\begin{figure}[!htbp]
    \centering
    \begin{subfigure}[!htbp]{0.49\textwidth}
        \includegraphics[width=\textwidth]{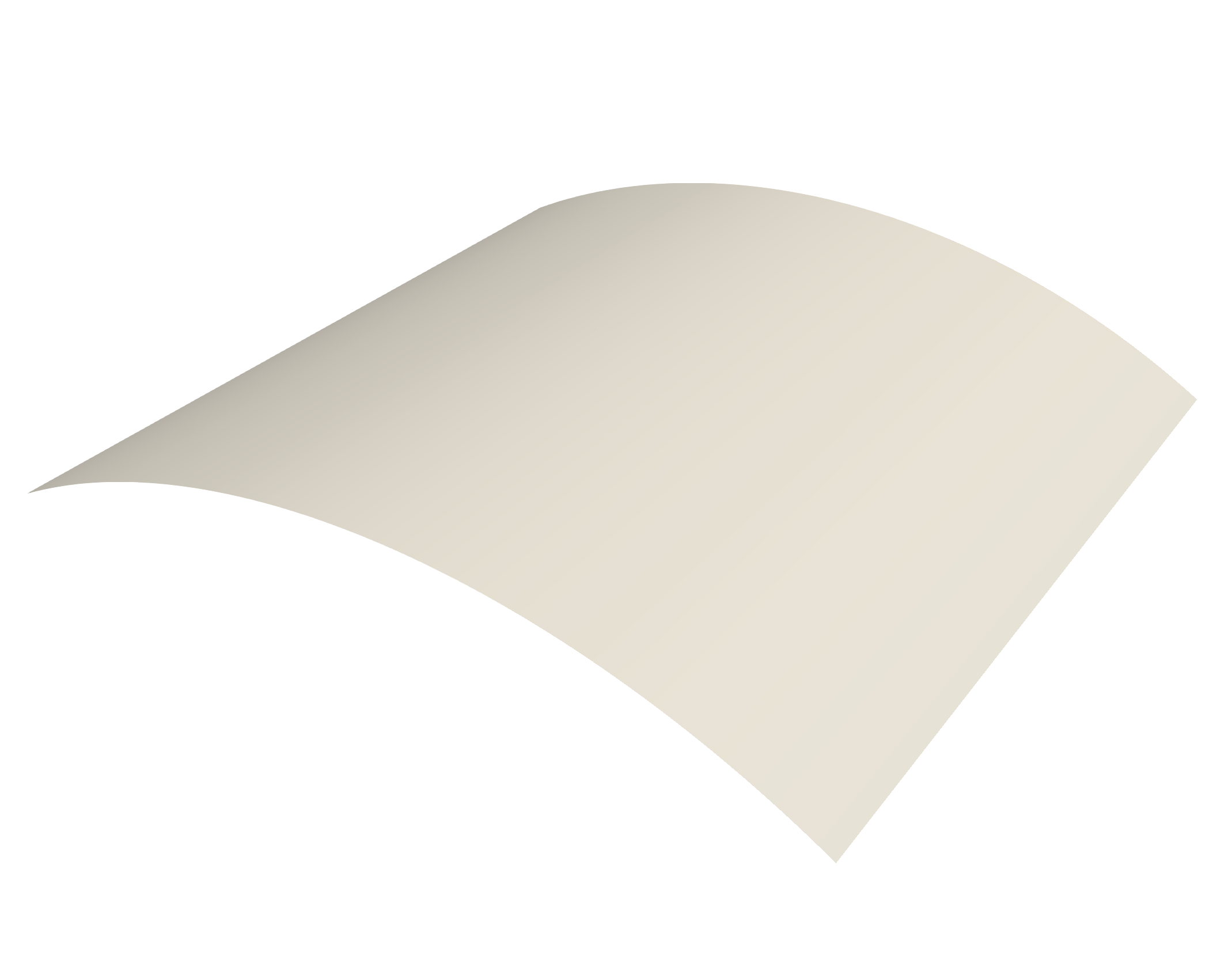}
        \caption{}
        \label{subfig:nonmatcing-shell-phy-init}
    \end{subfigure}
    \hfill
    \begin{subfigure}[!htbp]{0.49\textwidth}
        \includegraphics[width=\textwidth]{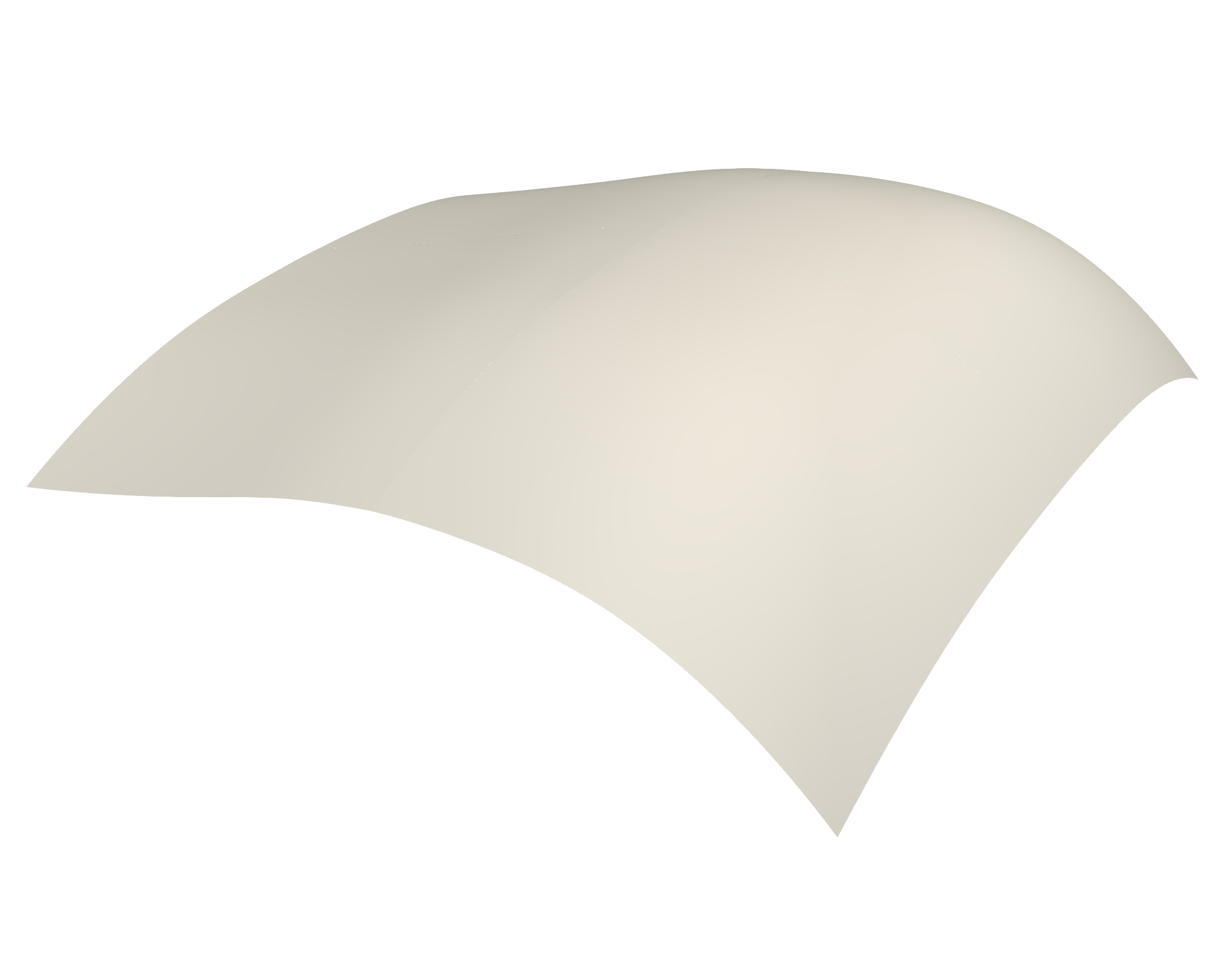}
        \caption{}
        \label{subfig:nonmatcing-shell-phy-new}
    \end{subfigure}
    \caption{(a) Initial configuration of the cylindrical roof geometry consisting of four non-matching NURBS patches. (b) Updated NURBS surfaces using FFD block.}
    \label{fig:nonmatching-shell-comparison}
\end{figure}

To simplify formulation and implementation, we use an identity mapping for the FFD block B-spline block, so that the parametric coordinate coincides with the physical coordinate,
\begin{align}
    \mathbf{V}(\boldsymbol{\theta}) = \mathbf{N}_\text{FFD}(\boldsymbol{\theta}) \mathbf{P}_\text{FFD} = \boldsymbol{\theta} \text{ .} \label{eq:ffd-block-identity-mapping}
\end{align}
Substituting \eqref{eq:Lagrage-NURBS-control-points-relation} into \eqref{eq:ffd-block-identity-mapping}, NURBS surfaces of the non-matching shells can be expressed using the FFD block basis functions and control points,
\begin{align}
    \mathbf{V}\left(\mathbf{S}^\text{I}(\boldsymbol{\xi})\right) = \mathbf{N}_\text{FFD}\left(\mathbf{S}^\text{I}(\boldsymbol{\xi})\right) \mathbf{P}_\text{FFD} = \mathbf{S}^\text{I}(\boldsymbol{\xi}) \text{ .} \label{eq:ffd-block-NURBS-surface}
\end{align}
In the continuous context of \eqref{eq:ffd-block-NURBS-surface},
shell patches will not separate in the final configuration as long as they are interconnected in the initial geometry. 
As the shape update of the FFD block is continuous, there is no relative movement between patches within the FFD block. 
In the discrete space, we can relate the NURBS control points of the shell patches to the control points of the FFD block,
\begin{align}
    \mathbf{N}_\text{FFD}\left(\mathbf{S}^\text{I}(\boldsymbol{\xi}^{\text{I,FE}})\right) \mathbf{P}_\text{FFD} = \mathbf{N}_\text{FFD}(\mathbf{P}^{\text{I,FE}}) \mathbf{P}_\text{FFD} = \mathbf{M}^\text{I} \mathbf{P}^{\text{I,IGA}} \text{ .} \label{eq:pffd-piga-relation}
\end{align}
The control points of the NURBS surface $\mathbf{P}^{\text{I,IGA}}$ of shell patches can be updated through the control points of the FFD block $\mathbf{P}_\text{FFD}$. Let $\mathbf{N}_\text{FFD}(\mathring{\mathbf{P}}^{\text{I,FE}})\coloneqq \mathbf{A}_{\text{FFD}}^{\text{I}}$, where $\mathring{\mathbf{P}}^{\text{I,FE}}$ denotes Lagrange control points of spline patch $\text{I}$ in the baseline configuration. Then $\mathbf{P}^{\text{I,IGA}}$ can be computed as
\begin{align}
    \mathbf{P}^{\text{I,IGA}} = \left( \left( (\mathbf{M}^{\text{I}})^T\mathbf{M}^\text{I} \right)^{-1}(\mathbf{M}^{\text{I}})^T \mathbf{A}_{\text{FFD}}^{\text{I}} \right) \mathbf{P}_{\text{FFD}} \text{ .} \label{eq:piga-pffd-equation}
\end{align}
It is noted that we need to solve the system using Moore--Penrose pseudo inverse due to the non-square nature of the extraction matrix $\mathbf{M}^\text{I}$, which has dimensions of $n^\text{I,FE}\times n^\text{I,IGA}$. For the extraction matrix, we have $n^\text{I,FE} > n^\text{I,IGA}$, which means that we are solving an overdetermined system. Therefore, $\mathbf{P}^{\text{I,IGA}}$ is considered as a least square fit in \eqref{eq:piga-pffd-equation} rather than an exact solution. The shape update strategy using FFD block is illustrated in the third step in Figure \ref{fig:nonmatching-FFD-update}, and the resulting shell patches with control net are depicted in the fourth step. A comparison between the initial non-matching cylindrical roof and updated NURBS surfaces is shown in Figure \ref{fig:nonmatching-shell-comparison}, where the surface--surface intersections keep overlapping within tolerance in the updated configuration.

The procedures to update control points of non-matching shells with $m$ patches are summarized as follows:
\begin{enumerate}
    \singlespacing
    \item In the preprocessing step, generate sparse matrices of evaluation of FFD block B-spline basis functions at shells' Lagrange control points in the initial configuration $\left\{\mathbf{A}_{\text{FFD}}^{\text{I}}\right\}$ and Lagrange extraction matrices $\left\{ \mathbf{M}^\text{I}\right\}\text{, for I} \in \{\text{A},\text{B},\ldots,\text{I}^m\}$.
    
    \item At optimization iteration step $i^\text{opt}$, obtain updated control points of the FFD block $\left( \mathbf{P}_{\text{FFD}} \right)^{i^\text{opt}}$. Compute updated Lagrange control points $\left(\mathbf{P}^{\text{I,FE}} \right)^{i^\text{opt}}$ for all shell patches,
    \begin{align}
        \mathbf{A}_{\text{FFD}}^{\text{I}} \left(\mathbf{P}_{\text{FFD}}\right)^{i^\text{opt}} = \left(\mathbf{P}^{\text{I,FE}} \right)^{i^\text{opt}} \text{ .} 
        \label{eq:update-Lagrange-control-points-step-i}
    \end{align}
    
    \item Solve NURBS control points $\left(\mathbf{P}^{\text{I,IGA}} \right)^{i^\text{opt}}$ at step $i^\text{opt}$ through Moore--Penrose pseudo inverse,
    \begin{align}
        \mathbf{M}^\text{I} \left(\mathbf{P}^{\text{I,IGA}} \right)^{i^\text{opt}} =  \left(\mathbf{P}^{\text{I,FE}} \right)^{i^\text{opt}} \text{ .} \label{eq:solve-NURBS-control-points-step-n}
    \end{align}

    \item Perform IGA with updated shell geometry, evaluate objective function and derivatives if needed, then proceed with optimization iteration.
    \singlespacing
\end{enumerate}

\begin{figure}[!htb]\centering
    \includegraphics[width=0.9\textwidth]{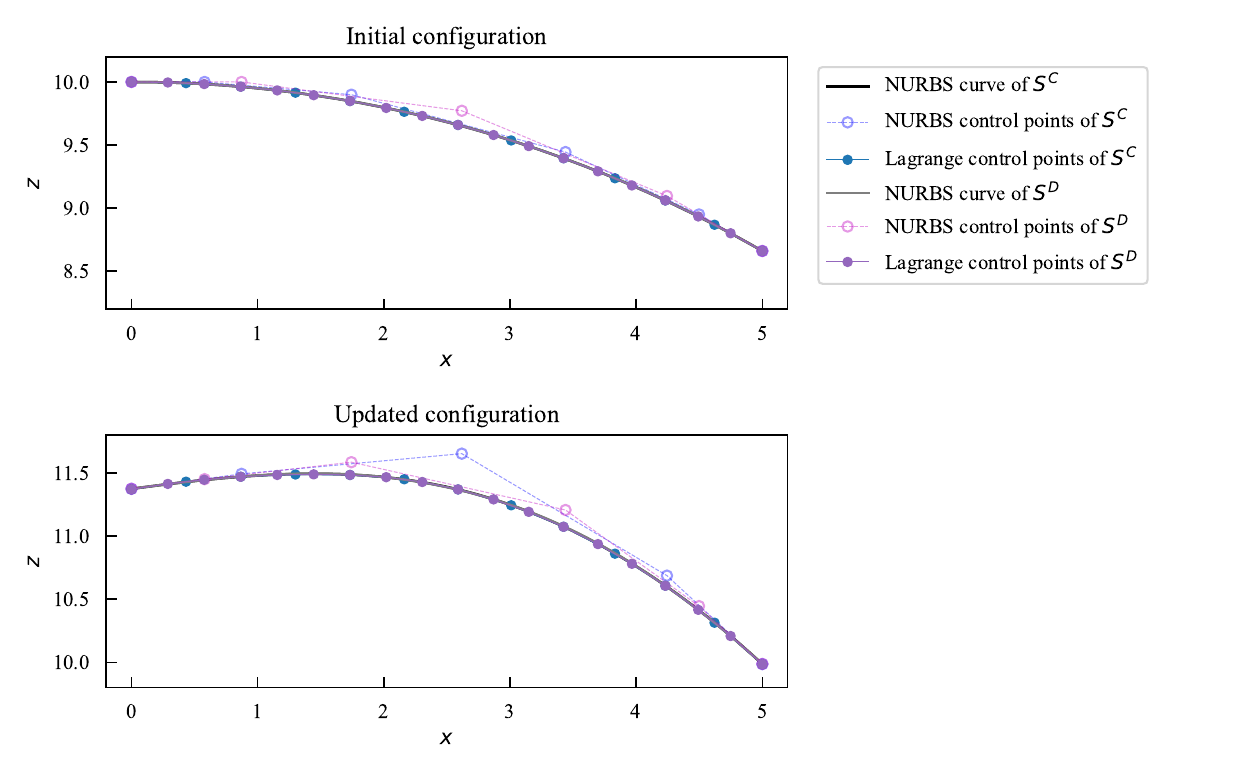}
    \caption{Sliced view of the intersecting edges between shell patches $\mathcal{S}^\text{C}$ and $\mathcal{S}^\text{D}$ of the cylindrical roof. The two edges remain overlapping in the updated configuration.}
    \label{fig:FFD-intersection-update}
\end{figure}

Though the control points of the shell patches are computed in the least square fit sense, the updated geometry can still retain the intersection with sufficient discretization. A sliced view of the intersection, in accompaniment with NURBS and Lagrange control points, between the right two spline patches $\mathcal{S}^\text{C}$ and $\mathcal{S}^\text{D}$ in the first step of Figure \ref{fig:nonmatching-FFD-update} is shown in Figure \ref{fig:FFD-intersection-update}, where we use coarser discretizations to make the comparison clearer. In the updated configuration, the two cubic intersecting edges are still overlapping even with only 5 and 6 NURBS control points.

\begin{remark}
    Since there is no relative movement between intersecting spline patches within the FFD block, which can be achieved with adequate control points in the discrete context, parametric coordinates of surface--surface intersections remain unchanged during shape updates. Therefore, transfer matrices introduced in Section \ref{subsubsec:integration-penalty-energy} can be reused to interpolate data from spline patches to quadrature meshes when integrating penalty energies in the optimization iteration. These matrices only need to be generated once at the preprocessing stage. \label{rem:transfer-matrices-preprocessing}
\end{remark}

\subsection{Sensitivities for shape optimization} \label{subsec:sensitivities-shape-opt}
By utilizing the capabilities of direct analysis for non-matching isogeometric shells and incorporating FFD-based shape updates, we are able to conduct shape optimization for the shell structures in a seamless manner. The problem that optimizes the shape of non-matching shells can be formulated as follows,
\begin{align}
\begin{split}
    \text{minimize} \quad &f^\text{obj}\left(\mathbf{P}_{\text{FFD}}\right) \\
    \text{subject to} \quad & g^\text{con}_{i_g}\left(\mathbf{P}_{\text{FFD}}\right) \leq \mathbf{0} \text{ , for }  i_g \in \{1,2,\ldots,n_g\} \\
    & h^\text{con}_{i_h}\left(\mathbf{P}_{\text{FFD}}\right) = \mathbf{0} \text{ , for } i_h \in \{1,2,\ldots,n_h\} \\
    & \mathbf{P}_{\text{FFD}l} \leq \mathbf{P}_{\text{FFD}} \leq \mathbf{P}_{\text{FFD}u} \text{ ,} \label{eq:shape-minimization} 
\end{split}
\end{align}
where control points of the FFD block $\mathbf{P}_{\text{FFD}}$ are design variables, $f^\text{obj}$ is the objective function, $g^\text{con}_{i_g}$ are inequality constraints, and $h^\text{con}_{i_h}$ are equality constraints. $\mathbf{P}_{\text{FFD}l}$ and $\mathbf{P}_{\text{FFD}u}$ are lower and upper limits for the design variables.

To perform gradient-based design optimization, we formulate total derivatives of the objective function with respect to design variables as
\begin{align}
    \frac{d f^\text{obj}}{d \mathbf{P}_{\text{FFD}}} = \left( \frac{\partial f^\text{obj}}{\partial \mathbf{P}_{}^{\text{IGA}}} + \frac{\partial f^\text{obj}}{\partial \mathbf{u}^\text{IGA}} \frac{d \mathbf{u}^\text{IGA}}{d \mathbf{P}_{}^{\text{IGA}}} \right) \frac{d \mathbf{P}_{}^{\text{IGA}}}{d \mathbf{P}_{\text{FFD}}} \text{ ,} \label{eq:shape-opt-total-derivative-init}
\end{align}
where $\mathbf{P}_{}^{\text{IGA}} = \begin{bmatrix}
    \mathbf{P}_{}^{\text{A,IGA}} & \mathbf{P}_{}^{\text{B,IGA}} & \ldots & \mathbf{P}_{}^{\text{I}^m\text{,IGA}}
\end{bmatrix}^T$ is the full vector of NURBS control points, and similarly, $\mathbf{u}^{\text{IGA}} = \begin{bmatrix}
    \mathbf{u}^{\text{A,IGA}} & \mathbf{u}^{\text{B,IGA}} & \ldots & \mathbf{u}^{\text{I}^m\text{,IGA}}
\end{bmatrix}^T$ is the full vector of shell displacements in IGA DoFs. 

Partial derivatives $\frac{\partial f^\text{obj}}{\partial \mathbf{P}_{}^{\text{IGA}}}$ and $\frac{\partial f^\text{obj}}{\partial \mathbf{u}^\text{IGA}}$ in \eqref{eq:shape-opt-total-derivative-init} can be computed and depend on the objective function combined with extraction matrices,
\begin{align}
    \frac{\partial f^\text{obj}}{\partial \mathbf{P}_{}^{\text{IGA}}} = \mathbf{M}^T \frac{\partial f^\text{obj}}{\partial \mathbf{P}_{}^{\text{FE}}} \quad \text{ and } \quad \frac{\partial f^\text{obj}}{\partial \mathbf{u}^{\text{IGA}}} = \mathbf{M}^T \frac{\partial f^\text{obj}}{\partial \mathbf{u}^{\text{FE}}} \text{ ,} \label{eq:derivative-fobj}
\end{align}
where $\mathbf{P}_{}^{\text{FE}} = \begin{bmatrix}
    \mathbf{P}_{}^{\text{A,FE}} & \mathbf{P}_{}^{\text{B,FE}} & \ldots & \mathbf{P}_{}^{\text{I}^m\text{,FE}}
\end{bmatrix}^T$, $\mathbf{u}^{\text{FE}} = \begin{bmatrix}
    \mathbf{u}^{\text{A,FE}} & \mathbf{u}^{\text{B,FE}} & \ldots & \mathbf{u}^{\text{I}^m\text{,FE}}
\end{bmatrix}^T$. $\mathbf{M} = \text{diag}(\mathbf{M}^\text{A}, \mathbf{M}^\text{B},\ldots,\mathbf{M}^{\text{I}^m})$ is a diagonal block matrix for global extraction. Calculation of partial derivatives in \eqref{eq:derivative-fobj} is automated using FEniCS. Formulation for total derivative $\frac{d \mathbf{P}_{}^{\text{IGA}}}{d \mathbf{P}_{\text{FFD}}}$ is introduced in \eqref{eq:piga-pffd-equation}. As for total derivative $\frac{d \mathbf{u}^\text{IGA}}{d \mathbf{P}_{}^{\text{IGA}}}$, we have the implicit relation between $\mathbf{P}_{}^{\text{IGA}}$ and $ \mathbf{u}^\text{IGA}$,
\begin{align}
    \mathbf{r}_\text{t} = \mathbf{R}_\text{t}(\mathbf{P}_{}^{\text{IGA}}, \mathbf{u}^\text{IGA}) = \frac{\partial W_\text{t}(\mathbf{u}^\text{IGA}, \mathbf{P}_{}^{\text{IGA}})}{\partial \mathbf{u}^\text{IGA}} = \mathbf{0} \text{ ,} \label{eq:uiga-piga-implicit-relation}
\end{align}
where $W_\text{t}$ is the total energy of the non-matching shells defined in \eqref{eq:total-virtual-energy}. Once an updated $\mathbf{P}_{}^{\text{IGA}}$ is obtained, the shell displacements need to be solved using \eqref{eq:total-virtual-energy} until the residual vector $\mathbf{r}_\text{t}$ reaches a tolerance. Thus, $\mathbf{r}_\text{t}$ is supposed to remain as $\mathbf{0}$ despite the change of $\mathbf{P}_{}^{\text{IGA}}$, and we have the following derivative
\begin{align}
    \frac{d \mathbf{r}_\text{t}}{d \mathbf{P}_{}^{\text{IGA}}} = \frac{\partial \mathbf{R}_\text{t}}{\partial \mathbf{P}_{}^{\text{IGA}}} + \frac{\partial  \mathbf{R}_\text{t}}{\partial \mathbf{u}^\text{IGA}} \frac{d \mathbf{u}^\text{IGA}}{d \mathbf{P}_{}^{\text{IGA}}}= \mathbf{0} \text{ ,} \label{eq:derivative-dr-dpiga-0}
\end{align}
and the total derivative $\frac{d \mathbf{u}^\text{IGA}}{d \mathbf{P}_{}^{\text{IGA}}}$ in \eqref{eq:shape-opt-total-derivative-init} is given by
\begin{align}
    \frac{d \mathbf{u}^\text{IGA}}{d \mathbf{P}_{}^{\text{IGA}}} = - \left( \frac{\partial  \mathbf{R}_\text{t}}{\partial \mathbf{u}^\text{IGA}} \right)^{-1} \frac{\partial \mathbf{R}_\text{t}}{\partial \mathbf{P}_{}^{\text{IGA}}} \text{ .} \label{eq:derivative-duiga-dpiga}
\end{align}
Partial derivative $\frac{\partial  \mathbf{R}_\text{t}}{\partial \mathbf{u}^\text{IGA}}$ is equivalent to $\frac{\partial^2 W_\text{t}}{\partial (\mathbf{u}^\text{IGA} )^2}$ and is the stiffness matrix defined in \eqref{eq:nonmatching-linear-system}. Analogously, we use a pair of shell patches to illustrate the formulation of partial derivative $\frac{\partial \mathbf{R}_\text{t}}{\partial \mathbf{P}_{}^{\text{IGA}}}$,
\begin{align}
    \setlength\arraycolsep{6pt}
    \renewcommand\arraystretch{1.5}
    \frac{\partial \mathbf{R}_\text{t}}{\partial \mathbf{P}_{}^{\text{IGA}}} = \begin{bmatrix}
        \frac{\partial^2 W_\text{t}}{\partial \mathbf{u}^\text{A,IGA} \partial \mathbf{P}_{}^{\text{A,IGA}}} & \frac{\partial^2 W_\text{t}}{\partial \mathbf{u}^\text{A,IGA} \partial \mathbf{P}_{}^{\text{B,IGA}}} \\ \frac{\partial^2 W_\text{t}}{\partial \mathbf{u}^\text{B,IGA} \partial \mathbf{P}_{}^{\text{A,IGA}}} & \frac{\partial^2 W_\text{t}}{\partial \mathbf{u}^\text{B,IGA} \partial \mathbf{P}_{}^{\text{B,IGA}}} 
        \end{bmatrix} \text{ .} \label{eq:derivative-pr-ppiga}
\end{align}
Partial derivatives in \eqref{eq:derivative-pr-ppiga} have identical expressions to \eqref{eq:ppartial-wt-ppartial-uaiga} and \eqref{eq:ppartial-wt-ppartial-uaiga-ubiga}.

Extend partial derivatives in \eqref{eq:derivative-pr-ppiga} to shell structures with an arbitrary number of patches, and substitute $\frac{d \mathbf{u}^\text{IGA}}{d \mathbf{P}_{}^{\text{IGA}}}$ in \eqref{eq:shape-opt-total-derivative-init} with \eqref{eq:derivative-duiga-dpiga}, we can obtain the total derivative of the shape optimization
\begin{align}
    \frac{d f^\text{obj}}{d \mathbf{P}_{\text{FFD}}} = \left( \frac{\partial f^\text{obj}}{\partial \mathbf{P}_{}^{\text{IGA}}} - \frac{\partial f^\text{obj}}{\partial \mathbf{u}^\text{IGA}} \left( \frac{\partial  \mathbf{R}_\text{t}}{\partial \mathbf{u}^\text{IGA}} \right)^{-1} \frac{\partial \mathbf{R}_\text{t}}{\partial \mathbf{P}_{}^{\text{IGA}}} \right) \frac{d \mathbf{P}_{}^{\text{IGA}}}{d \mathbf{P}_{\text{FFD}}} \text{ .} \label{eq:shape-opt-total-derivative-final}
\end{align}

\subsection{Sensitivities for thickness optimization} \label{subsec:sensitivities-thickness-opt}
The idea of FFD-based shape update can be applied to shell thickness optimization, where the shell thickness is treated as an extra field of the NURBS control points. We can use \eqref{eq:piga-pffd-equation} to build the relation of the thickness between shell patches and FFD block,
\begin{align}
    \mathbf{t}^{\text{I,IGA}} = \left( \left( (\mathbf{M}_{\text{s}}^{\text{I}})^T\mathbf{M}_\text{s}^\text{I} \right)^{-1}(\mathbf{M}_{\text{s}}^{\text{I}})^T \mathbf{A}_{\text{FFDs}}^{\text{I}} \right) \mathbf{t}_{\text{FFD}} \text{ ,} \label{eq:tiga-tffd-equation}
\end{align}
where $\mathbf{t}_{}^{\text{I,IGA}}$ is the thickness for shell $S^\text{I}$ in IGA DoFs, and $\mathbf{t}_{\text{FFD}}$ is the corresponding thickness field of the FFD block. Subscript $\text{s}$ in $\mathbf{M}_{\text{s}}^{\text{I}}$ and $\mathbf{A}_{\text{FFDs}}^{\text{I}}$ denotes matrices for scalar fields. Note that $\mathbf{t}_{\text{FFD}}$ is not the actual thickness of the B-spline solid geometry, but an extra set of the control points on the FFD block to update the thickness of the non-matching shells. Accordingly, the identical shape update strategy Section \ref{subsubsec:relating-cp-ffd-shells} is applicable to thickness update. FFD-based thickness optimization also offers the benefit that shell thickness remains continuous on the surface--surface intersections.

Replacing control points of the FFD block in \eqref{eq:shape-minimization} with $\mathbf{t}_{\text{FFD}}$, one can have the problem description of thickness optimization. Since both Kirchhoff--Love shell total work $W^\text{A}$ and $W^\text{B}$, and penalty energy $W^\text{AB}_\text{pen}$ involve shell thickness, the total derivative and associated partial derivatives of the thickness optimization problem can be acquired by replacing $\mathbf{P}_{\text{FFD}}$, $\mathbf{P}_{}^{\text{FE}}$, $\mathbf{P}_{}^{\text{IGA}}$ with $\mathbf{t}_{\text{FFD}}$, $\mathbf{t}_{}^{\text{FE}}$ and $\mathbf{t}_{}^{\text{IGA}}$ in Eqs. \eqref{eq:derivative-duiga-dpiga}--\eqref{eq:shape-opt-total-derivative-final}, respectively. The total derivative of FFD-based thickness optimization reads as
\begin{align}
    \frac{d f^\text{obj}}{d \mathbf{t}_{\text{FFD}}} = \left( \frac{\partial f^\text{obj}}{\partial \mathbf{t}_{}^{\text{IGA}}} - \frac{\partial f^\text{obj}}{\partial \mathbf{u}^\text{IGA}} \left( \frac{\partial  \mathbf{R}_\text{t}}{\partial \mathbf{u}^\text{IGA}} \right)^{-1} \frac{\partial \mathbf{R}_\text{t}}{\partial \mathbf{t}_{}^{\text{IGA}}} \right) \frac{d \mathbf{t}_{}^{\text{IGA}}}{d \mathbf{t}_{\text{FFD}}} \text{ .} \label{eq:thickness-opt-total-derivative-final}
\end{align}

In some applications, one may choose to have a constant thickness for each shell patch. This can be easily achieved within the current framework by relating the shell thickness in IGA DoFs to one scalar value $t^\text{I,const}$ as
\begin{align}
    \mathbf{t}^{\text{I,IGA}} = \mathbf{c}^\text{I} t^\text{I,const} \text{ ,} \label{eq:tiga-tconst-relation}
\end{align}
where $\mathbf{c}^\text{I} = \begin{bmatrix} 1 & 1 & \ldots & 1 \end{bmatrix}^T$ is a unit column vector contains $n^{\text{I,IGA}}$ entries. The total derivative for piecewise constant thickness optimization can be obtained by replacing $\frac{d \mathbf{t}_{}^{\text{IGA}}}{d \mathbf{t}_{\text{FFD}}}$ in \eqref{eq:thickness-opt-total-derivative-final} with the following derivative,
\begin{align}
    \frac{d \mathbf{t}_{}^{\text{IGA}}}{d \mathbf{t}^\text{const}} = \text{diag}(\mathbf{c}^\text{A}, \mathbf{c}^\text{B}, \ldots, \mathbf{c}^{\text{I}^m}) \text{ .} \label{eq:tiga-tconst-full}
\end{align}
The FFD block is not needed in piecewise constant only thickness optimization.

These two approaches can be combined together to achieve a more realistic design, where specific sections of the structure necessitate a continuous thickness distribution while constant thickness is better suited for other patches. In the implementation, shell patches can be separated into various groups. One group comprises the shell patches immersed within an FFD block allowing for a continuous thickness distribution. On the other hand, the shell patches not contained in an FFD block are assumed to have a constant thickness. Moreover, shell patches originating from different FFD blocks would exhibit discontinuous thickness at their intersections. Therefore, the combined thickness optimization approach provides more flexibility.


\subsection{Implementation}\label{subsec:implementation}

This section introduces the open-source implementation for the proposed optimization method in Sections \ref{subsec:nonmatching-shell-update-ffd}--\ref{subsec:sensitivities-thickness-opt}. We employ the open-source Python library PENGoLINS \cite{Zhao2022} for IGA of non-matching shell structures. PENGoLINS utilizes functionalities in pythonOCC \cite{paviot2018pythonocc} to approximate patch intersections in CAD geometries. Additionally, it leverages advanced code generation in FEniCS \cite{Logg2012} and extraction operators in tIGAr \cite{Kamensky2019} to perform IGA. CAD models in IGES or STEP formats can be imported into PENGoLINS, making them directly available for IGA. The code generation capabilities and computer algebra in FEniCS allow us to obtain partial derivatives in \eqref{eq:partial-wt-partial-uaiga}--\eqref{eq:ppartial-wt-ppartial-uaiga-ubiga}, as well as \eqref{eq:derivative-fobj} automatically. Integrating computed analytical derivatives with matrices $\mathbf{M}^{} = \text{diag}(\mathbf{M}^{\text{A}},\mathbf{M}^{\text{B}}, \ldots, \mathbf{M}^{\text{I}^m})$ and $\mathbf{A}_{\text{FFD}} = \text{diag}(\mathbf{A}_{\text{FFD}}^{\text{A}},\mathbf{A}_{\text{FFD}}^{\text{B}}, \ldots, \mathbf{A}_{\text{FFD}}^{\text{I}^m})$, we can calculate the total derivatives in \eqref{eq:shape-opt-total-derivative-final} and \eqref{eq:thickness-opt-total-derivative-final} for shape and thickness optimization.

With the availability of automated derivatives, we use the open-source framework CSDL \cite{gandarillas2022graph} for conducting gradient-based large-scale optimization. Each partial derivative computed from FEniCS, in combination with predefined matrices, is passed to a CSDL sub-model, which contains explicit or implicit operations. A system-level model is created to connect all sub-models to enable efficient and modularized design optimization. For the FFD-based shape optimization problem, which involves minimizing the internal energy of a non-matching shell structure, four essential sub-models are required and listed in Table \ref{tab:shape-opt-csdl-models}. All the derivatives required for these models are defined in Sections \ref{subsec:computation-nonmatching-shell} and \ref{subsec:sensitivities-shape-opt}.
\begin{table}[!htbp]
    \centering
    \setlength{\tabcolsep}{8pt} 
    \def\arraystretch{1.5}%
    \begin{tabular}{c c c c c}
        \toprule
        Model &  Input(s) & Output & Operation & Derivative(s)  \\
        \midrule
        \makecell{CP FFD to FE} & $\mathbf{P}_{\text{FFD}}$ &  $\mathbf{P}^{\text{FE}}$& Explicit & $\mathbf{A}_{\text{FFD}}^{}$ \\
        \makecell{CP FE to IGA} & $\mathbf{P}^{\text{FE}}$ & $\mathbf{P}^{\text{IGA}}$ & Implicit & $\mathbf{M}^{}$\\
        \makecell{Disp} & $\mathbf{P}^{\text{IGA}}$ & $\mathbf{u}^{\text{IGA}}$ & Implicit & $\frac{\partial \mathbf{R}_{\text{t}}}{\partial \mathbf{P}^{\text{IGA}}}$, $\frac{\partial \mathbf{R}_{\text{t}}}{\partial \mathbf{u}^{\text{IGA}}}$ \\
        \makecell{Int energy} & $\mathbf{P}^{\text{IGA}}$, $\mathbf{u}^{\text{IGA}}$ & $W_{\text{int}}$ & Explicit & $\frac{\partial W_{\text{int}}}{\mathbf{P}^{\text{IGA}}}$, $\frac{\partial W_{\text{int}}}{\mathbf{u}^{\text{IGA}}}$\\
        \bottomrule
    \end{tabular}
    \caption{Essential CSDL models for FFD-based shape optimization.}
    \label{tab:shape-opt-csdl-models}
\end{table}

The models listed in Table \ref{tab:shape-opt-csdl-models} and other constraint models are modularized in a Python library named GOLDFISH (Gradient-based Optimization, Large-scale Design Framework for Isogeometric SHells).  This library provides users with the option to use predefined models or create custom models to build the system-level model for their specific optimization needs. The library also includes an open-source optimizer, SLSQP \cite{kraft1988software}, which can be utilized for shape or thickness optimization of non-matching shells. For optimization problems involving a large number of design variables, the SNOPT \cite{gill2005snopt} optimizer is available for faster convergence. Furthermore, we have integrated the analysis framework with another optimization library OpenMDAO \cite{gray2019openmdao}, which is developed by NASA. This integration further enhances the accessibility of the GOLDFISH library, enabling users to benefit from the collective expertise of the OpenMDAO community. The GOLDFISH library is hosted on a GitHub repository \cite{goldfish-code} and open to the public.


\section{Benchmark Problems}\label{sec:benchmark}
A series of benchmark problems are considered to verify the effectiveness of the optimization method. Sections \ref{seubsec:arch-shape-opt}--\ref{subsec:T-beam-shape-opt} illustrate that baseline non-matching shell structures with arbitrary intersections are able to accurately converge to the analytical optimum. Section \ref{subsec:plate-th-opt} studies the capability and flexibility of the framework for thickness optimization.

\subsection{Arch shape optimization}\label{seubsec:arch-shape-opt}
\begin{figure}[!htbp]
    \centering
    \begin{subfigure}[!htbp]{0.49\textwidth}
        \includegraphics[width=\textwidth]{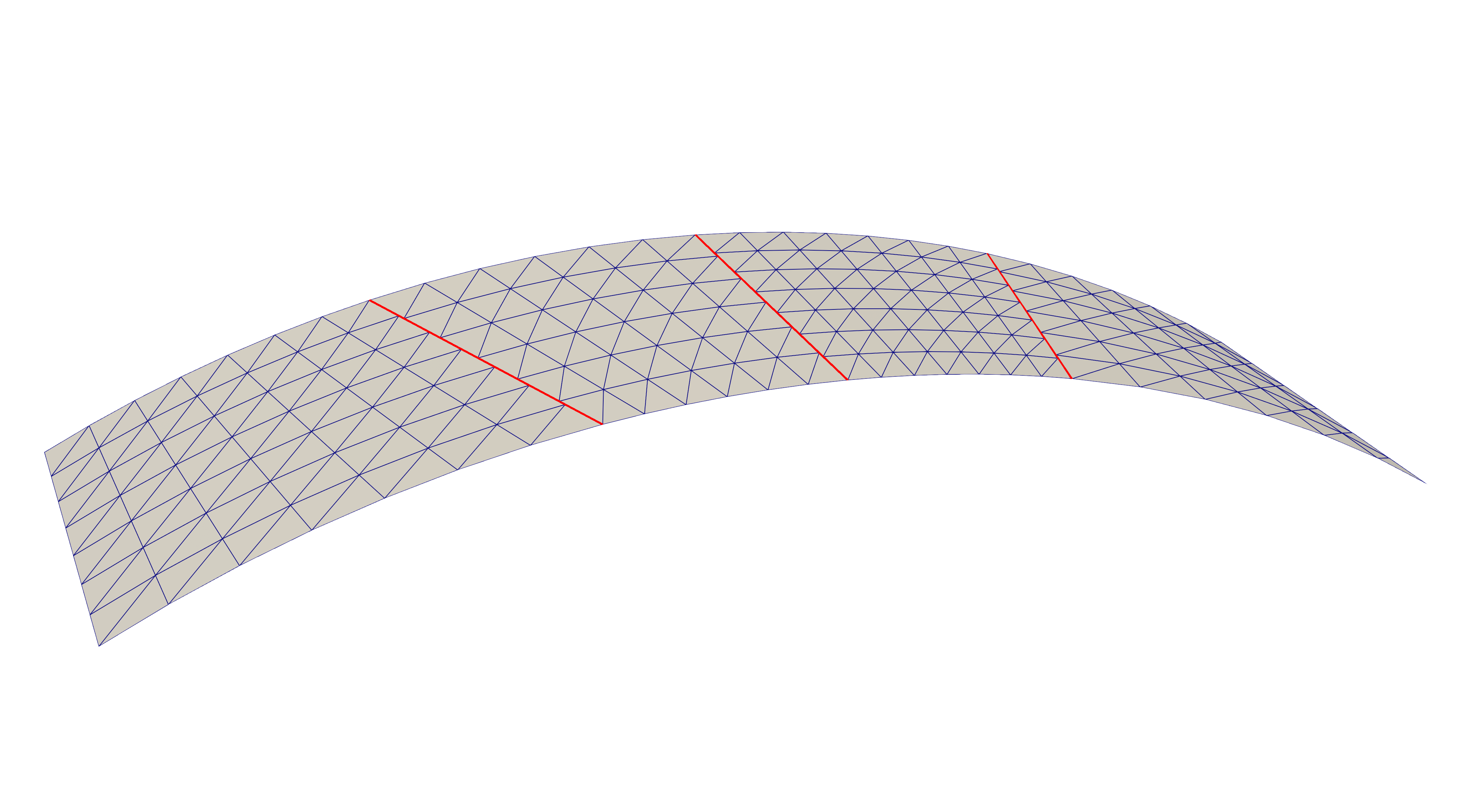}
        \caption{}
        \label{subfig:nonmatcing-arch-geom-init}
    \end{subfigure}
    \hfill
    \begin{subfigure}[!htbp]{0.49\textwidth}
        \includegraphics[width=\textwidth]{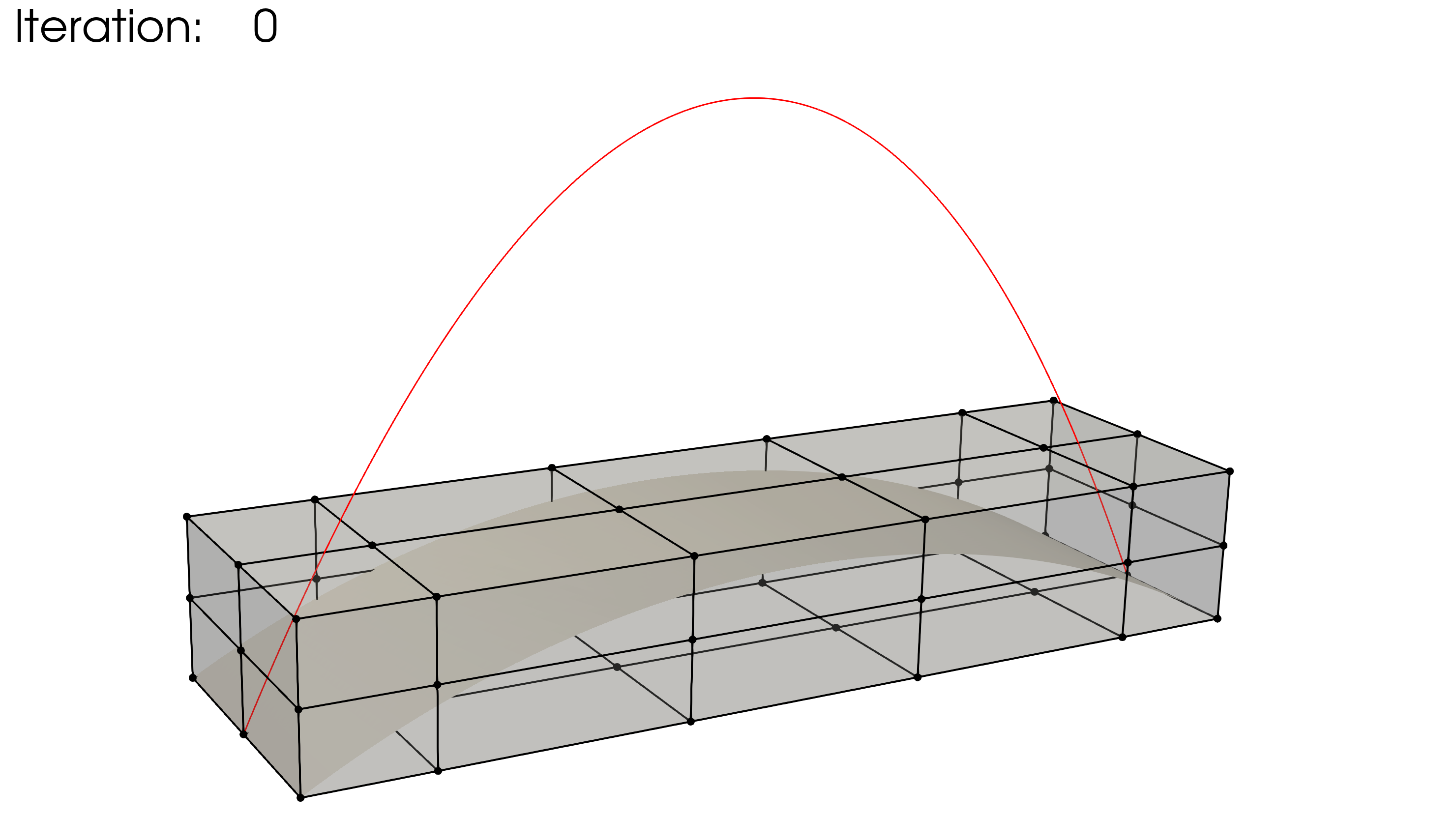}
        \caption{}
        \label{subfig:arch-shape-opt-ffd-init}
    \end{subfigure}
    \caption{(a) Baseline geometry of a non-matching arch consisting of four NRUBS patches, three surface--surface intersections are indicated with red lines. (b) Initial configuration of the arch immersed in an FFD block, where black lines and dots denote the control net. The analytical optimal shape is plotted with a red curve. }
    \label{fig:arch-shape-opt-init}
\end{figure}

An arch fixed at two ends and subjected to a constant downward load per unit horizontal length is modeled by a Kirchhoff--Love shell theory. Detailed problem definitions can be found in \cite[Section 8]{KIENDL2014148}. To test the effectiveness of FFD-based shape optimization for the non-matching shells approach, we model the arch using four NURBS patches with three intersections, where the arch geometry in the baseline configuration is shown in Figure \ref{subfig:nonmatcing-arch-geom-init}. We immerse the arch geometry in a trivariate B-spline block in the initial configuration, as is illustrated in Figure \ref{subfig:arch-shape-opt-ffd-init}. The analytical optimal solution is given by a quadratic parabola, where the ratio between the height of the arch and the horizontal distance of two fixed edges is 0.54779.

\begin{figure}[!htb]\centering
    \includegraphics[width=1.0\textwidth]{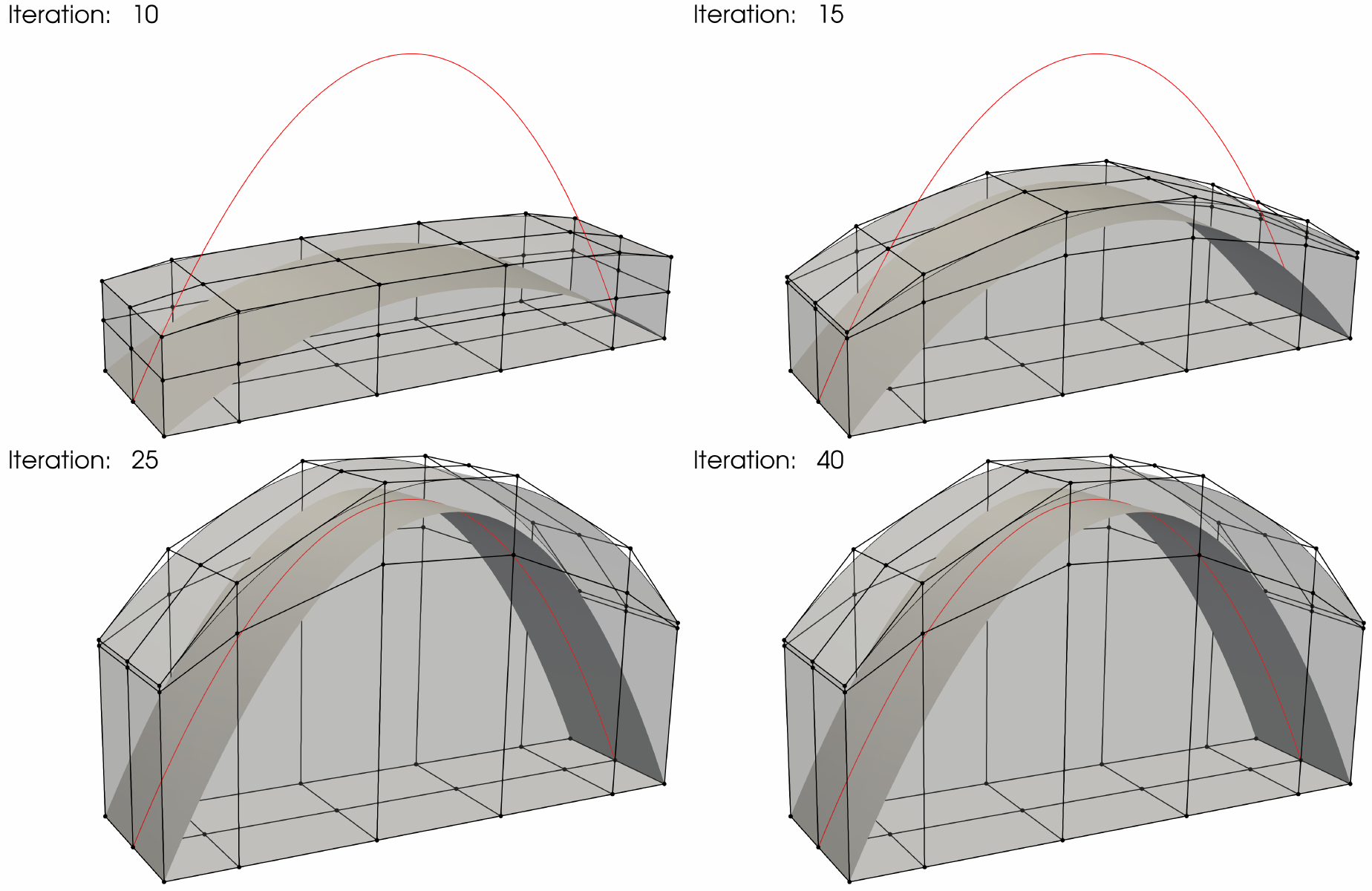}
    \caption{Snapshots of non-matching arch shape optimization.}
    \label{fig:arch-shape-opt-history}
\end{figure}

We use quadratic B-spline for the FFD block in all three directions, $p_\text{FFD}=2$. The arch shell patches are described by cubic NURBS surfaces, $p_\text{sh}=3$, with 1086 DoFs in total. This benchmark problem minimizes the internal energy of the shell structure, with vertical coordinates of the control points of the FFD block being the design variables. From the control net in Figure \ref{subfig:arch-shape-opt-ffd-init}, it can be observed that there are 54 design variables. Two constraints are applied to this problem. The first constraint ensures that the lines in the FFD control net are parallel to the axial direction of the arch, keeping the arch devoid of tilting or twisting during the optimization process. The second constraint fixes the bottom layer of FFD control points so that the two edges of the arch remain in the initial position. We use the SLSQP optimizer with a tolerance of $10^{-12}$ to perform the optimization, snapshots of the optimization iteration are demonstrated in Figure \ref{fig:arch-shape-opt-history}.

The arch converges to the analytical optimum shape after 40 iterations. The shape of the FFD block in the final configuration is shown in Figure \ref{fig:arch-shape-opt-history}. As anticipated, the optimized arch geometry is still contained in the FFD block. The height to base span ratio in this problem is measured as 0.54748, exhibiting a relative error of 0.057\% compared to the exact value. Considering the coarse discretization of the arch geometry, the results are encouraging.

\subsection{Tube shape optimization}\label{subsec:tube-shape-opt}
\begin{figure}[!htb]
    \centering
    \begin{subfigure}[!htb]{0.42\textwidth}
        \centering
        \includegraphics[height=1.8in]{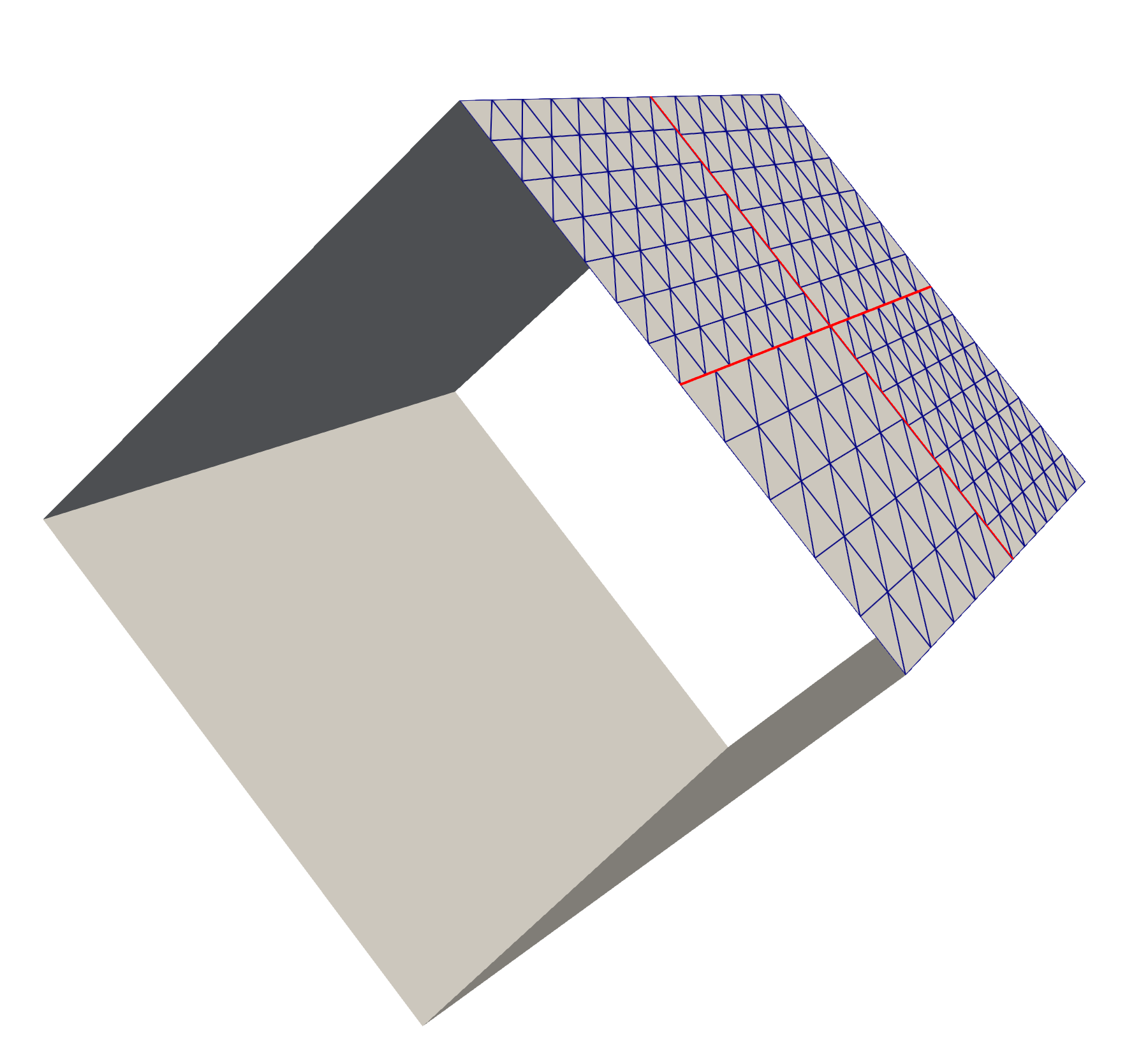}
        \caption{}
        \label{subfig:nonmatcing-tube-geom-init}
    \end{subfigure}
    \hfill
    \begin{subfigure}[!htb]{0.54\textwidth}
        \centering
        \includegraphics[height=1.8in]{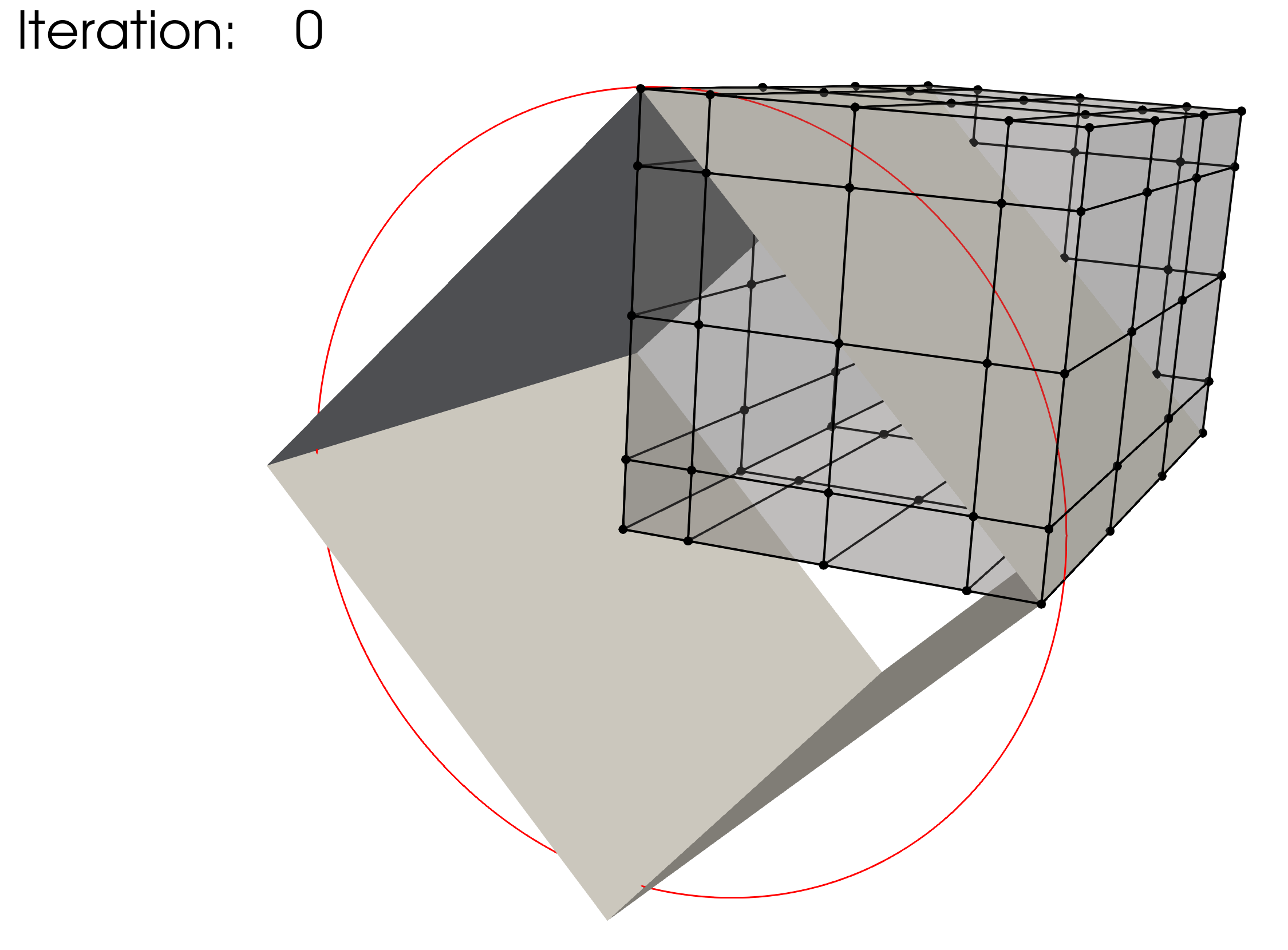}
        \caption{}
        \label{subfig:tube-shape-opt-ffd-init}
    \end{subfigure}
    \caption{(a) Baseline geometry of the square tube, a quarter of the tube is modeled using four non-matching B-spline patches with four intersections. (b) Initial configuration of the FFD block with control net. The optimal cross-section is depicted by a red circle.}
    \label{fig:tube-shape-opt-init}
\end{figure}

A square tube in the baseline configuration is subjected to an internal constant pressure. The optimal shape is given by a cylindrical tube \cite[Section 7]{KIENDL2014148}. We use four cubic NURBS surfaces to model one-quarter of the square tube, where the initial geometry is illustrated in Figure \ref{subfig:nonmatcing-tube-geom-init}. The square tube geometry contains four non-matching intersections with 1035 DoFs in total, and symmetric boundary conditions are applied in the analysis. The geometry is immersed in a cubic B-spline block to perform FFD-based shape optimization, as shown in Figure \ref{subfig:tube-shape-opt-ffd-init}, where the cross-section of the optimal shape is indicated by the red curve. In this example, the horizontal and vertical coordinates of the FFD control points are design variables, totaling 200 design variables. Similar constraints are employed in this problem as those in the arch shape optimization. Control points on the left and bottom layers are fixed, where the lines in the FFD control net that are parallel to the tube axis maintain their orientations. 

\begin{figure}[!htb]\centering
    \includegraphics[width=1.0\textwidth]{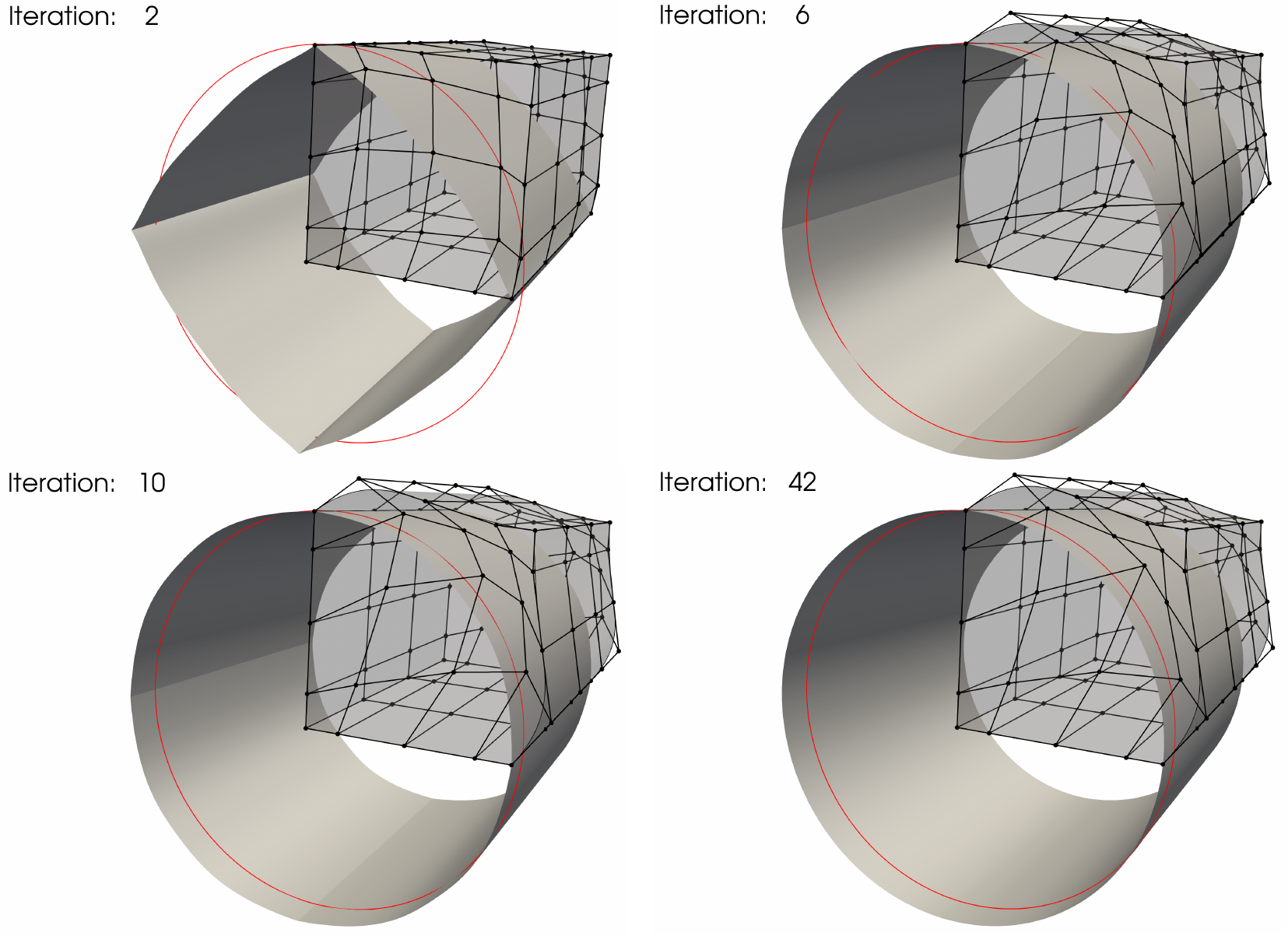}
    \caption{Iteration history for tube shape optimization under follower pressure.}
    \label{fig:tube-shape-opt-history}
\end{figure}

The optimization problem converges successfully after 42 iterations using SLSQP optimizer with a tolerance of $10^{-12}$, and snapshots are depicted in Figure \ref{fig:tube-shape-opt-history}. As expected, the initial square tube converges to the cylindrical tube.

\subsection{T-beam shape optimization}\label{subsec:T-beam-shape-opt}
\begin{figure}[!htb]
    \centering
    \begin{subfigure}[!htb]{0.45\textwidth}
        \centering
        \includegraphics[height=1.8in]{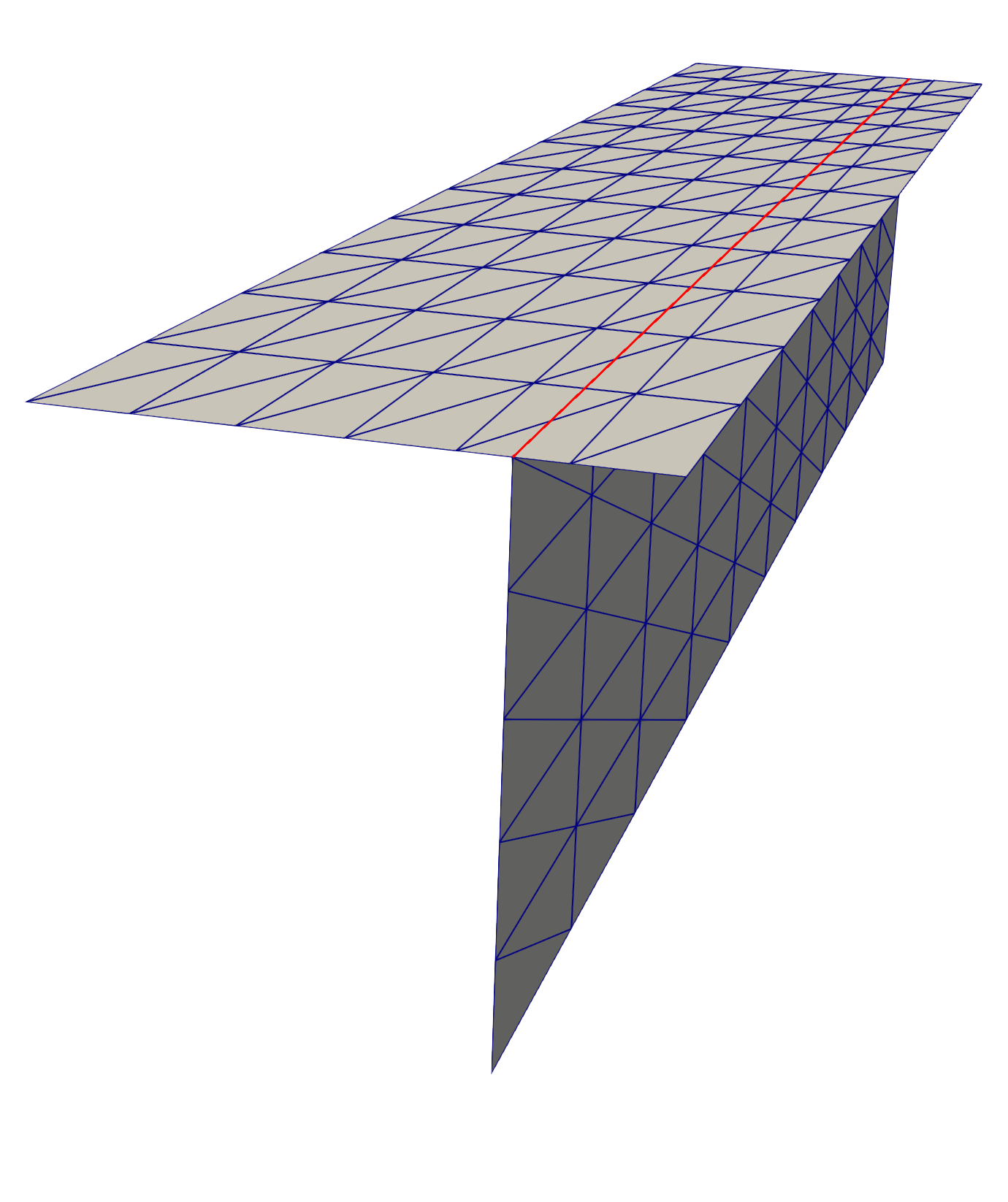}
        \caption{}
        \label{subfig:nonmatcing-tbeam-geom-init}
    \end{subfigure}
    \hfill
    \begin{subfigure}[!htb]{0.54\textwidth}
        \centering
        \includegraphics[height=1.8in]{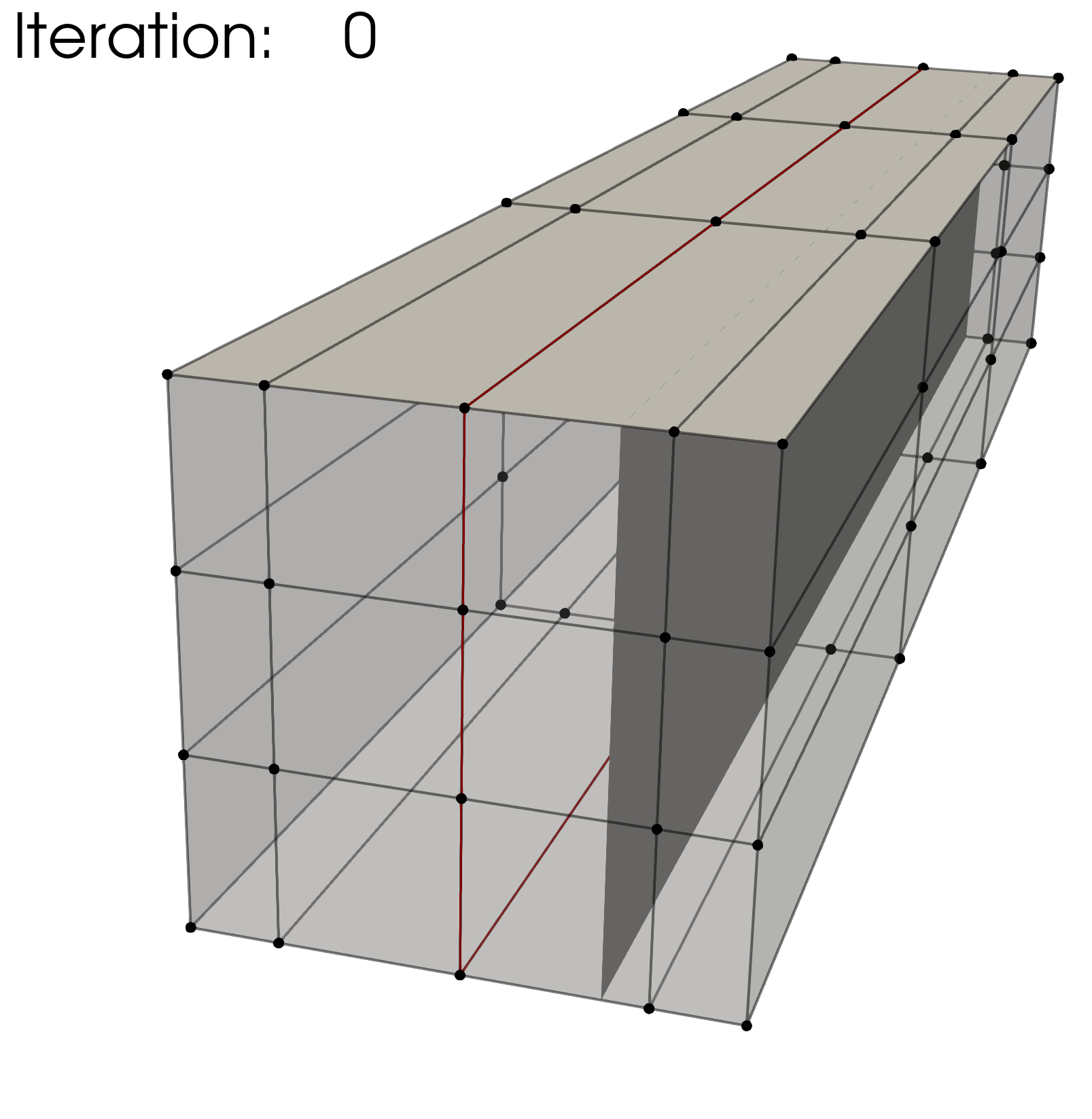}
        \caption{}
        \label{subfig:tbeam-shape-opt-ffd-init}
    \end{subfigure}
    \caption{(a) Baseline configuration of a T-beam whose vertical patch is at the three-quarter position of the horizontal patch. (b) The T-beam is placed in an FFD B-spline block. The optimal position of the vertical patch is depicted with red lines.}
    \label{fig:tbeam-shape-opt-init}
\end{figure}

In this section, a T-beam is considered to test the method for shell structures with intersections in the middle. We model a T-beam using two NURBS patches which are shown in Figure \ref{subfig:nonmatcing-tbeam-geom-init}. In the baseline design, the vertical patch in the T-beam is located at the three-quarter position, where the mismatched intersection is indicated with a red line.

\begin{figure}[!htb]\centering
    \includegraphics[width=0.8\textwidth]{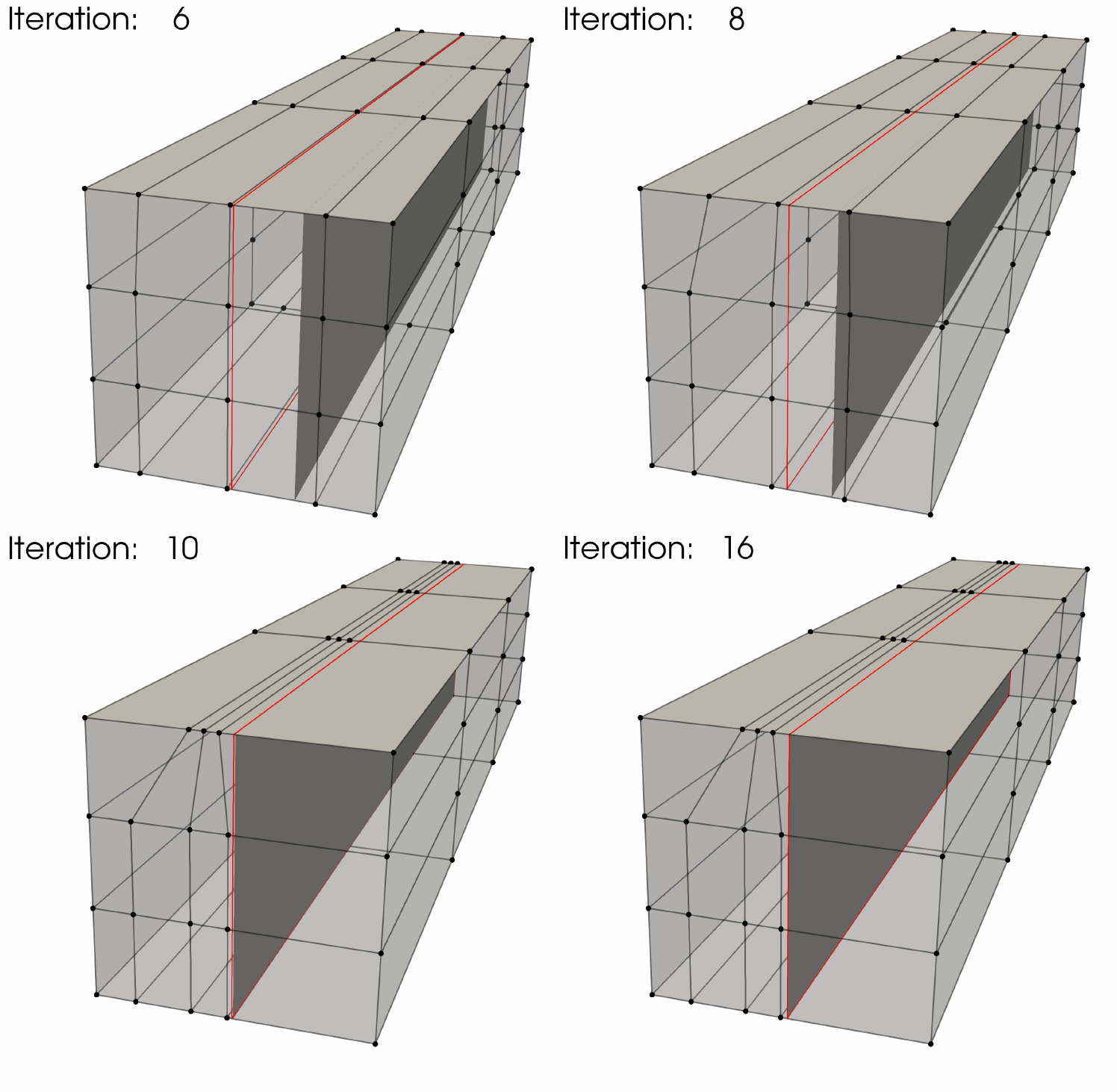}
    \caption{Screenshots of T-beam shape optimization process.}
    \label{fig:tbeam-shape-opt-history}
\end{figure}

In this benchmark problem, we aim to minimize the internal energy of the T-beam by updating the horizontal coordinates of shell patches' control points. Subsequently, the T-beam is placed in a cubic FFD block, where the horizontal coordinates of the FFD block's control points serve as design variables. The rear end of the T-beam is fixed. Given a downward distributed load, the optimal shape of the T-beam is expected to have a junction at the center of the horizontal patch under a constant volume constraint. The optimal position of the vertical patch is depicted in Figure \ref{subfig:tbeam-shape-opt-ffd-init} with red lines. The right and left layers of the FFD block's control points are fixed by employing equality constraints, and control points stay collinear in the axial direction. Using a tolerance of $10^{-15}$, the SLSQP optimizer converges successfully after 16 iterations, and the optimization process is shown in Figure \ref{fig:tbeam-shape-opt-history}.

In Figure \ref{fig:tbeam-shape-opt-history}, the T-beam converges to the optimal solution, indicating that the proposed method is effective for non-matching shell structures with intersections. A more complicated demonstration is presented in Section \ref{subsec:evtol-wing-shthopt}, where the geometry contains curved T-junctions.

\subsection{Thickness optimization of a clamped plate}\label{subsec:plate-th-opt}
As stated in Section \ref{subsec:sensitivities-thickness-opt}, the FFD-based optimization methodology can also be applied to thickness optimization. In the following section, we first demonstrate a piecewise constant thickness optimization for a clamped non-matching plate, in which the FFD block is not needed. Subsequently, we proceed to perform the variable thickness optimization for the same geometry.

\subsubsection{Piecewise constant thickness optimization}\label{subsubsec:plate-th-opt-const}
\begin{figure}[!htb]
    \centering
    \begin{subfigure}[!htb]{0.45\textwidth}
        \centering
        \includegraphics[height=1.8in]{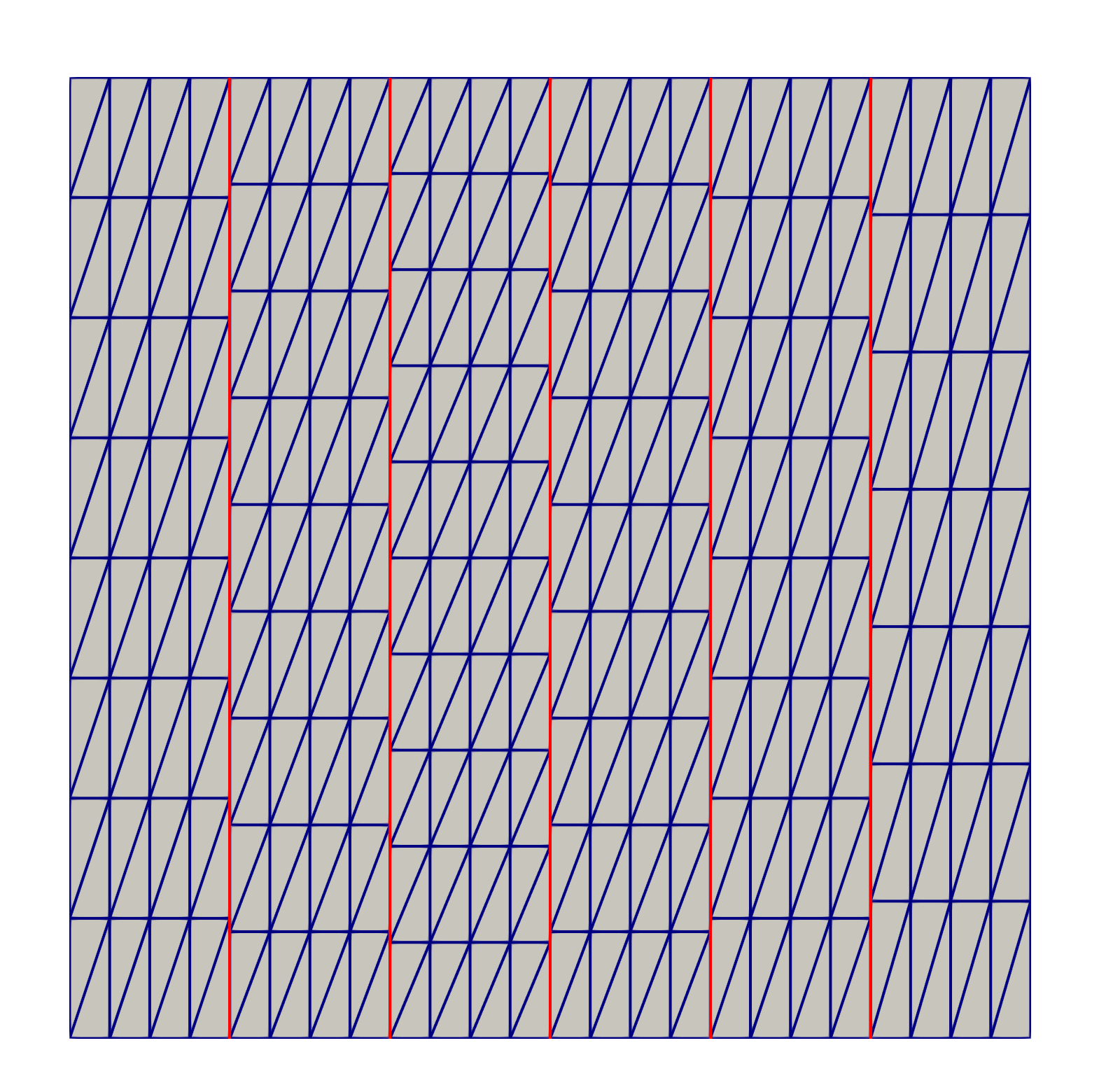}
        \caption{}
        \label{subfig:nonmatcing-plate-geom-init}
    \end{subfigure}
    \hfill
    \begin{subfigure}[!htb]{0.54\textwidth}
        \centering
        \includegraphics[height=1.8in]{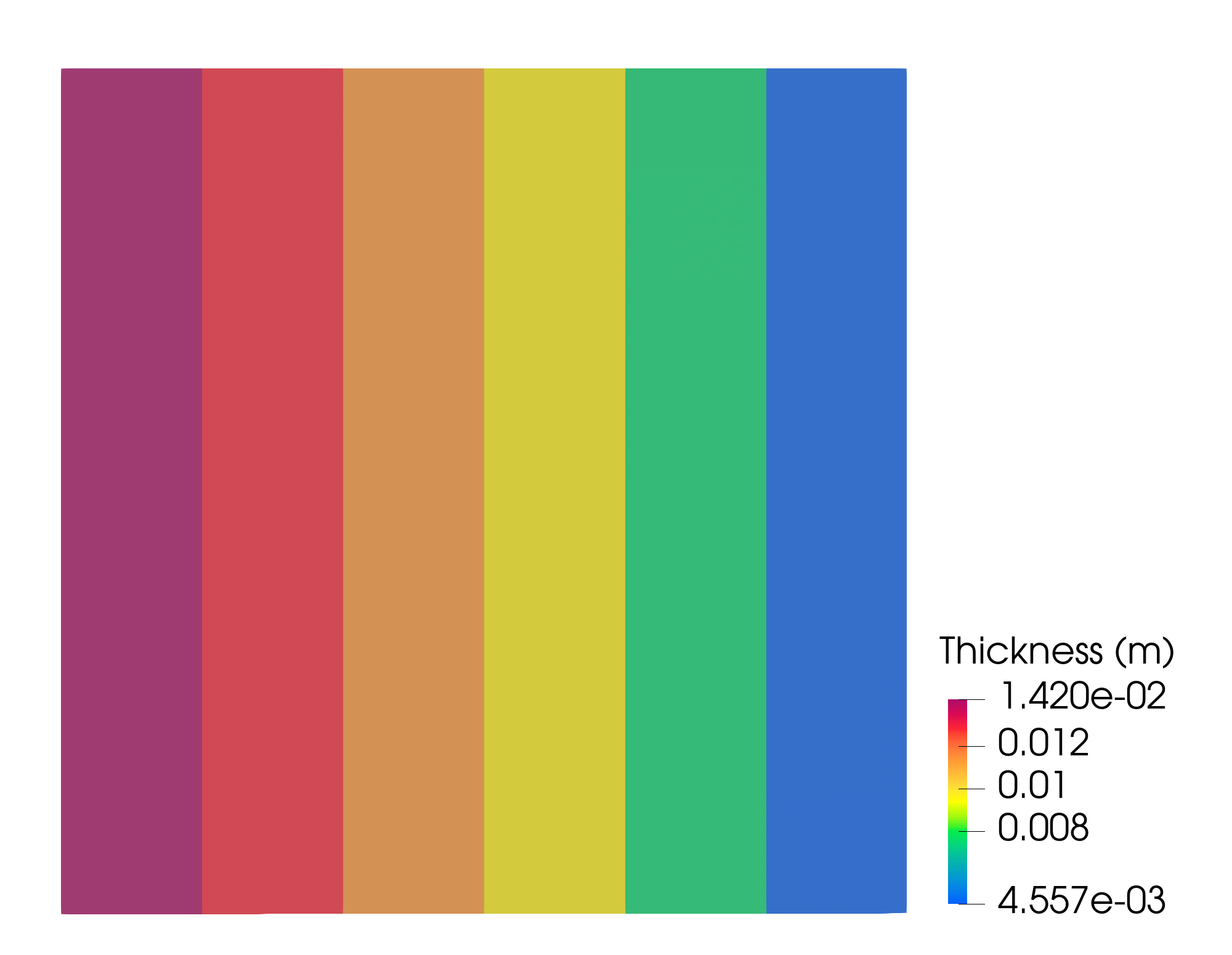}
        \caption{}
        \label{subfig:plate-const-thopt-final}
    \end{subfigure}
    \caption{(a) A unit square plate consisting of six non-matching patches, intersections are indicated with red lines. (b) Final plate thickness for piecewise constant thickness optimization.}
    \label{fig:plate-const-thopt}
\end{figure}

For the thickness optimization example, a unit square plate composed of six cubic non-matching NURBS surfaces is considered. The geometry, which is shown in Figure \ref{subfig:nonmatcing-plate-geom-init}, exhibits 5 intersections with a total of 1449 DoFs. We apply a clamped boundary condition on the left side and with a line force applied to the right side in the normal direction of the plate. All patches of the plate have an initial thickness of 0.01 m. Using the strategy introduced in Section \ref{subsec:sensitivities-thickness-opt}, we perform piecewise constant thickness optimization for the clamped plate to minimize the internal energy under the constant volume constraint. This problem only has 6 design variables. In this and the following sections, the SNOPT optimizer is used for faster convergence. The optimal thickness is plotted in Figure \ref{subfig:plate-const-thopt-final}.

The observed optimal piecewise constant thickness in Figure \ref{subfig:plate-const-thopt-final} shows material redistributes toward the clamped side, which provides enhanced support to the plate. The internal energy in the final configuration is 37.17\% less than the baseline configuration.

\subsubsection{Variable thickness optimization}\label{subsubsec:plate-th-opt-var}
\begin{figure}[!htb]
    \centering
    \begin{subfigure}[!htb]{0.45\textwidth}
        \centering
        \includegraphics[width=\textwidth]{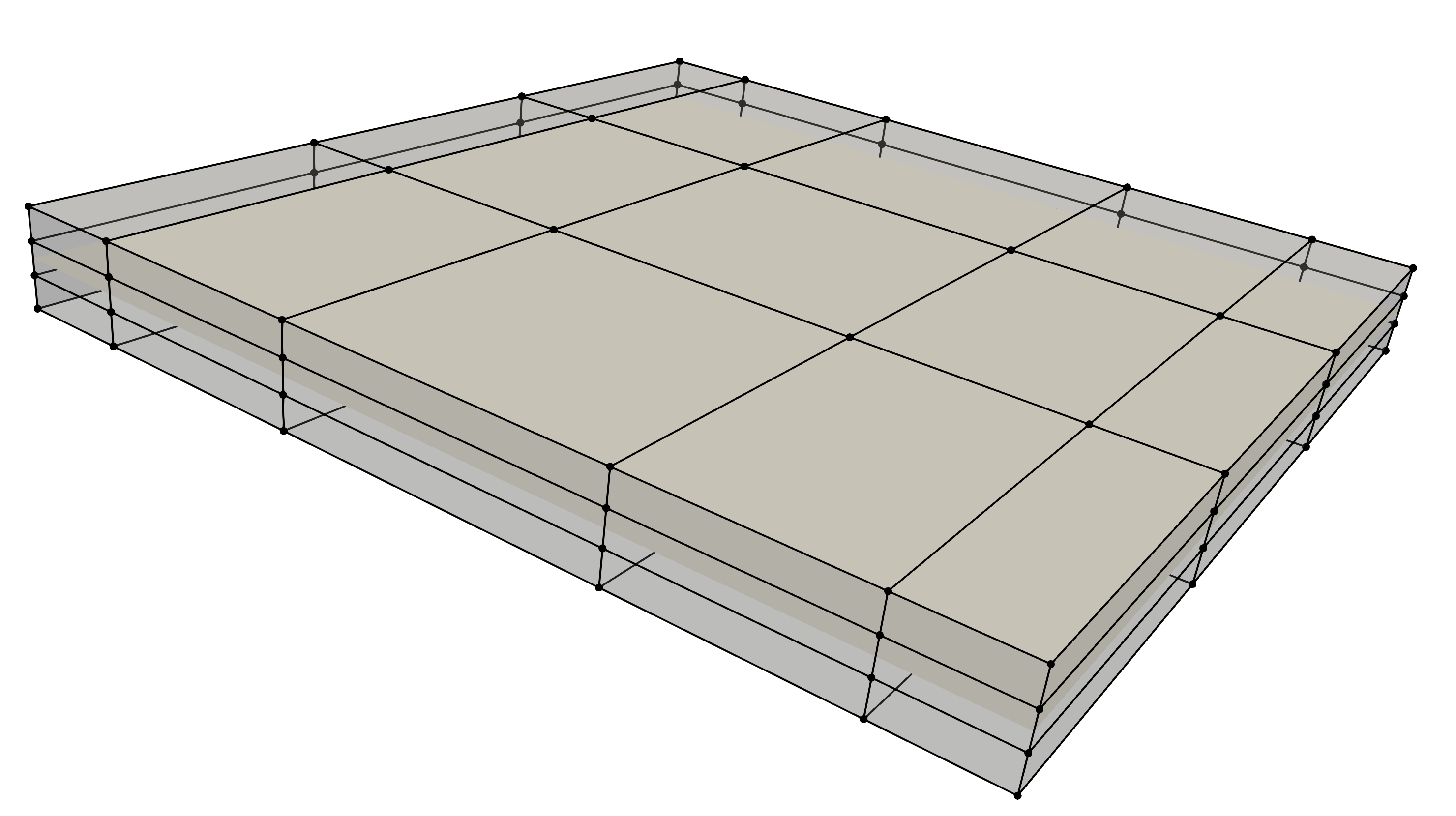}
        \caption{}
        \label{subfig:plate-var-thopt-ffd-init}
    \end{subfigure}
    \hfill
    \begin{subfigure}[!htb]{0.54\textwidth}
        \centering
        \includegraphics[height=1.8in]{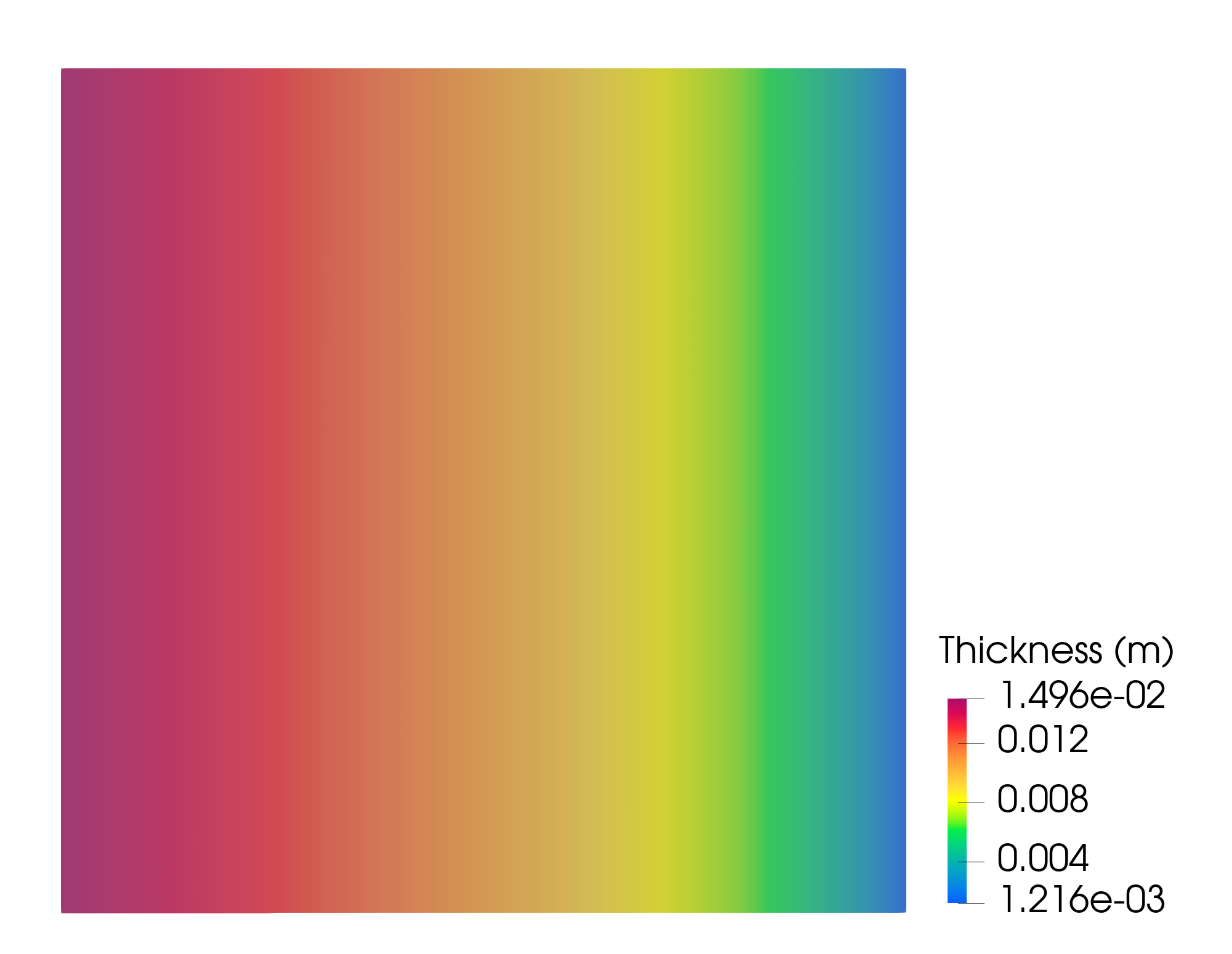}
        \caption{}
        \label{subfig:plate-var-thopt-final}
    \end{subfigure}
    \caption{(a) The non-matching plate is immersed in an FFD block. (b) Optimized thickness distribution using the FFD-based approach.}
    \label{fig:plate-var-thopt}
\end{figure}

We now perform variable thickness optimization for the non-matching plate. The plate is placed in a cubic B-spline block, which is shown in Figure \ref{subfig:plate-var-thopt-ffd-init}. Besides the constant volume constraint, we only optimize the plate thickness in one direction that is perpendicular to the intersections. 

\begin{figure}[!htb]
    \centering
    \begin{subfigure}[!htb]{0.49\textwidth}
        \centering
        \includegraphics[width=\textwidth]{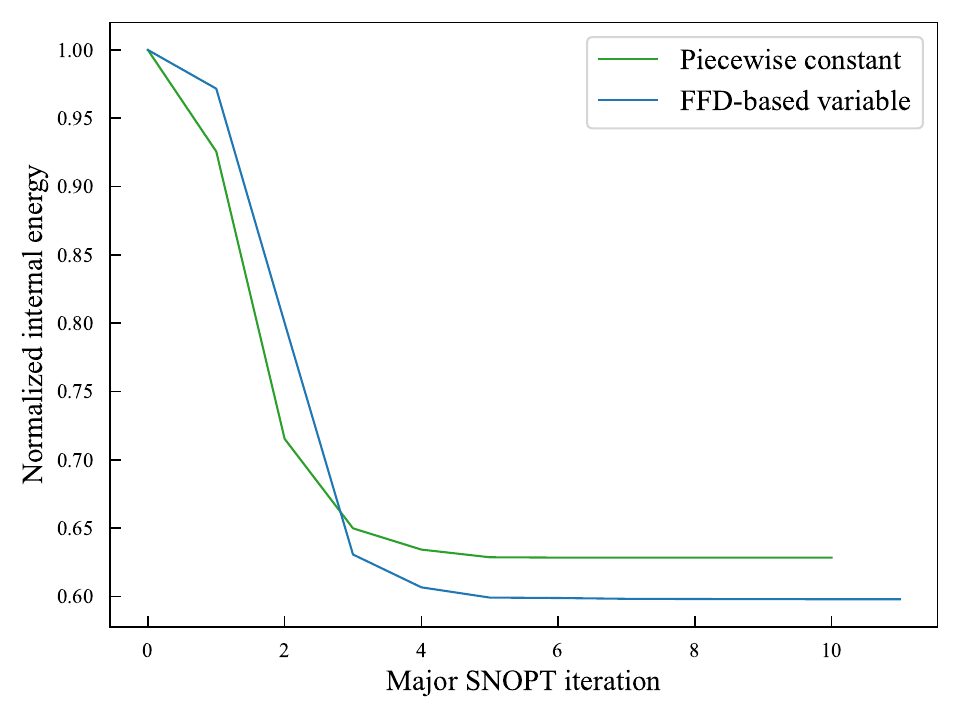}
        \caption{}
        \label{subfig:plate-thopt-wint-compare}
    \end{subfigure}
    \hfill
    \begin{subfigure}[!htb]{0.49\textwidth}
        \centering
        \includegraphics[width=\textwidth]{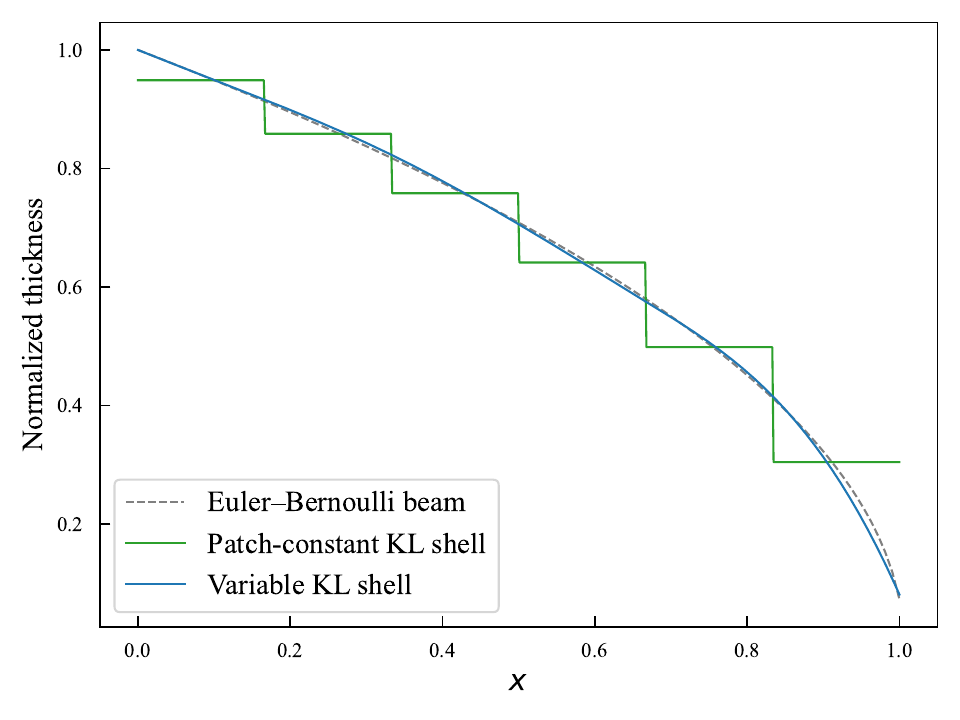}
        \caption{}
        \label{subfig:plate-thopt-thickness-compare}
    \end{subfigure}
    \caption{(a) Optimization process of normalized internal energy for two approaches. (b) Cross-sectional view of piecewise constant thickness and variable thickness, and comparison with Euler-Bernoulli beam thickness optimization.}
    \label{fig:plate-thopt-compare}
\end{figure}

With the FFD-based approach, the continuity of shell thickness is maintained at all intersections. Figure \ref{subfig:plate-var-thopt-final} depicts the converged solution, where a smooth thickness distribution is observed. The smooth thickness profile offers an improved design compared to the piecewise constant thickness approach, resulting in a 40.20\% reduction of the initial internal energy. A comparison of optimization iteration of normalized internal energy between the two methods is illustrated in Figure \ref{subfig:plate-thopt-wint-compare}. The FFD-based thickness optimization approach converges to a smaller internal energy.

To validate the proposed method, we compare the continuous thickness profile to an optimal thickness configuration of a cantilever beam \cite{beam-thickness-opt}. The cantilever beam is modeled using the Euler--Bernoulli beam theory, where a point load is applied to the free end. Since the Kirchhoff--Love shell is an extension of the Euler--Bernoulli beam, both models are expected to yield identical thickness distributions. The normalized thickness profiles of these two models, along with the piecewise constant thickness profile, are plotted in Figure \ref{subfig:plate-thopt-thickness-compare}. A good agreement is observed between the variable thickness of the Kirchhoff--Love shell at the center line and the Euler--Bernoulli beam. Meanwhile, the cross-sectional view of the piecewise constant thickness shows a similar trend to the Euler--Bernoulli beam, albeit with discontinuities at the intersections.

We then investigate the effect of basis function order of continuity in the FFD block. Using the same knots vectors as illustrated in Figure \ref{subfig:plate-var-thopt-ffd-init}, we increase the order of the B-spline basis functions from linear ($C^0$) to quartic ($C^3$) and compare the amounts of reduced internal energy relative to the baseline configuration. These data points are summarized in Table \ref{tab:plate-wint-reductione-pffd}. The results presented in Table \ref{tab:plate-wint-reductione-pffd} indicate that an FFD block with quadratic B-spline basis functions can achieve a better optimal thickness distribution for the clamped plate. The internal energy with quadratic FFD block only exhibits 0.27\% of relative difference compared to the quartic FFD block. Table \ref{tab:plate-wint-reductione-pffd} also suggests that elevating the order of continuity of the FFD block leads to better designs with lower internal energy, particularly when transitioning from linear to quadratic B-spline basis functions.

\begin{table}[!htbp]
    \centering
    \setlength{\tabcolsep}{8pt} 
    \begin{tabular}{c c c c c}
        \toprule
        $p_{\text{FFD}}$ &  1 & 2 & 3 & 4\\
        \midrule
        \makecell{Internal energy \\reduction (\%)} & 39.76 & 40.11 & 40.20 & 40.22 \\
        \bottomrule
    \end{tabular}
    \caption{Reduction of internal energy of the clamped plate for different degrees of the FFD block.}
    \label{tab:plate-wint-reductione-pffd}
\end{table}
The plate example demonstrates that both piecewise constant and variable thickness optimization can be conducted in the proposed framework. One can select desired thickness distribution, or a mixed approach of these two, demonstrated in Section \ref{sec:application}, based on the physical conditions and problem requirements to obtain an optimal design.


\section{Application to aircraft wings}\label{sec:application}
Two aircraft wing design optimization problems are performed. Section \ref{subsec:pegasus-wing-thopt} considers a Parallel Electric-Gas Architecture with Synergistic Utilization Scheme (PEGASUS) wing thickness optimization problem, and Section \ref{subsec:evtol-wing-shthopt} demonstrates a simultaneous shape and thickness optimization for an Electric Vertical Takeoff and Landing (eVTOL) aircraft wing. Both applications illustrate that a design-analysis-optimization workflow can be achieved with the proposed framework.

\subsection{PEGASUS wing thickness optimization}\label{subsec:pegasus-wing-thopt}
For the PEGASUS wing problem, we first verify the accuracy of the structural analysis using PENGoLINS for a shell structure with a large number of patches and intersections. Then two types of thickness design optimizations are performed in the following section.

\subsubsection{Structural analysis of the PEGASUS wing}\label{subsubsec:pegasus-wing-structural-analysis}
\begin{figure}[!htb]\centering
    \includegraphics[width=0.8\textwidth]{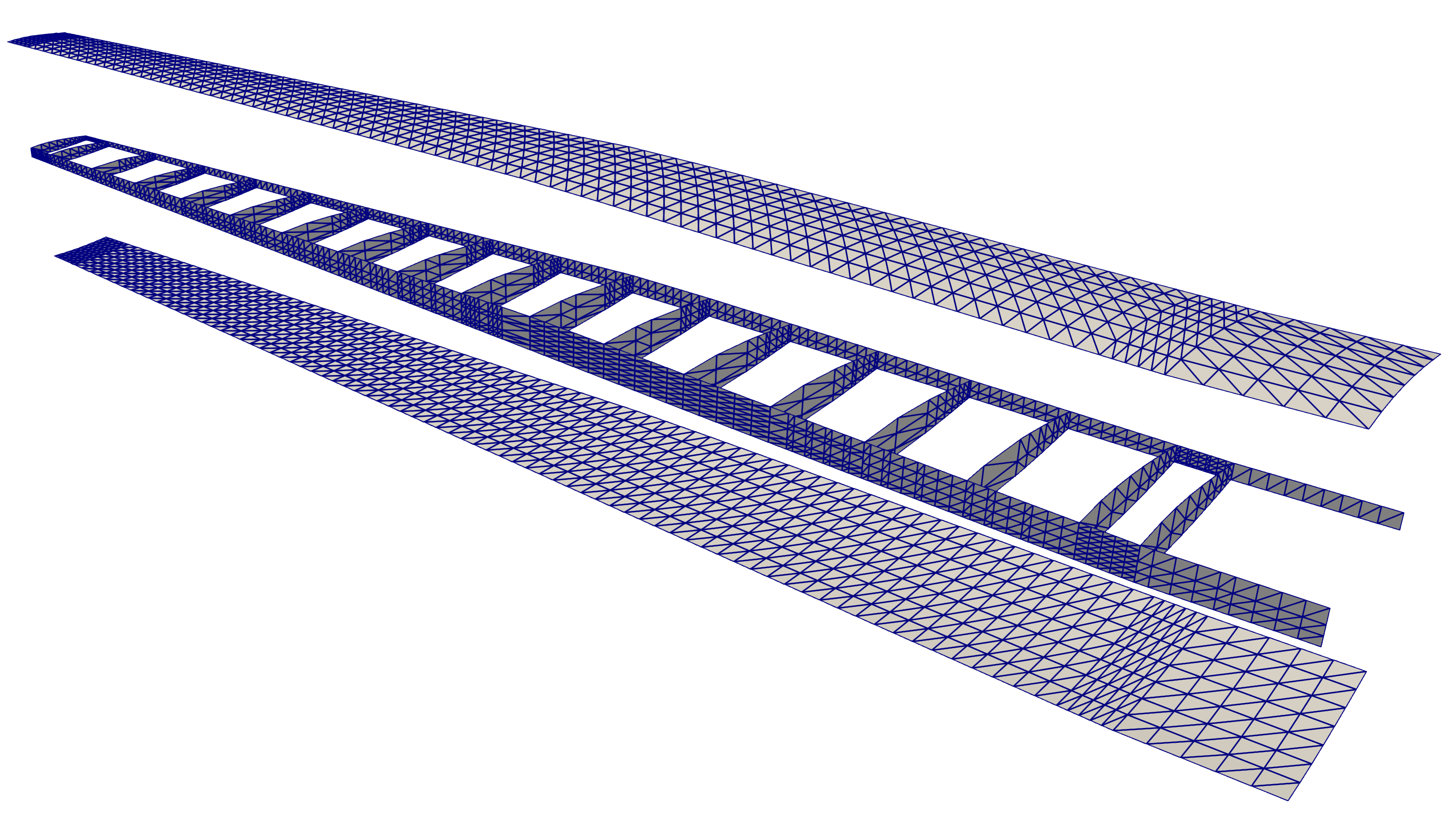}
    \caption{CAD geometry of the PEGASUS wing which is composed of 90 NURBS patches with 280 intersections, totaling 19572 DoFs.}
    \label{fig:pegasus-wing-geom}
\end{figure}

The PEGASUS wing CAD model is created using a customized geometry engine, an exploded view of the wing with IGA discretization is shown in Figure \ref{fig:pegasus-wing-geom}. The CAD model consists of 90 NURBS patches, two outer skins (one lower skin and one upper skin) and two spars (one front spar and one rear spar) connecting two adjacent ribs. The NURBS surfaces in the PEGASUS wing are represented using cubic basis functions with maximal continuity, resulting in 19572 DoFs in total. Additionally, 280 patch intersections are formed in the wing structure. 

\begin{figure}[!htb]
    \centering
    \begin{subfigure}[t]{0.49\textwidth}
        \centering
        \includegraphics[width=\textwidth]{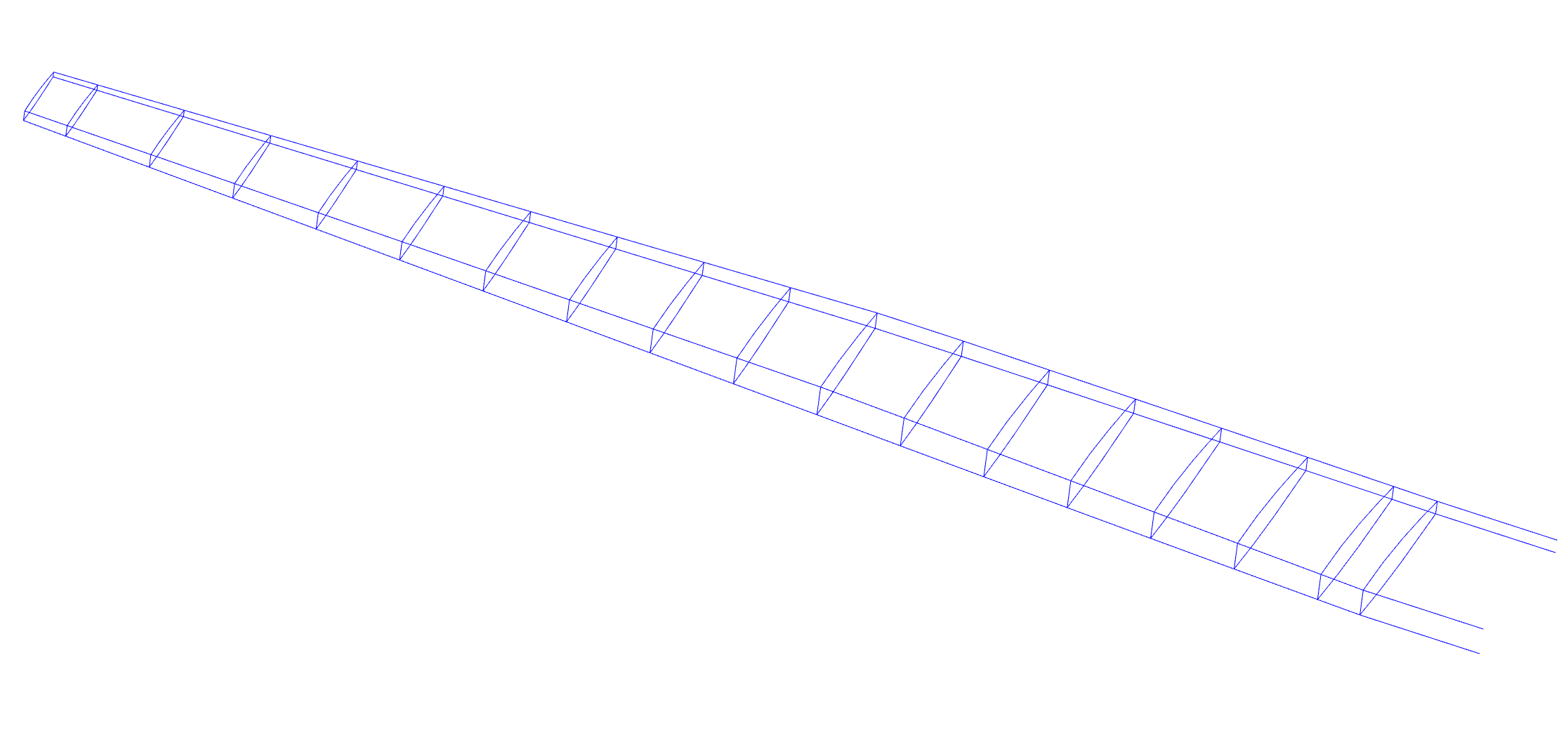}
        \caption{Illustration of surface--surface intersections in the PEGASUS wing geometry.}
        \label{subfig:pegasus-wing-intersections}
    \end{subfigure}
    \hfill
    \begin{subfigure}[t]{0.49\textwidth}
        \centering
        \includegraphics[width=\textwidth]{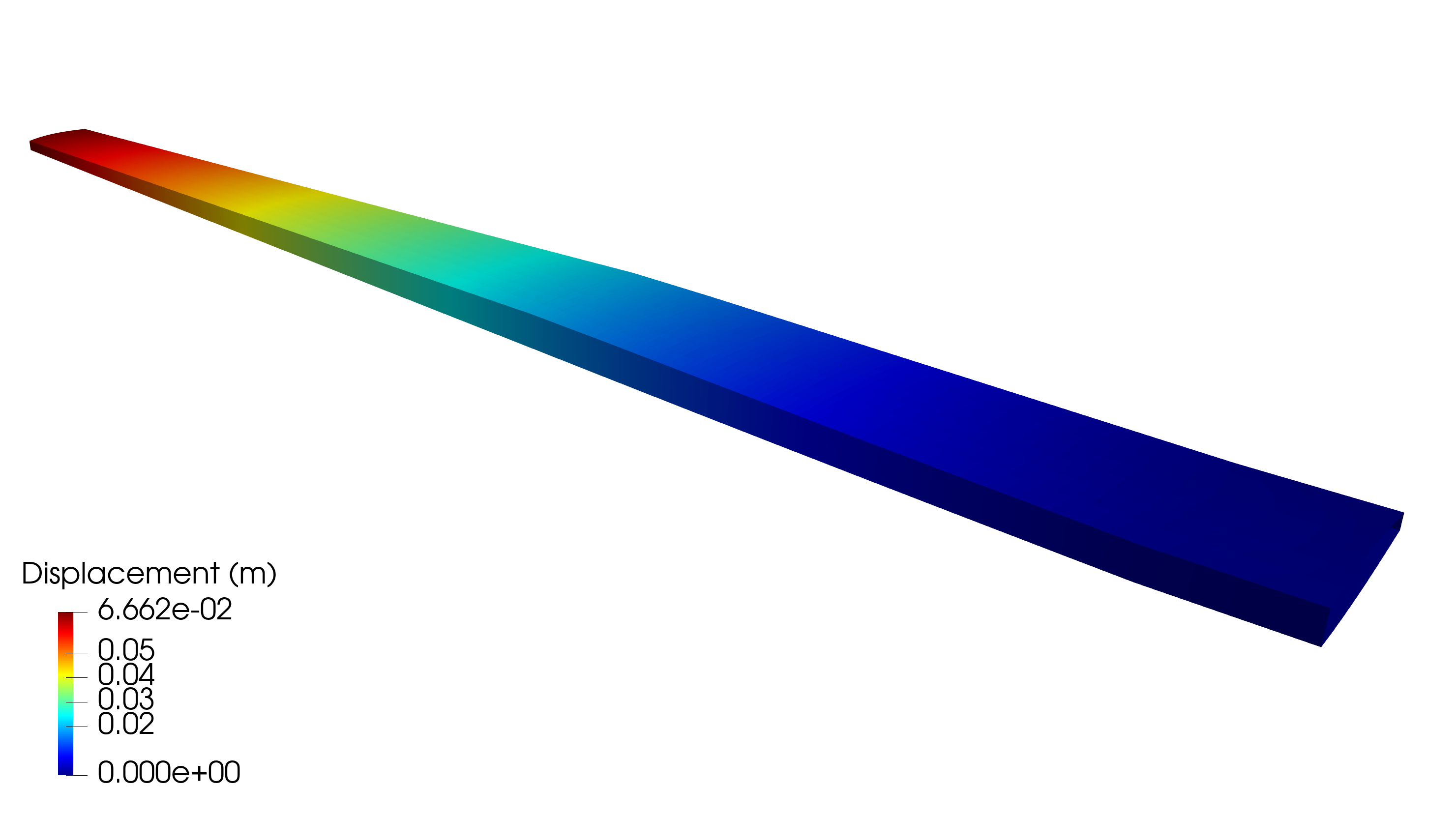}
        \caption{Visualization of displacement magnitude solved by PENGoLINS in baseline design.}
        \label{subfig:pegasus-wing-disp-mag-pengolins}
    \end{subfigure}
    \begin{subfigure}[t]{0.49\textwidth}
        \centering
        \includegraphics[width=0.9\textwidth]{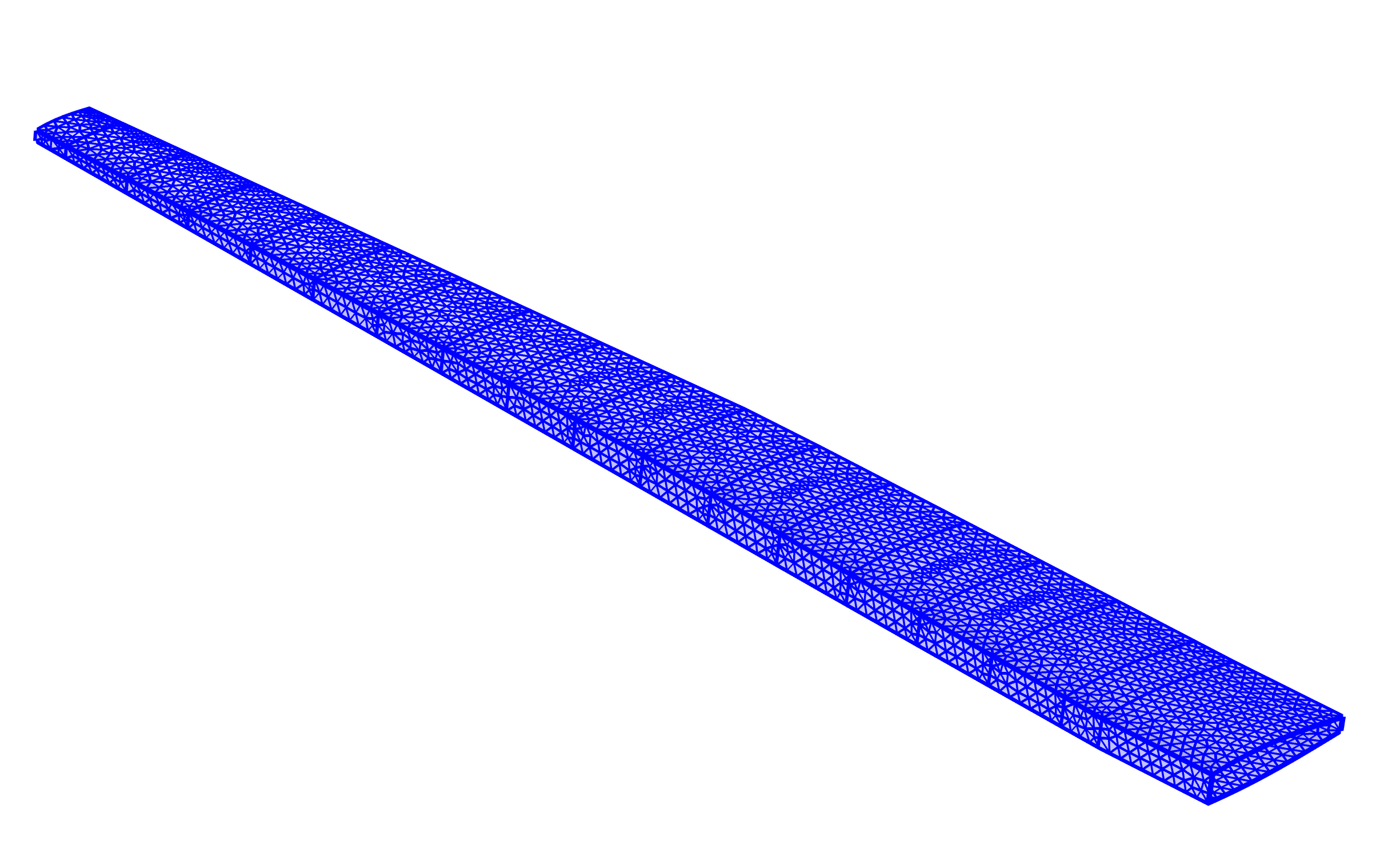}
        \caption{Finite element mesh of PEGASUS wing generated in COMSOL.}
        \label{subfig:pegasus-wing-mesh-comsol}
    \end{subfigure}
    \hfill
    \begin{subfigure}[t]{0.49\textwidth}
        \centering
        \includegraphics[width=\textwidth]{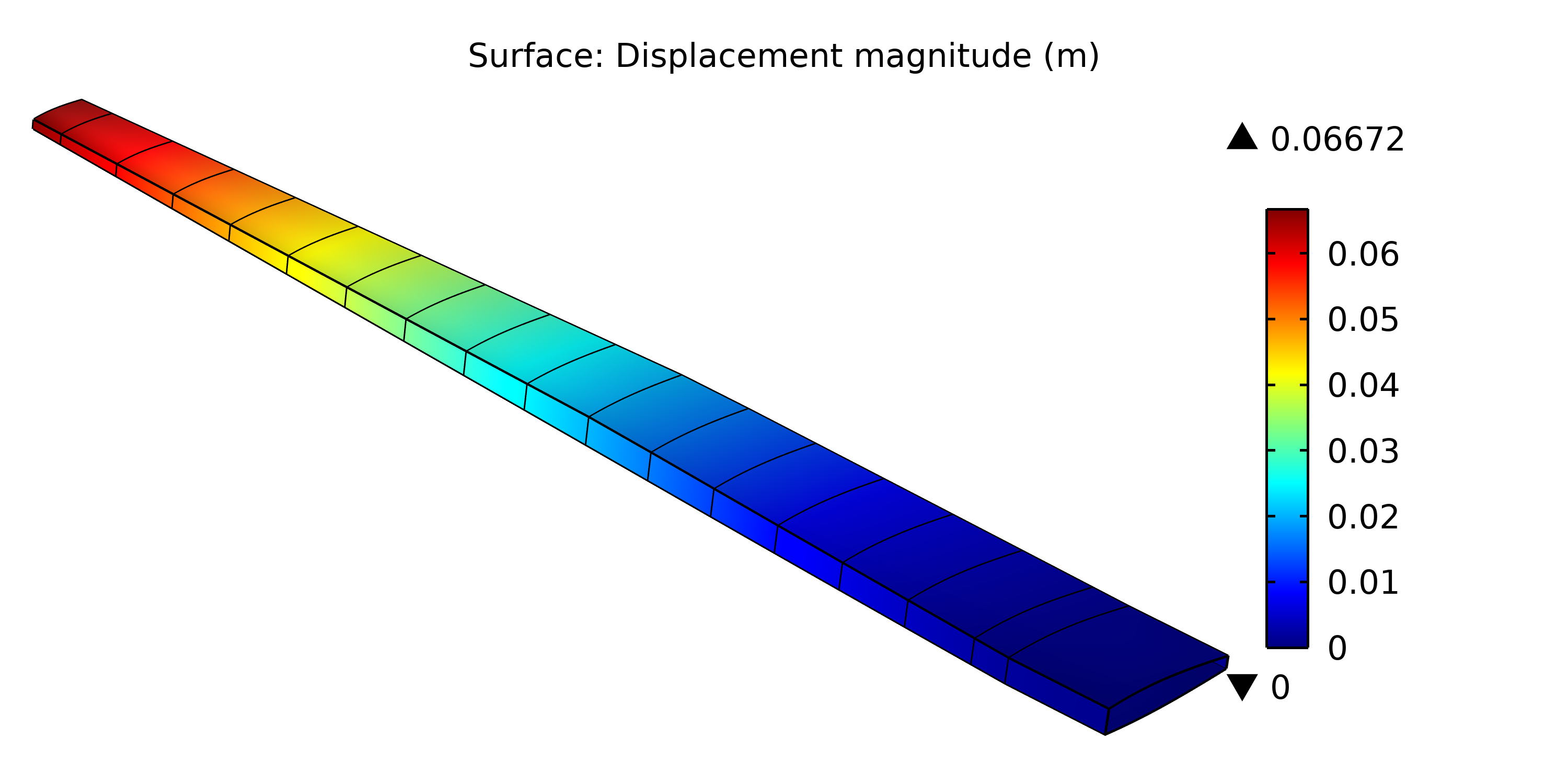}
        \caption{Analysis result obtained from COMSOL using Reissner--Mindlin shell theory.}
        \label{subfig:pegasus-wing-disp-mag-comsol}
    \end{subfigure}
    \caption{Structural analysis of the PEGASUS wing using PENGoLINS, and the resulting displacement magnitude is compared with the corresponding output from COMOSL.}
    \label{fig:pegasus-wing-analysis}
\end{figure}

The PEGASUS wing is made of material with Young's modulus 69 GPa and Poisson's ratio 0.35, and the wing span is 12.22 m. At the wing root, the chord is 1.52 m and the airfoil thickness is 0.37 m. A uniform initial thickness is 5 mm for all patches. Considering an aircraft take-off weight of 9000 kg, a distributed upward pressure with a magnitude of 132.5 $\mathrm{N/m^2}$ is determined by dividing half of the take-off weight by the surface area of the wing. Clamped boundary conditions are imposed at the wing root. Importing the PEGASUS wing geometry in IGES format to PENGoLINS, we perform structural analysis directly without finite element mesh generation. Given an analysis-suitable CAD file, the only required geometry preprocessing is to approximate surface--surface intersections, which is a much easier effort than finite element mesh generation and is automated in PENGoLINS using the functionality provided by pythonOCC. Figure \ref{subfig:pegasus-wing-intersections} shows all the intersections presented in the PEGASUS wing, while the displacements computed by PENGoLINS are visualized in Figure \ref{subfig:pegasus-wing-disp-mag-pengolins}.

To validate the calculated displacements, we conduct FE analysis for the PEGASUS wing using COMSOL. Figure \ref{subfig:pegasus-wing-mesh-comsol} displays an extensively refined finite element mesh for the COMSOL FE analysis. Displacements solved in COMSOL, utilizing the Reissner--Mindlin shell theory, are depicted in Figure \ref{subfig:pegasus-wing-disp-mag-comsol}. In this analysis, quadratic triangular elements with 118644 DoFs are used. Figure \ref{fig:pegasus-wing-analysis} indicates that the displacements obtained from PENGoLINS closely match the results from COMSOL. The maximum displacement magnitude in PENGoLINS is 0.06662 m, which has a relative difference of 0.15\% compared to the corresponding value of 0.06672 m in COMSOL. This aligns well with the findings of a numerical experiment presented in \cite[Section 5.3]{Zhao2022}, where a relative difference of 0.76\% is observed for the converged vertical displacement at the wingtip of an eVTOL wing between PENGoLINS and an open-source Reissner--Mindlin shell solver \cite{shell_analysis_fenicsx_code} using formulations from \cite{Campello2003}. The simulation results for the PEGASUS wing indicate that PENGoLINS provides good accuracy for complex shell structures, which is crucial for the subsequent design optimization.

\subsubsection{Thickness optimization of the PEGASUS wing}\label{subsubsec:pegasus-wing-thopt}

\begin{figure}[!htb]\centering
    \includegraphics[width=0.6\textwidth]{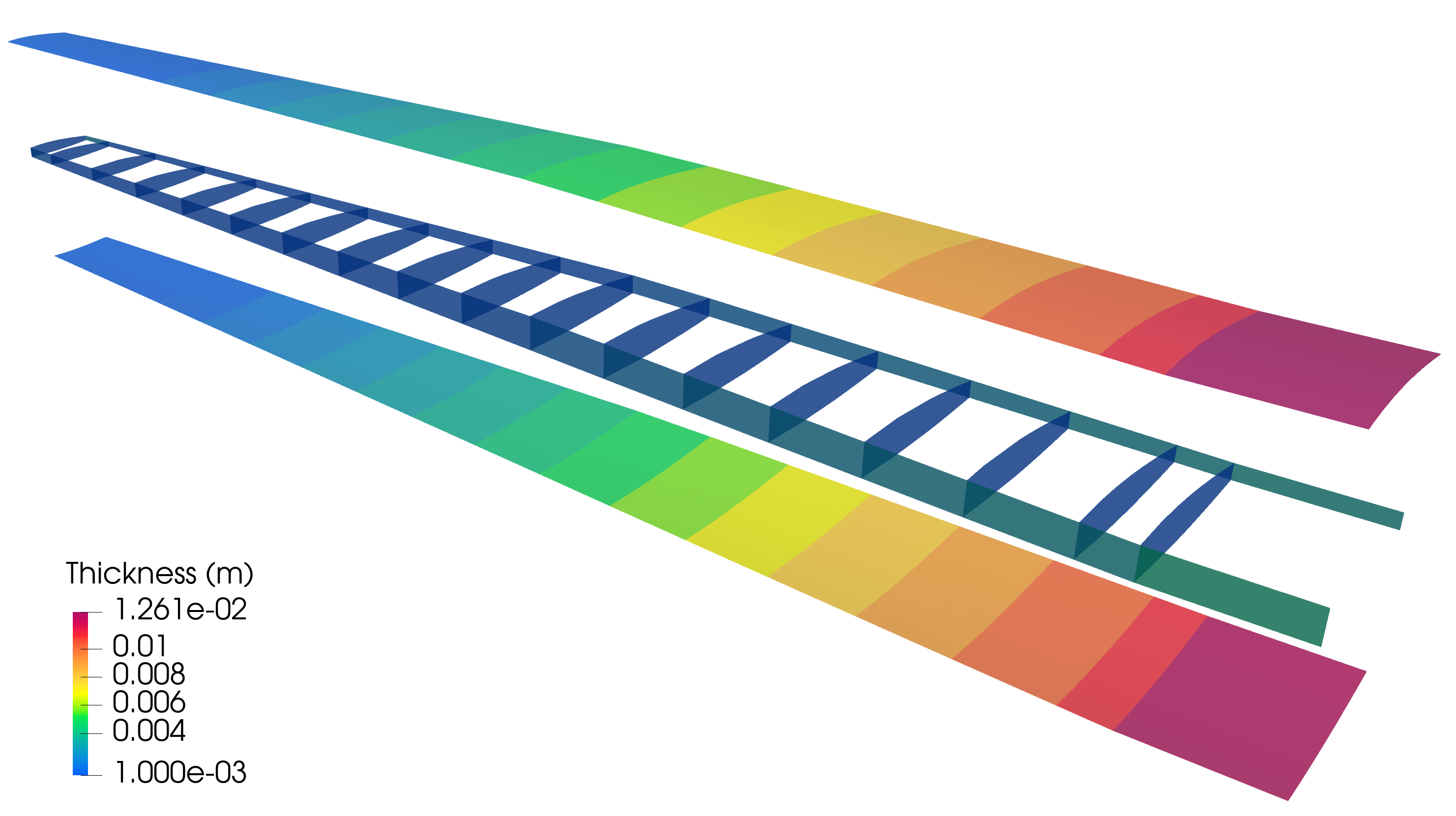}
    \caption{Optimization result of the PEGASUS wing with piecewise constant thickness.}
    \label{fig:pegasus-wing-const-thopt}
\end{figure}

Similar to the thickness optimization of the clamped plate discussed in Section \ref{subsec:plate-th-opt}, we apply the same methodology to the PEGASUS wing for piecewise constant thickness optimization. The same boundary conditions are employed throughout the optimization. In the piecewise constant thickness optimization case, a total of 90 design variables are included with lower and upper bounds of 1 mm and 100 mm, respectively. The initial thickness for all patches is taken as 5 mm. A constant volume constraint is employed in the optimization. The resulting shell thickness with minimum internal energy is depicted in Figure \ref{fig:pegasus-wing-const-thopt}. The shell patch with the maximum thickness is observed at the wing root in the outer skins, while the thickness decreases consistently along the span direction for both the outer skins and spars, following a pattern similar to that of the clamped plate. The thicknesses of the internal ribs and the wingtip are close to the lower bound since the bending moments are mainly carried by the lower and upper skins given the distributed upward load. Therefore, the majority of material is redistributed towards the clamped root of the outer skins. The optimized design in Figure \ref{fig:pegasus-wing-const-thopt} gives an internal energy 38.17\% less than that of the baseline configuration.

\begin{figure}[!htb]
    \centering
    \begin{subfigure}[t]{0.49\textwidth}
        \centering
        \includegraphics[width=\textwidth]{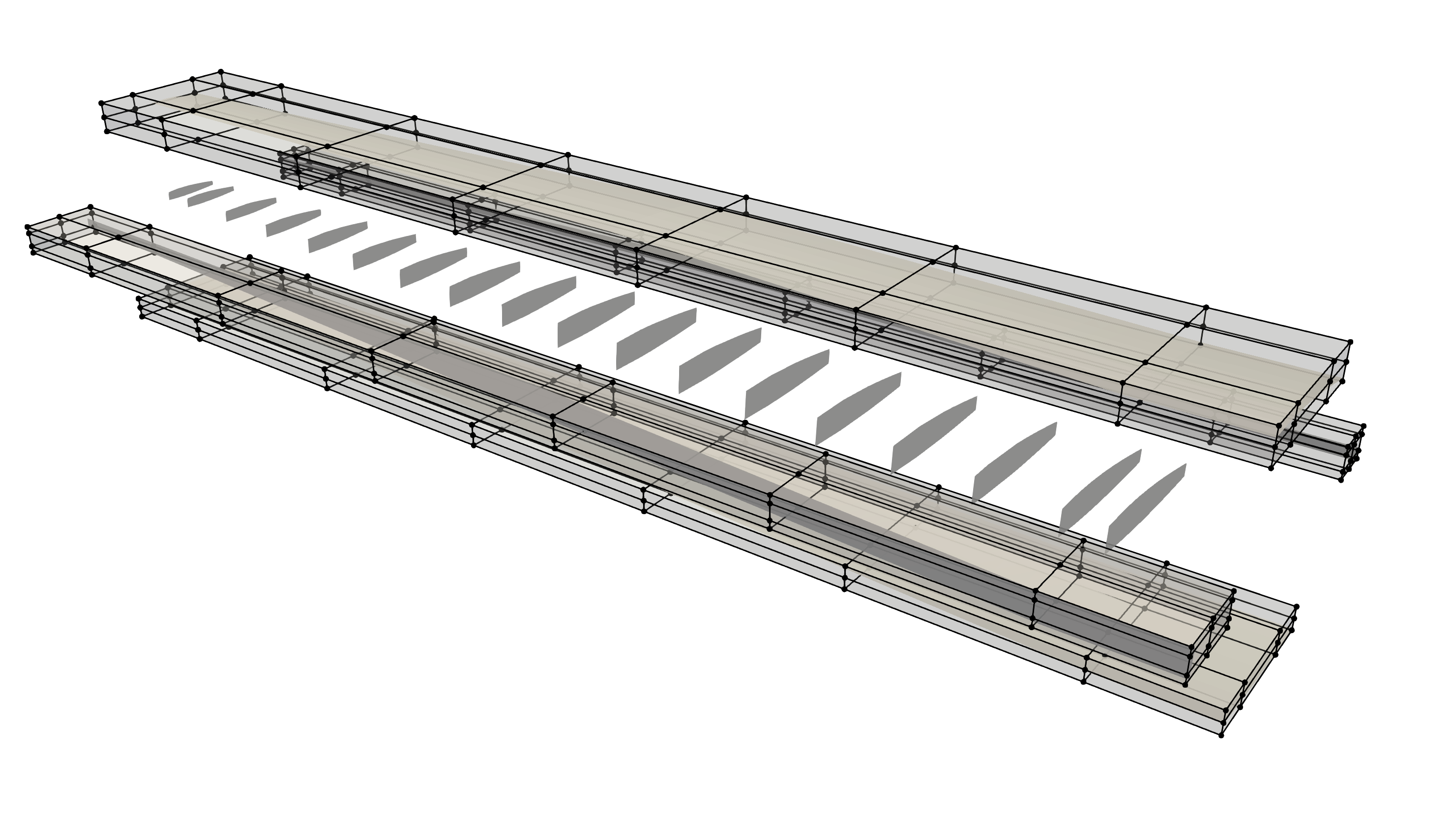}
        \caption{}
        \label{subfig:pegasus-wing-var-thopt-ffd}
    \end{subfigure}
    \hfill
    \begin{subfigure}[t]{0.49\textwidth}
        \centering
        \includegraphics[width=\textwidth]{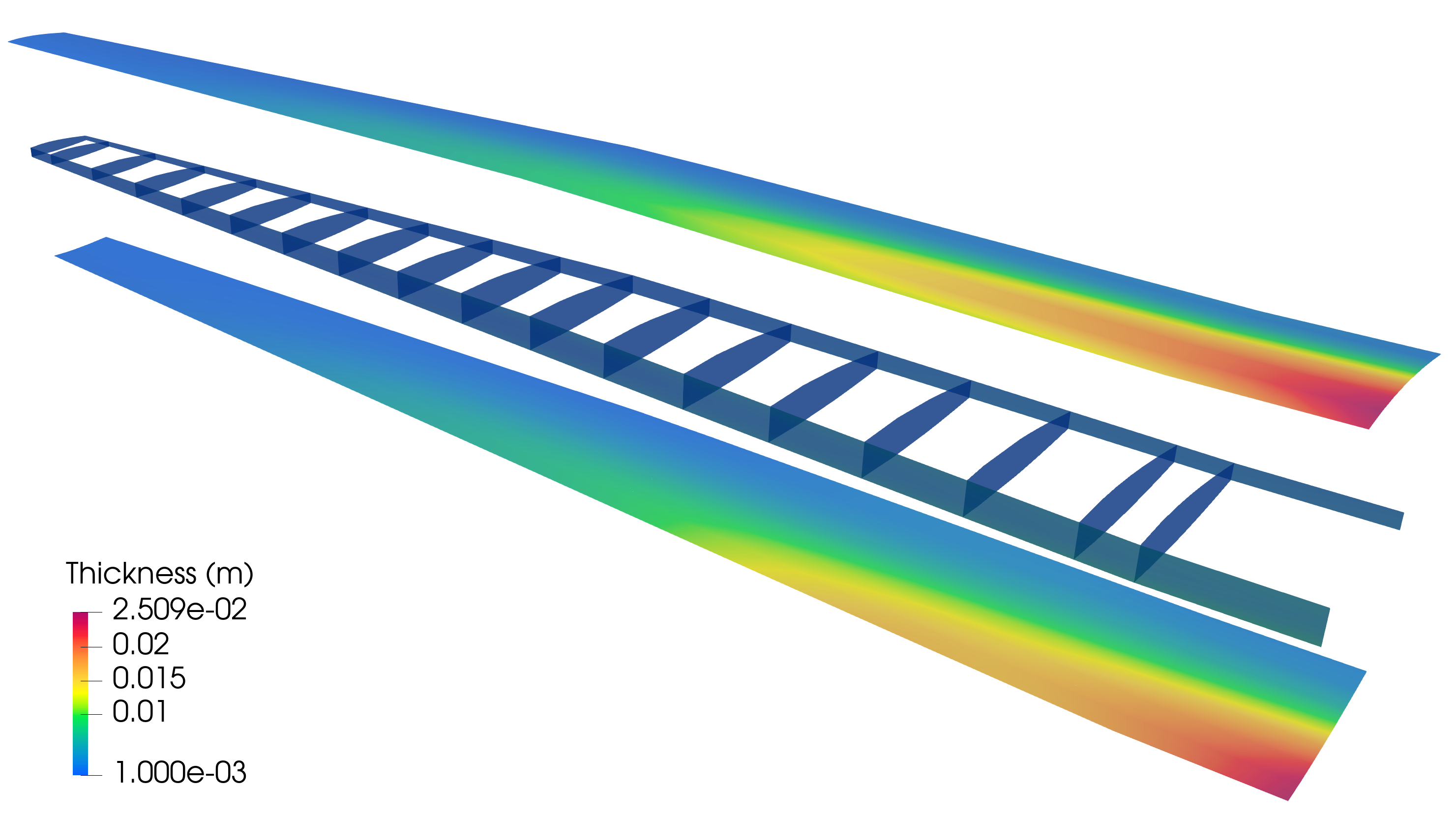}
        \caption{}
        \label{subfig:pegasus-wing-var-thopt}
    \end{subfigure}
    \caption{(a) Configuration of the combined thickness optimization. Each group of outer skins and spars is placed in one FFD block, and the remaining internal ribs have a piecewise constant thickness. (b) Optimal thickness distribution of PEGASUS wing.}
    \label{fig:pegasus-wing-thopt}
\end{figure}

To achieve an improved design, we consider variable thickness in the outer skins and spars of the PEGASUS wing, while keeping the internal ribs with piecewise constant thickness. The configuration of the FFD blocks for variable thickness optimization is illustrated in Figure \ref{subfig:pegasus-wing-var-thopt-ffd}. Four FFD blocks with quadratic bases are created to allow for variation in thickness within a single spline patch while ensuring continuity at patch intersections, with a total of 402 design variables used for this problem. By minimizing the internal energy again, the optimal thickness distribution is obtained as shown in Figure \ref{subfig:pegasus-wing-var-thopt}. Both applications in this section use the SNOPT optimizer with a tolerance of $10^{-4}$. In the optimized design, the outer skins of the PEGASUS wing at the clamped edge still have the largest thickness, where the thickness distribution is smooth at the surface--surface intersections within each FFD block. The decreasing thickness trend in outer skins and spars remains consistent with the optimal piecewise constant thickness case. Comparing the combined optimization strategy with the piecewise constant method, the maximum thickness in the former is higher while maintaining the same volume. Additionally, the internal energy is reduced by 44.71\% compared to the baseline design. These observations indicate that the combined thickness optimization method demonstrates a more efficient utilization of material than the piecewise constant method. 

\subsection{Simultaneous optimization for eVTOL wing}\label{subsec:evtol-wing-shthopt}
\begin{figure}[!htb]
    \centering
    \begin{subfigure}[t]{0.49\textwidth}
        \centering
        \includegraphics[width=\textwidth]{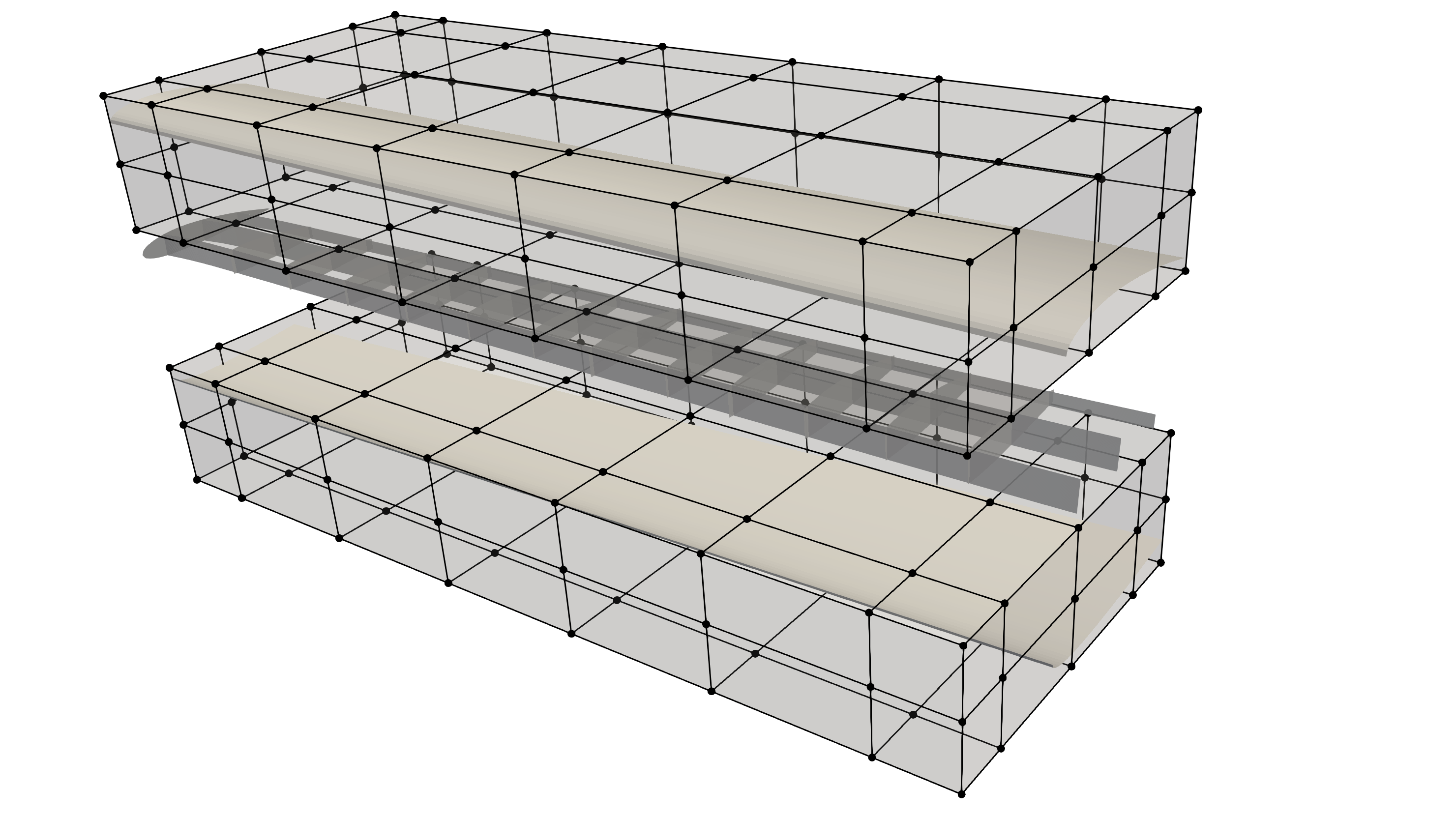}
        \caption{}
        \label{subfig:evtol_wing_thopt_ffd}
    \end{subfigure}
    \hfill
    \begin{subfigure}[t]{0.49\textwidth}
        \centering
        \includegraphics[width=\textwidth]{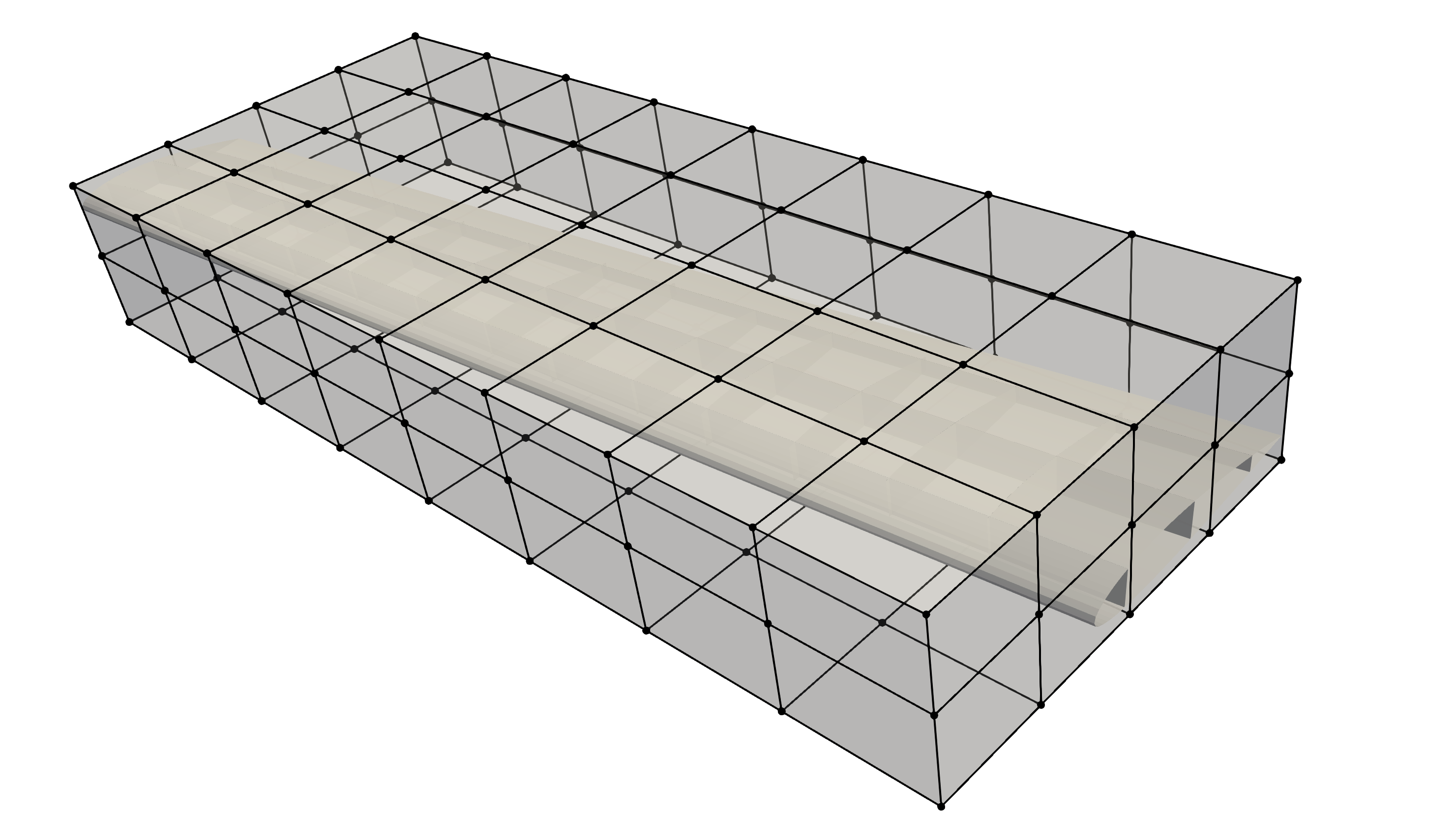}
        \caption{}
        \label{subfig:evtol_wing_shopt_ffd}
    \end{subfigure}
    \caption{(a) FFD blocks for eVTOL wing thickness optimization, where the lower and upper skins have a variable thickness. Wingtips and internal ribs and spars have a piecewise constant thickness. (b) FFD block for shape optimization.}
    \label{fig:evtol_wing_shthopt_ffd}
\end{figure}

With the continuous advancements in aviation battery technology \cite{liu2022-evtol}, eVTOL aircraft have emerged as a promising solution for cost-effective urban mobility \cite{Polaczyk2019}. In this section, we use a more advanced eVTOL wing to demonstrate the capabilities of the FFD-based optimization approach, where both the thickness and shape control points are considered as design variables simultaneously. By incorporating both thickness and shape coordinates in the design of shell structures, we can utilize the material more efficiently than considering thickness optimization only. The CAD geometry of the eVTOL wing, material parameters, and the corresponding structural analysis can be found in \cite[Section 5]{Zhao2022}. For the optimization problem, we implement the same clamped boundary conditions and distributed upward pressure as those in the previous application. The magnitude of pressure is set to 120 $\mathrm{N}/\mathrm{m}^2$, approximated by dividing the take-off weight by the surface area of the eVTOL wing. The baseline configuration of the eVTOL wing, which is composed of 27 cubic NURBS patches with 87 intersections, is illustrated in Figure \ref{fig:evtol_wing_shthopt_ffd}. We note that we perform the shape optimization for the eVTOL wing without including an aerodynamic model (which would be needed for a well-posed wing design problem) purely to provide a demonstration of the FFD-based method for complex aerospace structures.

To achieve a meaningful design for the eVTOL wing, we use two quadratic B-spline FFD blocks for thickness optimization. This configuration allows for varying thicknesses in the lower and upper skins of the eVTOL, while using piecewise constant thicknesses for the internal stiffeners and wingtips. The arrangement of the thickness FFD blocks is illustrated in Figure \ref{subfig:evtol_wing_thopt_ffd}. Furthermore, a cubic B-spline FFD block is created for shape optimization, as depicted in Figure \ref{subfig:evtol_wing_shopt_ffd}. Only the vertical coordinates of control points for the shape FFD block, denoted as $\mathbf{P}_{\text{FFD}3}$, are updated. In total, there are 642 design variables involved in this optimization process, and a constant volume constraint is applied as well. Regarding the thickness design variables, the lower and upper limits are selected as 1 mm and 50 mm, respectively. All shell components have initial thicknesses of 3 mm.

\begin{figure}[!htb]\centering
    \includegraphics[width=0.9\textwidth]{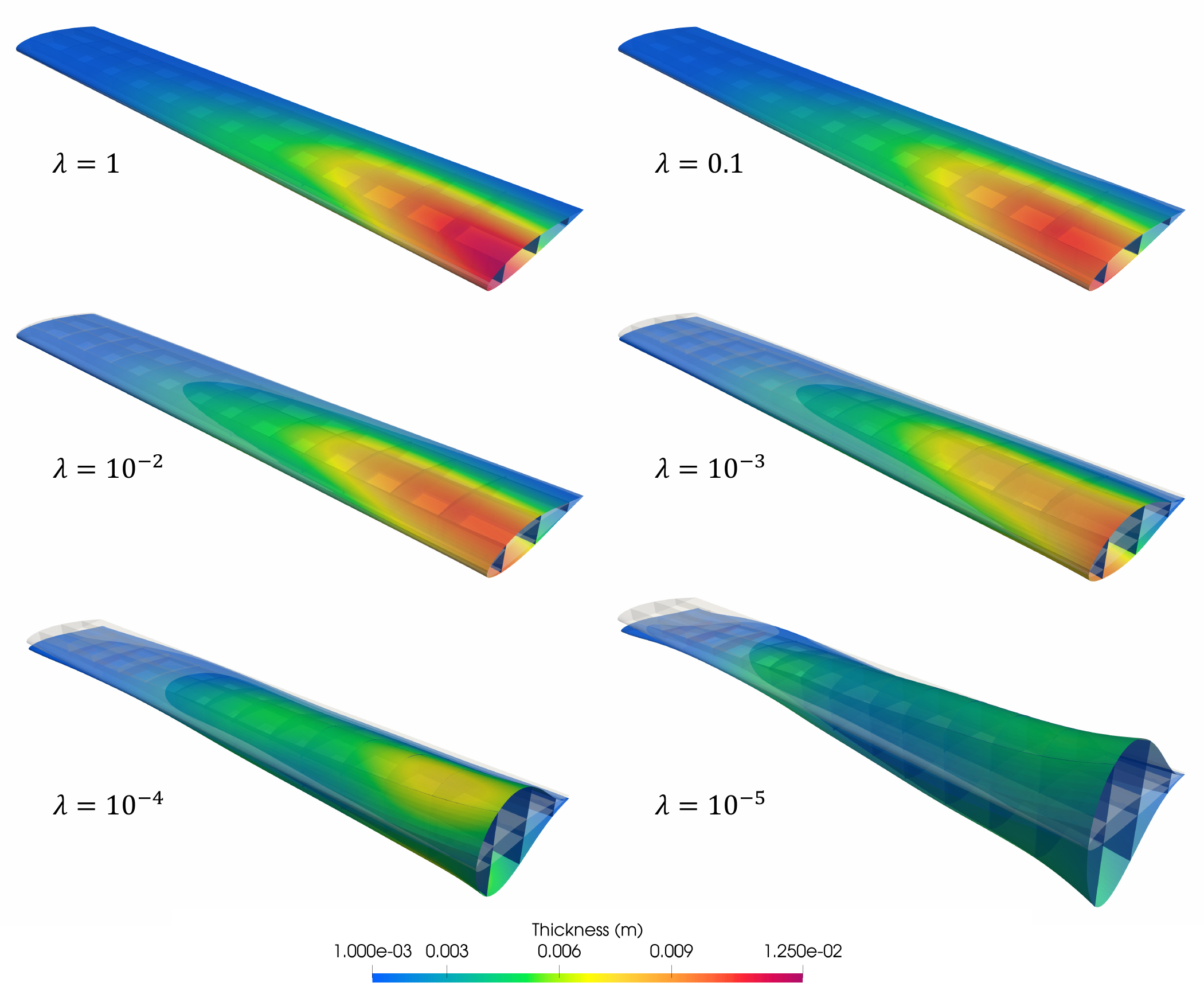}
    \caption{Optimized solutions of the eVTOL wing with varying regularization coefficient.}
    \label{fig:evtol-wing-shthopt-regu-test}
\end{figure}

One challenge encountered during shape optimization of complex shell structures, such as eVTOL wings, is the potential occurrence of oscillatory or highly distorted geometries in the updated designs. These distorted shapes can lead to poor element quality, resulting in spurious energy and affecting the accuracy of the analysis. To mitigate oscillation or radical change of shell components, an additional term is introduced in the objective function to regularize the gradient of the shape. The objective is formulated as
\begin{align}
    f^\text{obj} = \sum_{\text{I}=1}^{m} \left(W_{\text{int}}^{\text{I}} + \lambda \frac{E (\mathring{t}^\text{I})^3}{12 \mathring{h}^\text{I}_A (1-\nu^2)} \int_{\mathcal{S}^{\text{I}}} \Vert \nabla \mathbf{P}^\text{I}_{3} - \nabla \mathring{\mathbf{P}}^\text{I}_{3} \Vert^2 \ d\mathcal{S}\right) \text{ ,} \label{eq:wint-shape-regu}
\end{align}
where $\text{I}$ is the index of shell patches, $\lambda$ is a dimensionless regularization coefficient, $h^I_A$ is the element area of shell patch $\text{I}$ in the physical space, $\mathbf{P}^\text{I}_{3}$ is the vertical component of the shape variable for shell patch $\text{I}$, and $\mathcal{S}^{\text{I}}$ is the midsurface of shell patch $\text{I}$. $\mathring{(\cdot)}$ denotes quantities in the baseline configuration. The regularization term in \eqref{eq:wint-shape-regu} serves as an additional artificial bending energy associated with the curvature of the eVTOL wing. Therefore, the shape oscillation can be eliminated by adjusting the regularization coefficient $\lambda$. The SNOPT optimizer is employed with a tolerance of $10^{-3}$. Figure \ref{fig:evtol-wing-shthopt-regu-test} demonstrates the optimized eVTOL wing designs achieved with varying $\lambda$ values, providing insights into the influence of regularization on the final shape and thickness distribution. It can be seen that the patch intersections of the optimized shapes in Figure \ref{fig:evtol-wing-shthopt-regu-test} are still connected using the FFD-based approach. The reductions of internal energy for different $\lambda$ are listed in Table \ref{tab:wint-reduction-lambda}.

\begin{figure}[!htb]\centering
    \includegraphics[width=0.9\textwidth]{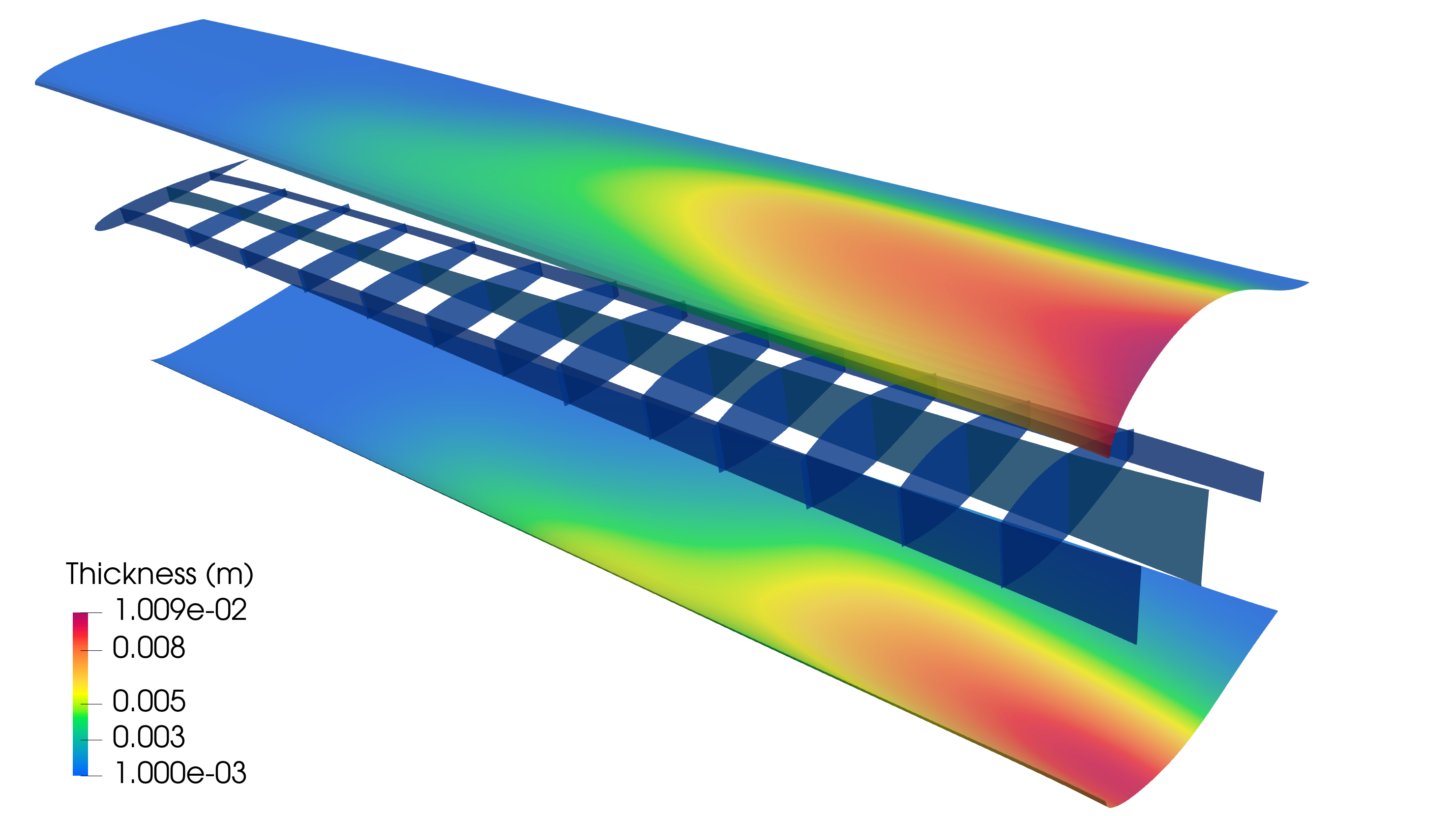}
    \caption{Exploded view of optimal design of eVTOL wing with regularization coefficient $\lambda=10^{-3}$, which results in 83.23\% reduction of internal energy compared to the baseline design.}
    \label{fig:evtol-wing-shthopt-regu1e-3-explosive}
\end{figure}

\begin{table}[!htbp]
    \centering
    \setlength{\tabcolsep}{8pt} 
    \begin{tabular}{c c c c c c c}
        \toprule
        $\lambda$ &  1 & 0.1 & $10^{-2}$ & $10^{-3}$ & $10^{-4}$ & $10^{-5}$\\
        \midrule
        \makecell{Internal energy \\reduction (\%)} & 49.07 & 51.58 & 66.45 & 83.23 & 92.82 & 95.09 \\
        \bottomrule
    \end{tabular}
    \caption{Reduction of internal energy of the eVTOL wing after simultaneous optimization with varying regularization coefficients.}
    \label{tab:wint-reduction-lambda}
\end{table}
For this eVTOL wing optimization problem, different regularization coefficients yield distinct outcomes. A regularization coefficient of $1$ and $0.1$ provide optimal designs dominated by thickness update, resembling the patterns observed in the combined thickness optimization of the PEGASUS wing, refer to Figure \ref{subfig:pegasus-wing-var-thopt}. The optimized shapes in these two cases are similar to the baseline configuration since the artificial bending energies are the major contributor to the objective function due to the large regularization coefficients. Even slight variations in the shape variables result in substantial increases in the objective function. Conversely, small regularization coefficients, i.e., $10^{-4}$ and $10^{-5}$, lead to considerable changes in the shape of the eVTOL wing and reduction of internal energy, amounting to $92.82\%$ and $95.09\%$, respectively. However, these cases exhibit noticeable oscillations in the wingtip area, which can lead to ill-conditioned systems. On the other hand, employing regularization coefficients with values of $10^{-2}$ and $10^{-3}$ yields balanced solutions, where both the thickness redistribution and shape updates contribute to the optimal design. The material moves towards the clamped area, accompanied by a widening of the cross-section to provide increased wing support. An exploded view of the optimal design with $\lambda=10^{-3}$ is shown in Figure \ref{fig:evtol-wing-shthopt-regu1e-3-explosive}. This optimal design achieves an internal energy reduction of 83.23\% compared to the baseline configuration. 

\section{Conclusion}\label{sec:conclusion}
In this paper, we have introduced an FFD-based shape and thickness optimization approach for shell structures composed of separately-parametrized NURBS surfaces. The integration of this method with the Lagrange extraction technique enables IGA with existing FE toolkits and provides a connection between the FFD block and non-matching shell patches. By employing the FFD block approach, the updated shell geometry and thickness remain properly connected at patch intersections throughout the optimization process. This feature prevents undesired shape discontinuities in shell structures. 

We made use of the penalty-based coupling formulation for Kirchhoff--Love shells to perform IGA in the optimization. The automation of analytical derivative computations is achieved through code generation in FEniCS, enabling gradient-based multidisciplinary design optimization. This automation streamlines the optimization process and allows for efficient exploration of design spaces. The unified NURBS representation shared by both the design geometry and analysis model enhances accuracy per DoF in the analysis and precise design updates. Moreover, the proposed framework circumvents FE mesh generation and streamlines design-analysis-optimization workflow for complex shell structures. Consequently, the automated workflow is expected to accelerate the conceptual design of novel eVTOL aircraft with minimal manual effort.

A suite of benchmark problems is adopted to verify the effectiveness of the FFD-based optimization approach proposed in this paper. Both the shape optimization and thickness optimization results agree well with analytical solutions or other established references. Furthermore, we have applied the framework to two different aircraft wings. This demonstration highlights the potential of the proposed method in exploring complex design spaces and obtaining superior designs for innovative aircraft structures. Future works can extend the shape optimization for aircraft wings with more realistic loads by coupling the structural solver with an aerodynamic solver \cite{van2023solver}. Additionally, practical constraints such as von Mises stress can be incorporated into the optimization process to ensure that the obtained designs meet the requirements of real-world applications. To promote code transparency and potential contributions to the shell optimization community, we make the Python library GOLDFISH openly accessible in a GitHub repository \cite{goldfish-code}.

\section*{Acknowledgements}\label{sec:Ackonwledgements}
This work is supported by National Aeronautics and Space Administration grant number 80NSSC21M0070. We thank Michael A. Warner (UCSD) and Darshan Sarojini (UCSD) for providing the CAD geometry of the PEGASUS wing in Figure \ref{fig:pegasus-wing-geom}.






\bibliography{main}


\end{document}